\documentclass[a4paper,12pt,twoside,openright,
titlepage=false]{scrreprt}
\newcounter{app1}
\newcounter{app2}
\newcounter{app3}
\newcounter{app4}
\newcounter{app5}

\usepackage{graphicx}
\usepackage{scrpage2}
\usepackage{amssymb}
\usepackage{exscale}
\usepackage{nicefrac}
\usepackage[thmmarks]{ntheorem}  
\usepackage{amsmath}
\usepackage{amssymb}
\usepackage{ifthen}
\usepackage[hyperref,usenames,dvipsnames]{xcolor}
\usepackage{bold-extra}
\usepackage[pdftex,pdfborderstyle={/S/U/W 1},plainpages=false,citebordercolor={0 0 1},linkbordercolor={0 1 0}]{hyperref}
\usepackage{cleveref}
\newtheorem{lemma}{Lemma}[chapter]
\newtheorem{korollar}[lemma]{Corollary}
\newtheorem{theorem}[lemma]{Theorem}
\newtheorem{prop}[lemma]{Proposition}
\newtheorem{definition}[lemma]{Definition}
\newtheorem{assumption}[lemma]{Assumption}

\crefname{lemma}{Lemma}{Lemmas}
\crefname{korollar}{Corollary}{Corollaries}
\crefname{theorem}{Theorem}{Theorems}
\crefname{prop}{Proposition}{Propositions}
\crefname{definition}{Definition}{Definitions}
\crefname{assumption}{Assumption}{Assumptions}
\Crefname{lemma}{Lemma}{Lemmas}
\Crefname{korollar}{Corollary}{Corollaries}
\Crefname{theorem}{Theorem}{Theorems}
\Crefname{prop}{Proposition}{Propositions}
\Crefname{definition}{Definition}{Definitions}
\Crefname{assumption}{Assumption}{Assumptions}

\newcommand{\hei}{10}
\newcommand{\barG}{\begin{picture}(9,0) $G$ \put(-10.5,\hei){\linethickness{0.03em}
  \line(3,0){5.4}}\end{picture}}
\newcommand{\barGamma}
   {\begin{picture}(9,0) $\Gamma$\put(-10.5,\hei){\linethickness{0.03em}
     \line(3,0){5.2}}\end{picture}}
\newcommand{\barD}{\begin{picture}(10.5,0)\put(-0.6,\hei){\linethickness{0.03em} 
  \line(3,0){5.4}}$D$\end{picture}}
\newcommand{\barA}{\begin{picture}(9,0) $A$ \put(-8.9,\hei){\linethickness{0.03em}
  \line(3,0){4.5}}\end{picture}}
\newcommand{\barnabla}{\begin{picture}(10,0)\put(-3,\hei){\linethickness{0.03em} 
  \line(1,0){7}}$\nabla$\end{picture}}

\newcommand{\gt}{\mathcal T}
\newcommand{\operp}{$\begin{picture}(8.5,4)\put(4,2.8){\thinlines\circle{7}\put(-10.2,
-1){\line(1,0){6.5}}\put(-7.1,-1){\line(0,1){4.3}}}\end{picture}$}
\newcommand{\mue}[1]{\mu\begin{picture}(4,5)\put(-1,-3){\scriptsize{$#1$}}\end{picture}}
\newcommand{\chif}[1]{\chi\begin{picture}(4,5)\put(-1,-5){\scriptsize{$#1$}}
\end{picture}}
\newcommand{\chinullf}[1]{\chi\begin{picture}(4,5)\put(-1,-5){\scriptsize{$#1$}}
\end{picture}\begin{picture}(0,0)\put(-8.5,6){\tiny$\circ$}\end{picture}}
\newcommand{\nf}[2]{^{#1}/_{#2} }

\newcommand{\ric}{\mbox{\rm ric}}
\newcommand{\inj}{\mbox{\rm{inj}}}
\newcommand{\dist}{\mbox{dist}}
\newcommand{\grad}{\mbox{grad} }

\newcommand{\prn}[1]{\partial_{\R^n}^{#1} }
\newcommand{\pt}{\partial_t}
\newcommand{\norm}[4]
  {\ifthenelse
    {\equal{#4}{}}
    {\lVert#1\rVert_{\stackrel{\scriptstyle~}{\scriptstyle#3}}^{\stackrel{\scriptstyle#2}
                                                              {\scriptstyle~}}}
    {\ifthenelse
      {\equal{#4}{1}}
            {\|#1\|_{\stackrel{\scriptstyle~}{\scriptstyle#3,[\gamma=1]}}
               ^{\stackrel{\scriptstyle#2}{\scriptstyle~}}}
            {\ifthenelse{\equal{#4}{0}}{
                  \|#1\|_{\stackrel{\scriptstyle~}{\scriptstyle#3,[\gamma>#4]}}
                     ^{\stackrel{\scriptstyle#2}{\scriptstyle~}}}
	          {\|#1\|_{\stackrel{\scriptstyle~}{\scriptstyle#3,[\gamma\geq#4]}}
	             ^{\stackrel{\scriptstyle#2}{\scriptstyle~}}}
	    } 
    }
  }
\newcommand{\Norm}[1]{\lVert#1\rVert}
\newcommand{\bb}{\bf Proof:~\rm}
\newcommand{\eb}{$\hfill{}_{_\blacksquare}$\\}
\newcommand{\bs}{\begin{picture}(5,5)
           \put(449,20){$~_{_\blacksquare}$}\end{picture}}
\newcommand{\er}{$\hfill_{\Diamond}\\$}
\newcommand{\is}{\int\limits_{\Sigma}}
\newcommand{\iz}{\!\!\int\limits_{[\g>0]}\!\!}

\newcommand{\cf}{\begin{picture}(10.5,0)\put(0.9,0.7){$\scriptstyle\circ f\,\,$}
  \end{picture}}
\newcommand{\pis}{\begin{picture}(10.5,0)\put(0.9,0.7){$\scriptstyle
   \circ\pi_\Sigma\,\,$}\end{picture}}
\newcommand{\pii}{\begin{picture}(10.5,0)\put(0.9,0.7){$\scriptstyle
   \circ\pi_I\,\,$}\end{picture}}
\newcommand{\pa}{\partial_\alpha}
\newcommand{\pb}{\partial_\beta}
\newcommand{\pg}{\partial_\gamma}
\newcommand{\pd}{\partial_\delta}
\renewcommand{\c}{\!\cdot\!}
\newcommand{\sigsig}{\Sigma_\sigma}
\newcommand{\sigsigrho}{\Sigma_{\sigma,\varrho}}
\newcommand{\sigrho}{\Sigma_{\varrho}}
\newcommand{\betrag}[2]{\lvert #1\rvert^{#2}}
\newcommand{\iss}{\!\!\int\limits_{\sigsig}}
\newcommand{\issr}{\!\!\int\limits_{\sigsigrho}}
\newcommand{\isr}{\!\!\int\limits_{\Sigma_\varrho}}
\newcommand{\kreis}[1]{\put(-9,14){\circle{8}}\put(-11.1,11.6){\scriptsize{#1}}}
\newcommand{\nr}[6]{\put(#1,#2){\circle{#3}}\put(#4,#5){\scriptsize{#6}}}
\renewcommand{\a}{\alpha}
\renewcommand{\b}{\beta}
\newcommand{\g}{\gamma}
\renewcommand{\d}{\delta}
\newcommand{\e}{\varepsilon}

\renewcommand{\t}{\tau}
\renewcommand{\r}{\varrho}
\newcommand{\wg}{\widetilde\gamma}
\renewcommand{\sec}{{\,\rm sec\put(-0,-3){\tiny{(M,g)}}\;\;\;\;\;\,}\,}
\newcommand{\eps}{\varepsilon}

\newcommand{\W}{\mathcal W}

\newcommand{\bet}[1]{\lvert #1\rvert}
\newcommand{\R}{\mathbb R}
\newcommand{\N}{\mathbb N}
\newcommand{\T}{\mathcal T}

\newcommand{\mm}{d\mu}
\newcommand{\no}{\nonumber}

\newcommand{\go}{\bar g}

\newcommand{\goophi}[2]{\go_{\varphi_x}^{\,jk}}
\newcommand{\ka}{\big(}
\newcommand{\kz}{\big)}

\newcommand{\p}{\partial\hspace{0.05em}}
\newcommand{\s}{\!\ast\!}
\renewcommand{\*}{f^*\!}
\newcommand{\dipl}[1]{}
\newcommand{\w}{\widetilde}
\newcommand{\na}{\nabla}
\newcommand{\Q}[3]{\ifthenelse{\equal{#3}{x}}
                      {Q^{#1,#2}_{}}
                      {\ifthenelse{\equal{#3}{y}}
                          {Q^{#1,#2}_{}}
                          {\ifthenelse{\equal{#3}{R}}
                              {Q^{#1,#2}_{\!R\,\star R}}
                              {Q^{#1,#2}_{#3}}
                          }
                      }
                   }
\newcommand{\pe}[2]{P^{#1}_{#2}(A)}
\newcommand{\eck}[1]{\langle #1 \rangle}
\newcommand{\ds}{\mm}
\newcommand{\info}[1]{}
\renewcommand{\text}[1]{(#1) }
\newcommand{\notiz}[1]{}

\newcommand{\absch}[1]{\chi\begin{picture}(18,2)\put(0,-3.5){$\scriptstyle B_{#1}$} 
  \end{picture}}
\newcommand{\chiunten}[1]{\chi\begin{picture}(18,2)\put(0,-5){$\scriptstyle {#1}$} 
  \end{picture}}
\newcommand{\gm}[1]{\Gamma_{\!#1}}
\newcommand{\wgm}[1]{\widetilde\Gamma_{\!#1}}

\newcommand{\lam}[1]{\Lambda_{#1} }

\newcommand{\ricci}{\mbox{ricci}}

\newcommand{\tildeg}{\,\widetilde{\!g}\,}

\newcommand{\2}{\!\!}
\newcommand{\3}{\!\!\!}
\newcommand{\4}{\!\!\!\!}
\newcommand{\5}{\!\!\!\!\!}
\newcommand{\6}{\!\!\!\!\!\!}
\newcommand{\7}{\!\!\!\!\!\!\!}
\newcommand{\bg}{15}
\newcommand{\bgminusdrei}{12}
\newcommand{\bgpluseins}{16}
\newcommand{\RB}[1]{R^{#1}_{\!\bot}\!}

\newcommand{\kle}[1]{#1}
\newcommand{\gro}[1]{#1}

\renewcommand{\gg}{\Gamma}

\newcommand{\np}[1]{\bet{\na^{#1}\phi}}
\newcommand{\strich}[1]{\begin{picture}(9.3,1)\put(0,-1.8){\Big\lvert}
  \put(3,-8){\footnotesize#1}\end{picture}}
  
\newcommand{\strichklein}[1]{\begin{picture}(9.3,1)\put(0,-1.8){\big\lvert}
  \put(3,-5){\footnotesize{#1}}\end{picture}}
\newcommand{\abint}[1]{\!\!\int\limits_{[\gamma\geq#1]}}
\newcommand{\rg}{{}^{\varrho} g}
\newcommand{\bbet}[1]{\bet{#1}}

\clearscrheadfoot
\rofoot[\pagemark]{\pagemark}
\lefoot[\pagemark]{\pagemark}
\begin{document}
\newcounter{satz}
\newcounter{beh}
\newcounter{eqnarray}

\setlength{\topmargin}{0cm}
\title{\vspace{-3em}\rm\LARGE\scshape{Albert-Ludwigs-Universit\"at Freiburg\\
Fakult\"at f\"ur Mathematik und Physik}\\\vspace{3em}\rm
\scalebox{0.07}{\includegraphics{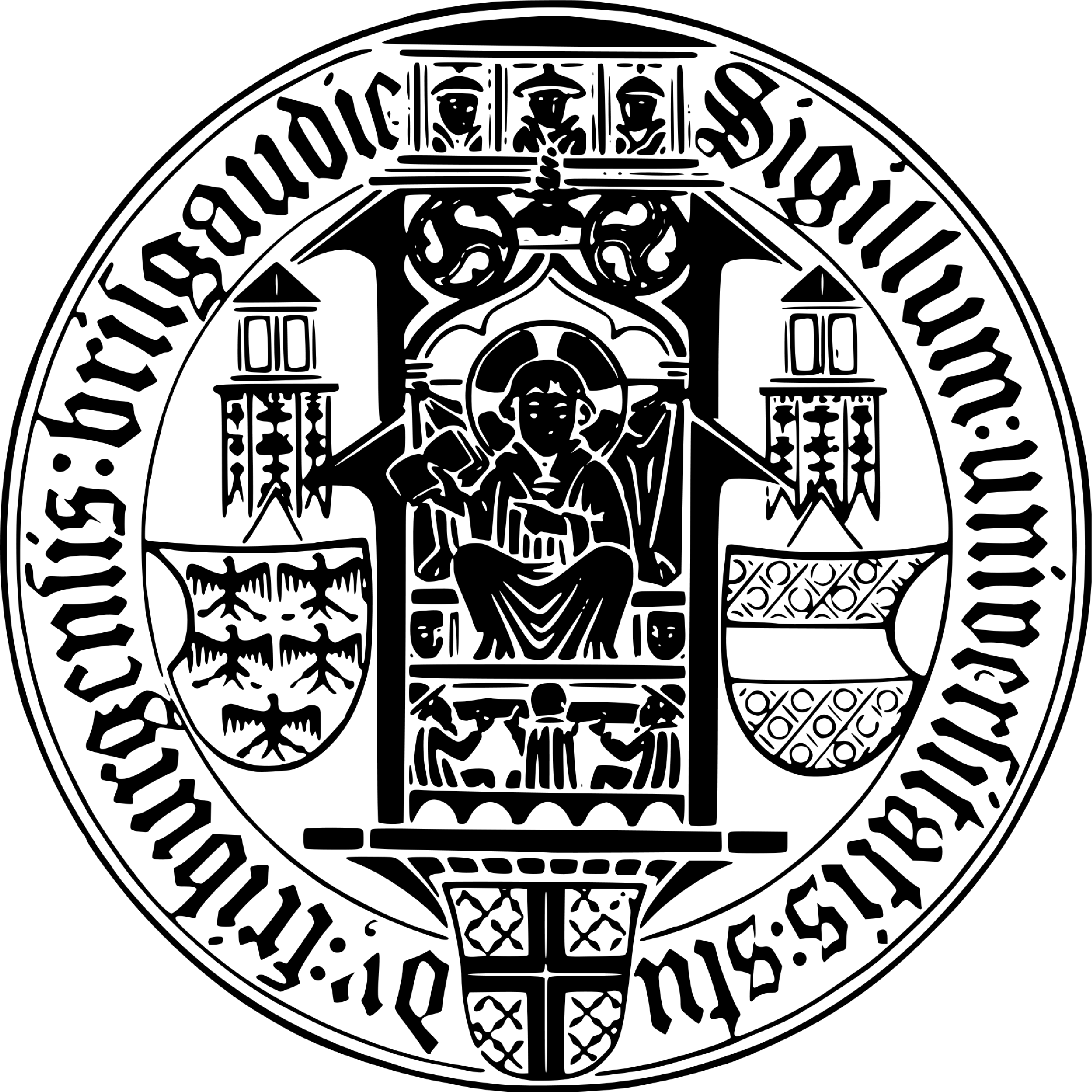}}\\\vspace{3em}
Gradient Flow for the Willmore Functional\\ in Riemannian Manifolds of bounded
Geometry\\\vspace{2em}
\date{\large Dissertation zur Erlangung des Doktorgrades\\ der Fakult\"at f\"ur
Mathematik und Physik\\ der Albert-Ludwigs-Universit\"at Freiburg im Breisgau\\
~\\~\\vorgelegt von Florian Link\\~\\betreut durch Prof. \!Dr.\! Ernst Kuwert\\
~\\August 2013}}\vspace{-15em}
\author{}
\maketitle
\thispagestyle{empty}
\newpage\thispagestyle{empty}
\setlength{\topmargin}{0cm}
\section*{\centering{\sc Abstract}}
\begin{abstract}
  \noindent We consider the $L^2$ gradient flow for the Willmore functional in Riemann\-ian
  manifolds of bounded geometry. In the euclidean case E.\;Kuwert and R.\;Sch\"atzle
  [\textsl{Gradient flow for the Willmore functional,} Comm. Anal. Geom., 10: 307-339,
  2002] established a lower bound of a smooth solution of such a flow, which depends only
  on how much the curvature of the initial surface is concentrated in space. In a second
  joint work [\textsl{The Willmore flow with small initial energy,} J.\;Differential
  Geom., 57: 409-441, 2001] the aforementioned authors proved that a suitable blow-up
  converges to a nonumbilic (compact or noncompact) Willmore surface. In the lecture 
  notes of the first author [\textsl{The Willmore Functional,} unpublished lecture
  notes, ETH Z\"urich, 2007] the blow-up analysis was refined. In the present work we
  intend to generalize the results mentioned above to the Riemannian setting.
  
\end{abstract}

\lefoot{}
\tableofcontents\pagenumbering{roman}
\newpage~
\pagestyle{scrheadings}

\chapter*{Introduction}\lehead{Introduction}\lefoot{\pagemark}
\addcontentsline{toc}{chapter}{Introduction}
In 2002, E.\:Kuwert and R.\:Sch\"atzle considered in their work \it
Gradient Flow for the Willmore Functional\rm~\cite{KS02} two-dimensional compact
immersed surfaces $\Sigma$ in $\R^n$ moving by the gradient of the Willmore functional
$f\mapsto\mathcal W_\circ(f)=\int_\Sigma\bet{A^\circ}^2\mm$, that is, solutions of
$\p_tf=-\grad_{L^2}\mathcal W_\circ$ with $f\rvert_{t=0}=f_0,$ called Willmore flow.
Here, $A^\circ$ is the trace-free
part of the second fundamental form and $\mu$ is the induced area measure on $\Sigma$.
The authors gave a lower bound on the lifespan of a smooth
solution, which depends only on how much the curvature of the initial surface is
concentrated in space. Moreover, they showed that the
curvature concentrates in space if a finite time singularity develops.
In a second joint work \cite{KS01} Kuwert and Sch\"atzle continued studying those 
singularities in greater detail by setting up a blow-up procedure. More precisely, 
they proved apart from other results that a suitable blow-up converges to a nonumbilic
Willmore surface. In his lecture notes \cite{lecKuw}, \mbox{Kuwert} not only gave
a broad summary of the above subjects and other works, but also performed a
somewhat refined analysis concerning the behaviour of the flow near an assumed
singularity.\\

The aim of this thesis is, roughly speaking, to generalize \cite{KS02}, the blow-up
procedure in \cite{KS01} and parts of \cite{lecKuw} in the sense that we 
study the Willmore functional and the corresponding Willmore flow on Riemannian
manifolds
$(M,g)$ of bounded geometry as the target. That is, manifolds with $W^{k,\infty}$ bounds
on the Riemannian curvature tensor and with a strictly positive injectivity radius.
We denote by $H$ and $A$ the mean curvature vector and second fundamental form
respectively.
In $\R^n$ with the standard metric tensor the respective gradient flows of the
functionals induced by the energies $\W_H(f)=\frac 12\int_\Sigma\bet{H}^2\mm,$  $\mathcal
W_A(f)=\frac 12\int_\Sigma\bet{A}^2\mm$ and $\mathcal W_\circ$ are all equivalent since
they differ only by a topological constant. This is in
general no longer true for non-flat targets as it is described in
Section \ref{wfdef} in detail. We decided to study $\W_H,$ that is, the flow
$\p_t f=-\grad_{L^2}\W_H,$ where the $L^2$ gradient refers to the metric tensor $g$ on
$M$.\\

In the first chapter we give a detailed overview of the notation used in this work. In
particular, we adopted the so-called star notation $A\s B$ for tensors $A$ and $B$ used
for example in \cite{HAM82}, \cite{huisken} and \cite{KS02} to describe
nonlinearities in the evolution equations to be developed. In our setting of a non-flat
ambient manifold it is suitable to extend such a notation for terms involving the
Riemannian curvature tensor as it was used in \cite{SF02} for the study of the evolution
of elastic curves in Riemannian manifolds. Later on, we collect basic facts such as
classical differential geometric identities as well as variational formulas for various
geometric objects as in \cite{KS02}. \\\\
Analogously to the approach in \cite{KS02} we prove a lower bound for the maximal time
span of smooth existence of the flow. Chapter two establishes this lower bound in four
steps: First, we compute the evolution of derivatives of the curvature
\begin{eqnarray}
  \na_{\pt}(\na^m\!A)+\Delta^2(\na^m\!A)
    =\pe{m+2}3+\pe m5+\Q {m+2}1{}+\Q m1R,\label{eveq}
\end{eqnarray}
where $\pe kl=\na^{i_1}\!A\s\ldots\s\na^{i_l}\!A~(i_1+\ldots+i_l=k)$
denotes a universal
sum of terms depending multilinearly on derivatives of the
curvature (see Proposition 2.4, \cite{KS02}). The $Q$-terms basically denote
universal sums of terms depending multilinearly on derivatives of the curvature and the
Riemannian curvature tensor $D^{r\!}R.$ The precise algebraic structure is not
being used here. Second, this information will then be transformed into
localized $L^2$-type integral estimates (cf. \cref{lemma34}). To further estimate the
right hand side, i.e. the nonlinearities in (\ref{eveq}), we use interpolation and
absorption techniques. To do so, it is essential to assume that the
concentration
of the curvature is locally small. Third, to interpolate the higher
nonlinearities we need to generalize a Sobolev inequality originally due to Michael and
Simon \cite{MS73} for the case $M=\R^n.$ Up to a little worsening, we prove a
generalization of the aforementioned Sobolev inequality for isometrically immersed
manifolds of any codimension provided the Ricci tensor of $(M,g)$ is bounded in
$L^\infty$ and the injectivity radius is strictly
positive (see \cref{mssi}). In a last step, we employ a multiplicative Sobolev
inequality, which is likewise based on the Sobolev inequality above, to
obtain interior $W^{k,\infty}$-estimates for the curvature still assuming small curvature
concentration
(see \cref{prop47}). Controlling the growth of the maximal local curvature concentration
then allows us to bound from below the first time this smallness condition is
violated.\\\\ 
To state the main theorem of chapter two, that is, the lifespan estimate, we want to
introduce some preliminary definitions. First, we want to define the following
non-local quantity used to control the maximal local energy concentration.
\begin{definition}
  Let $\Sigma$ be a closed manifold, $(M,g)$ be an open or closed Riemannian
  manifold and $I\subset\R$ be an interval. Let further 
  $f_t:(\Sigma,\tildeg_{\!t})\to(M,g)$ be the one-parameter family of isometric
  immersions induced from $f\in C^{2,0}(\Sigma\times I,M),$ where $f_t:=f(\,\c\,,t).$
  If $A_{f_t}$ denotes the second fundamental
  form of $f_t,$ and $\mu_{f_t}$ the induced area measure of $(\Sigma,\tildeg_{\!t}),$ we
  let
  \begin{eqnarray*}
    \chif{f}(r,t):=\sup_{p\in M}
    \int\limits_{f^{-1}_t(\overline{B_r^g(p)})}\!\!\!\!\bet{A_{f_t}}^2d\mue{f_t}\:,
  \end{eqnarray*}
  where $B_r^g(p)\subset M$ denotes the geodesic ball of radius $r$
  and centre p.
\end{definition}
We further define $$\W(f_0):=\W_{H,g}(f_0):=\frac 12\is\bet{H_{f_0}}^2d\mue{f_0}$$
to be the \it{Willmore energy}\rm~ of an isometric $C^2$-immersion 
$f_0:(\Sigma^2,\tildeg_{\!0})\to(M^n,g).$ Let further $R$ be the Riemannian curvature tensor
induced by the Levi-Civit\`a connection of $(M,g),$ and $\mbox{inj}(M,g)$ the injectivity
radius of $(M,g).$ We call $(M,g)$ a Riemannian manifold of \it bounded geometry of order
$k,$\rm~if 
$$\sum_{i=0}^k\norm{D^{i\!}R}{}{L^\infty(M,g)}{}+\mbox{inj}(M,g)^{-1}<\infty.$$ Clearly,
this is automatic (for any $k\in\N_0$) in case $\Sigma$ happens to be compact.
For a more detailed introduction to the notation, we refer to section
\ref{notation}. We prove
the following theorem (see section \ref{lifespanEst}).
\begin{theorem}\rm~(Lifespan estimate for the Willmore flow in
  Riemannian manifolds).\it\label{lifespan}~Given
  an isometric $C^{4+\a}$-immersion $f_0:(\Sigma,\tildeg)\to(M,g)$ of a
  closed surface into a Riemannian n-manifold of bounded geometry (of order 5), 
  we let $f:\Sigma\times [0,T)\to (M,g)$ be a maximal Willmore flow with initial
  datum $f_0.$ For $\varrho>0$ satisfying 
  \begin{eqnarray}
    \varrho\big(\sum_{i=0}^2\norm{D^{i\!}R}{\frac1{i+2}}{L^\infty(M,g)}{}
    +\inj_{(M,g)}^{-1}\big)<c(n)\label{bg2small}
  \end{eqnarray}
  let
  \begin{eqnarray*}
    t^+(\varrho):=\sup\{t\geq 0:~\chi(\varrho,\,\c\,)<\eps_0^2
    \mbox{ on }[0,t)\},
  \end{eqnarray*}
  where $\e_0>0$ and $c>0$ are small universal constants depending only on n. 
  Then either $T=t^+(\varrho)=\infty,$ or for a small constant $C=C(n)>0$
  \begin{eqnarray}
    T>t^+(\varrho)\geq C\varrho^4\log\dfrac{C\eps_0^2}
    {\chi(\varrho,0)
    +\varrho^4\norm{DR}{2}{L^\infty(M,g)}{}
    \big(\mue{f_0}\:(\Sigma)+\varrho^2\W(f_0)\big)}.\label{lifespan2}\bs
  \end{eqnarray}
\end{theorem}\eb\\\notiz{warum besch geo der ord 5? kurze bem machen}
\notiz{concentration in space}
In the third chapter we aim to study singularities of the Willmore
flow in Riemannian manifolds of bounded geometry. Numerical examples
of Mayer and Simonett \cite{MS} indicate the existence of finite
time singularities in euclidean space.
In the same geometric setting Blatt \cite{Blatt} has
shown analytically that singularities may occur at finite or infinite time. However, to
our knowledge there is no analytical evidence for the existence of finite time
singularities. As has been done in \cite{KS01} and
\cite{lecKuw}, we perform a blow-up procedure at an assumed singularity. Using the result
of the lifespan estimate we show that the rescaled flows (sub-)converge to a static
Willmore surface, that is, a static critical
point of $\W.$ For this, we use a compactness theorem that goes
back to Langer \cite{JL} and has then been generalized by Breuning \cite{PB}. It
turns out to be an essential assumption that the area is uniformly bounded along the
flow. Namely, to obtain compactness there are basically two conditions to be
fulfilled. First, the curvature has to be bounded in
$W^{k,\infty}$
and second, a mass-density estimate has to be obeyed. To satisfy the first condition,
we prove a second version of interior $W^{k,\infty}$-estimates for the curvature similar
to those of Kuwert and Sch\"atzle \cite{KS01} using a localization in time (see
\cref{interiorestimates}). The second version allows to perform a blow-up at an assumed
singularity at infinity because the first version (see
\cref{prop47}) depends on initial curvature estimates in contrast to
the second depending on uniform area estimates. The latter dependence is due to an extra
area term in the above mentioned $L^2$-estimates calculated in the second chapter
(cf. also the Remark on page \pageref{simId}). To satisfy the second
condition of the above compactness theorem we establish a mass-density estimate
for immersed surfaces in manifolds of bounded geometry of order one,
at least for radii that are small enough depending on the geometry of the ambient space.
For large radii we were again forced to make the assumption that the area is bounded so
that the large scale estimate is then trivially fulfilled.\\\\
We now want to formulate the main theorem of the third chapter.
To prove this Theorem, a building block is to make the assumption that for a maximal
Willmore flow the total area is uniformly bounded by a number $\mathfrak M_f<\infty$ as
mentioned above.
Fortunately, this requirement
is satisfied for finite singularities and for ambient manifolds of strictly negative
curvature (see the Remark below \cref{blowupthm}).
\notiz{evtl. bem warum die $L^1$-norm von A auftritt $->$ simons inequality?!
warum im $\R^n$ alles ok?}

\begin{theorem}\rm~(Existence of a blow-up).\it\label{blowupthm}~Let 
  $f:\Sigma\times[0,T)\to (M,g)$ be a maximal Willmore flow on a closed surface $\Sigma$
  into a Riemannian manifold of bounded geometry (of order \bg) with the property that
  the total area of $(\Sigma,\tildeg(t))$ is uniformly bounded on $(0,T),$ that is,
  $$\mathfrak M_f:=\sup_{t\in(0,T)}\mue f(\Sigma)<\infty.$$
  Assume that the Willmore flow concentrates at time  $T\in(0,\infty]$ in the sense
  that
  \begin{eqnarray}
    \e_T^2=\lim_{\r\searrow 0}\Big(\limsup_{t\nearrow T}\chi(\r,t) \Big)>0.\label{epst}
  \end{eqnarray}
  Then there exist sequences $t_j\nearrow T$ and $r_j\searrow 0$
  such that the rescaled flows 
  \begin{equation}
    f_j:(\Sigma,\tildeg_{\!j})\times [-r_j^{-4}t_j,r_j^{-4}(T-t_j))\to(M,g_j),
    \qquad f_{j,t}(p):=f_j(p,t):=f(p,t_j+r_j^{4}t),\label{blowupIntro}
  \end{equation}
  where $g_j:=r_j^{-2}g$ and $\tildeg_{\!\!j}(t)=(f_{j,t})^*(g_j),$
  converge after reparametrization locally on
  $\widehat\Sigma\times\R$ to a static solution given by a properly
  immersed Willmore surface $\widehat F_0:\widehat\Sigma\to\R^n$ with
  \begin{eqnarray}
    \int\limits_{\widehat F^{-1}(\overline{B_1(0)})}\bet{A_{\widehat F_0}}^2\mm
    _{\widehat F_0}\geq c(n)>0.\label{nontrivialIntro}
  \end{eqnarray}
  Moreover 
  \begin{eqnarray}
    \e_T^2\geq \int_{\widehat\Sigma}\bet{A_{\widehat F_0}}^2\mm_{\widehat F_0}
    \geq\inf\Big\{\int\bet{A_f}^2\,d\mue f:f\in\mathcal C_n,~f
    \mbox{\rm~is not a union of planes}\Big\}>0,\label{epstIntro}
  \end{eqnarray}
  where $\mathcal C_n$ is the class of all
  properly immersed Willmore surfaces in $\R^n$.\\\\
  Additionally, adopting the terminology from the mean
  curvature flow, a finite time singularity of a Willmore flow of a compact surface
  into a Riemannian manifold of bounded geometry (of order \bg) is always of type II.\eb
\end{theorem}
\notiz{warum wir hier besch. geo der ordnung 15 brauchen.}
For a detailed formulation of the convergence result we refer to
page \pageref{partII}.\\\\
\bf{Remark:}\rm~The total area $\mue f(\Sigma)$ is uniformly bounded on $[0,T),$ if
\begin{itemize}
  \item $T<\infty$ is finite. In this case, it is
  $\mue f(\Sigma)\le\sqrt{2T}\mathcal W(f_0)+\mue{f_0}\:(\Sigma).$ 
  \item the sectional curvature $K^M$ of $(M,g)$ is uniformly negative, i.e. 
  $K^M\le-\hat\kappa^2<0.$
  In this case, it is 
  $\mue f(\Sigma)\le\frac{\mathcal W(f_0)-4\pi\chi(\Sigma)}{2\hat\kappa^2},$
  where $\chi(\Sigma)$ denotes the Euler-characteristic of $\Sigma.$
\end{itemize}
These area estimates hold for any ambient smooth Riemannian manifold $(M,g)$
not necessarily of bounded geometry as we prove in \cref{massbound3}. \er

Finally, in the appendix we have outsourced lengthy computations in coordinates and
have stated well known results concerning harmonic coordinates for convenience of the
reader. Moreover, the interpolation inequalities in \cite{KS02} could almost be carried
over to the Riemannian setting up to minor modifications.\vfill

\LARGE\bf{Acknowledgement}\rm\normalsize\\\\
First, I would like to express my grateful thanks to Prof.\:Dr.\:E.\,Kuwert for his help,
support and not least because of his confidence and patience he brought to me when
drawing up my diploma and doctoral thesis. Further gratitude is dedicated to 
all members of the working group Geometric Analysis in Freiburg for the fruitful
discussions and the pleasant working atmosphere. Also, I would like to thank Dr.\,Maria
Busl, Javier Sabio, PhD, and Felix Klock, PhD, for their help finding orthographical and
typographical errors. Last but not least, I want to thank my girlfriend, writing a
doctoral thesis in sports science herself, for her love and understanding she brought
to me.\\\\
\mbox{~~~~}From January 2007 to December
2010 I was member of the DFG-Forschergruppe Non-linear Partial Differential Equations:
Theoretical and Numerical Analysis. I would like
to thank the DFG for the support.\newpage
\thispagestyle{empty}
\setcounter{equation}{1}

\pagestyle{scrheadings}
\automark[section]{chapter}
\chapter{Foundations}

\pagenumbering{arabic}
\section{Notation and conventions}\label{notation}
In this section we introduce to the notation and conventions used in the present work.
For further reference concerning algebraic facts of Differential Geometry we recommend
\cite{GHV}, or any other good textbook on this topic.

\lehead{\leftmark}
\rohead{\rightmark}
\subsection*{Vector bundles and sections.}
A manifold is understood to be a second countable Hausdorff topological
space together with a smooth ($C^\infty$) structure. Let $B,\Sigma,M,N$ be
manifolds and $I\subset\R$ be a (time) interval. Let further $\xi=(E,\pi_E,B,\R^k)$ and
$\eta=(F,\pi_F,B,\R^m)$ be
vector bundles with total spaces $E$ and $F,$ footpoint projections $\pi_E$ and $\pi_F,$
common Base $B,$ and typical fibres $\R^k$ and $\R^m$ respectively. For simplicity, we
use the symbol $E$ to refer to the bundle $\xi$ etc. if there is no ambiguity. We
let $E_p:=\pi_E^{-1}(\{p \}),$ and with $\xi^{r,s}$ we denote the r-times contravariant
and
s-times covariant tensor product bundle induced by $\xi.$ Given an m-dimensional
manifold $\Sigma^m,$ we denote with
$\T_\Sigma=(T\Sigma,\pi_{\Sigma},\Sigma,\R^m)$ and 
$\T^*_\Sigma=(T^*\Sigma,\pi_{\Sigma},\Sigma,\R^m)$ the tangent and cotangent
bundle respectively, and let $T_p\Sigma:=(T\Sigma)_p$ and
$T^{r,s}\Sigma:=(T\Sigma)^{r,s}$. If $f:\Sigma\to B$ is
smooth, we denote by
$f^*\xi$ the pullback bundle of $\xi$ along f
with total space $\{(p,v)\in M\times E:v\in E_{f(p)} \}$.\\\\
With $\Gamma^k(\xi):=\{\sigma\in C^k(B,E)\mbox{ with }\pi_E\circ\sigma=Id_B\}$ we denote
the $C^k(B)$-module of $C^k$-sections in $\xi.$ We
abbreviate $\Gamma(\xi):=\Gamma^\infty(\xi)$ and
denote with $\Gamma^k_c(\xi)$ the class of sections $\sigma\in\Gamma^k(\xi)$ with compact
support. A $p$-times $C^\infty(B)$-multilinear  (not necessarily alternating) map
$\Gamma(\gt_B)\times\ldots\times \Gamma(\gt_B)\to \Gamma(\xi)$ is called a
\it{$\xi$-valued p-form.}\rm~If $\sigma\in\Gamma(\xi)$ is a section, we let
$f^*\sigma\in\Gamma(f^*\xi),$ defined by $(f^*\sigma):=(id_\Sigma,\sigma\cf),$ be the \it
pullback\rm~of $\sigma$,\notiz{hier weiter suchen ob sections wie $\xi$ richtig}
where $f^*\sigma$ and $\sigma\cf$ will be identified. If $B=\Sigma\times I$ we write
$\sigma_t:=j_t^*\sigma$ for $\sigma\in\Gamma(\xi),$ where $j_t:\Sigma\to\Sigma\times I$
is the inclusion map $j_t(x)\mapsto (x,t)$ opposite $t.$
For arbitrary but fixed $t,$ a section $\sigma\in\Gamma(\xi)$ may always be
considered as a section $\sigma_t\in\Gamma(j_t^*\xi)$ and vice versa provided it is
smooth. Occasionally, the index $t$ will be omitted.\\\\
Clearly, there exists an isomorphism 
$\T_{\Sigma\times I}\simeq \T_\Sigma\times\mathcal T_I.$
To obtain a similar decomposition for $\Gamma(\T_{\Sigma\times I})$ we remark
that also $\T_{\Sigma\times I}\simeq \pi_\Sigma^*\T_\Sigma\oplus\pi_I^*\T_I$ and
$\Gamma(\T_{\Sigma\times I})\simeq\Gamma(\pi_\Sigma^*\T_\Sigma\oplus\pi_I^*\T_I)
\simeq\Gamma(\pi_\Sigma^*\T_\Sigma)\oplus\Gamma(\pi_I^*\T_I),$ where $\pi_\Sigma$
and $\pi_I$ denote the projections onto the first and second factor of
$\Sigma\times I$ respectively. Also,
we may consider $\Gamma(\T_\Sigma)\subsetneq\Gamma(\pi_\Sigma^*\T_\Sigma)\subsetneq
\Gamma(\T_{\Sigma\times I}).$ For example, for any $V\in\Gamma(\T_{\Sigma\times I})$ one
can locally write $V=x^i\p_i\pis\;\,+\sigma\p_t\pii\;\,$ where 
$x^i,\sigma\in C^\infty(\Sigma\times I),$ and $\{\p_i\}$ and $\{\p_t\}$ are the
respective product coordinate frames of $\Sigma$ and $I$ respectively. 
\subsection*{Riemannian metrics, scalar product and Riemannian distance function.}
If $(\xi,g_\xi)$ and
$(\eta,g_\eta)$ are
Riemannian vector bundles, we always understand the tensor product bundle
$\xi\otimes\eta$ to be equipped with the product metric $g_{\xi\otimes\eta}$ defined by
$$g_{\xi\otimes\eta}\strichklein p(z\otimes w,z'\otimes w'):=
g_\xi\strichklein p(z,z')g_\eta\strichklein p(w,w')$$ for $p\in B,z,z'\in E_p,w,w'\in
F_p.$ For a Riemannian vector bundle $(\xi,g_\xi)$ and $T\in\Gamma(\xi)$ we have the
pointwise norms
$\bet{T}^2:=g_\xi(T,T).$ Here and in many other formulas we omit the dependence on
the geometric structure (e.g. the metric $g_\xi$) when there is no
ambiguity.
If $f:(\Sigma,\tildeg)\to(M,g)$ and $T\in\Gamma(T\Sigma\otimes f^*(TM))$ we also
write $\bet{T}^2_f=\bet{T}^2_{\tildeg;g}=T^\a_i T^\b_j \tildeg^{ij}g_{\a\b}\cf$ to
refer to the metrics.

If the basis of $\xi$ is a Riemannian manifold $(B,h)$, we let
for $S,T\in\Gamma^0(\xi),$ and $U\subset B$ measurable
\begin{eqnarray}
  \eck{S,T}_{L^2(\mu_{h})}:=\int\limits_\Sigma g_\xi(S,T)d\mu_{h}
  \mbox{\quad and \quad}
  \norm{T}{2}{2,U}{}:=\int\limits_U g_\xi(S,T)d\mu_{h},\label{scalProdDef}
\end{eqnarray}
where $\mu_h$ is denotes the induced area measure on $B.$
Further for $U\subset B$ $$\norm{T}{}{\infty,U}{}:=\sup\limits_U g_\xi(T,T)^{1/2}$$
provided the right-hand side is finite. 
We write $d_g(p,q)$  for the Riemannian distance induced by $g$ between
$p,q\in (M,g).$\\ 
Let $\Sigma$ be a manifold, $(M,g)$ be Riemannian manifold, a time interval $I$ and
a differentiable map $f:\Sigma\times I\to M.$ If $f_t:\Sigma\to
M,$ where $f_t(x):=f(x,t),$ is an immersion for all $t\in I$, we define
\begin{eqnarray}
  \tildeg_{\!t}:=f_t^*g\mbox{~~~for all~}t\in I\label{fdef}
\end{eqnarray}
making $f_t:(\Sigma,\tildeg_{\!t})\to
(M,g)$ to an isometric immersion.
Here, $f_t^*g\in\Gamma(T^{0,2}\Sigma)$ is the pullback of the tensor field 
$g\in\Gamma(T^{0,2}M)$ defined by $$(f^*_tg)(X,Y):=g\cf_{\!t}(Df_t\c X,Df_t\c Y)$$
for all $X,Y\in\Gamma(T\Sigma).$ More generally, we let
$$\tildeg:= f^*g\strichklein {$\pi_\Sigma^*\T_\Sigma\oplus\pi_\Sigma^*\T_\Sigma$}
  \qquad\qquad\;\,\oplus\delta\strichklein{$\pi_I^*\T_I\oplus\pi_I^*\T_I$}
  \qquad\qquad\;\in\Gamma(\T_{\Sigma\times I}^{0,2}),$$
i.e. $\tildeg(\:\c\:,t)\strichklein {$T^{0,2}\Sigma$}\quad\;\:=\tildeg_{\!t}.$
\subsection*{Connections, adjoint map, projectors, and normal bundle.}
For vector bundles
$(\xi,D')$ and $(\eta,D'')$ with
connections $D^\prime$ and $D^{\prime\prime}$ respectively, we always understand the
product $\xi\otimes\eta$ to be equipped with the
product connection $D^\otimes$ uniquely determined by
$$D^\otimes_X(\sigma\otimes\tau)=(D_X^{\prime}
\sigma)\otimes\tau+\sigma\otimes(D''_X\tau)$$
for all $X\in\Gamma(TB),$ $\sigma\in\Gamma(\xi)$ and $\tau\in\Gamma(\eta).$
In this work, various connections will be denoted with the same symbol. However, it
will be clear from the context which connection is meant.
For $f:\Sigma\times I\to(M,D^M)$ we let $D:=f^*D^M:\Gamma(f^*(TM))\to
\Gamma(T^*(\Sigma\times I)\otimes f^*(TM))$ the \it pullback connection with respect to
$f,$\rm~and write $\na:=P^\bot D$ for the induced connection on the normal
bundle $N_{f}.$ $D$ is the uniquely determined connection on $f^*(TM)$ such that 
$D_X(f^*\sigma)\rvert_{(x,t)}=((x,t),(D_{Df\,\c\, X}\sigma)\rvert_{f(x,t)})$
for all $X\in T_x(\Sigma\times I)$ and $\sigma\in\Gamma(TM).$
When a tensor product connection is used, the factor connections should arise from
the context and/or the symbol itself.\\\\
Let $I\subset\R$ be a (time) interval with standard metric $\delta$ and Christoffel
symbols ${^I}\Gamma\equiv 0.$ Let ${^{\Sigma_t}}\na$ be the Levi-Civit\`a connection on
$\Sigma$ induced by $\tildeg_{\!t}$ with Christoffel symbols ${^{\Sigma_t}}\Gamma$
at time $t$.
We will use the connection $\na$ on $\T_{\Sigma\times I}$ uniquely determined by the
Christoffel symbols $\Gamma(x,t):={^{\Sigma_t}}\Gamma(x)\oplus{^I}\Gamma(t).$
More concretely, in local (product-)coordinates we have for
$V=V^i\p_i\pis\;\,+\sigma\p_t\pii\;\,$ and $W=W^i\p_i\pis\;\,+\tau\p_t\pii\;\,$
$$\na_VW=(VW^k+V^iW^j\;{^{\Sigma_.}}\Gamma_{ij}^k\pis\;\:)\p_k\pis\;+V\tau\pt\pii\;.$$
Clearly, for $f:(\Sigma\times I,\tildeg)\to(M,g)$
$\varphi:=Df\rvert_{\pi^*_\Sigma(T\Sigma)}:\pi^*_\Sigma(T\Sigma)\to f^*(TM)$
is a
strong bundle map and each of the bundles $\pi^*_\Sigma(T\Sigma)$ and $f^*(TM)$ is dual
to itself with respect to the scalar products $\tildeg$ and $g\cf$
respectively.
The Riesz representation theorem guarantees the existence of a
strong bundle map $\Psi:\pi^*_\Sigma(T\Sigma)\leftarrow f^*(TM)$ with the property
$\tildeg(\Psi\tau,\widehat X)=g\cf(\tau,\varphi\,\c\widehat X)$ for all
$\tau\in f^*(TM)$ and $\widehat X\in\T_{\Sigma\times I}.$ Since $\tildeg$ is
non-degenerate, $\Psi^*$ is uniquely
determined and thus we may define $\varphi^*:=\Psi$ to be the dual (or, more
precisely, \it the adjoint\,\rm)~of $\varphi.$
Since $\varphi^*$ is a strong bundle map, this
induces a homomorphism $\varphi^*:\Gamma(f^*(TM))\to\Gamma(\pi^*_\Sigma(T\Sigma)).$
Moreover, we have the decomposition
$$f^*(TM)=\mbox{Im\,}\varphi\oplus \mbox{ker\,}\varphi^*=:T_{\!f}\oplus N_{\!f},$$
where $T_f$ and $N_f$ is the tangent and normal bundle respectively.
Analogously,
we also have $f_t^*(TM)=:T_{f_t}\oplus N_{f_t}.$
Further, we have the tangent and normal
projectors $P:=\varphi\circ\varphi^*\in\mathcal End(f^*(TM))$ and $P^\bot:=Id-P$ and in
particular, using that $\varphi_t=Df_t$, we have that
$P_t=Df_t\circ (Df_t)^*\in\mathcal End(f_t^*(TM)).$
Now for $V\in N_{f_t}$ we may consider $\na V:\Gamma(T\Sigma)\to\Gamma(f_t^*(TM))$
that analogously induces for any $t\in I$ the adjoint map
$(\na V)^*:\Gamma(T\Sigma)\leftarrow\Gamma(f_t^*(TM))$ with the respective
property $\tildeg_{\!t}((\na V)^*\Psi_t,X)=g\cf_t(\Psi_t,\na_{X}V)$
for any $\Psi\in f^*(TM)$ and $X\in\Gamma(T\Sigma).$\\\\

With $\na^*$ we denote the formal adjoint of the operator $\na$ with respect to the
scalar products defined in (\ref{scalProdDef}).
More precisely, for $\phi\in\Gamma(T^{0,s+1}\Sigma\otimes N_{f_t})$ we
let $\na^*T$ be the unique tensor field with
$$\int\limits_Mg(\na^*\phi,\psi)d\mu_g=\int\limits_Mg(\phi,\na\psi)d\mu_g$$
for all $\psi\in\Gamma(T^{0,s}\Sigma\otimes N_{f_t}).$ It
is $\na^*\psi=-(\na_{e_i}\psi)(e_i,\ldots).$
We define the (normal) Laplacian by $\Delta:=-\na^*\na,$ i.e. for $\psi$ as above it 
is $\Delta\psi=(\na^2\psi)(e_i,e_i,\ldots).$ Terms involving operators have to be read
from the right to the left, e.g. for a normal valued 1-form $\phi$ and
$X,Y,Z\in\gg(T\Sigma)$ we write
$\na_X\na_Y\phi(Z)\equiv\na_X\big(\na_Y(\phi(Z))\big),$
\mbox{$\na_X Df\c Y\equiv\na_X(Df\c Y),$} or
for $\xi\in\gg(TM)$ we write $D_X \xi\cf\equiv D_X(\xi\cf)=(D\xi)\cf\c Df\c X.$

\subsection*{Curvature.}
Let $f$ be as in (\ref{fdef}) and $(M,g,D)$ be equipped with
the Levi-Civit\`a connection $D.$
We write \mbox{$A_f:=D^{2\!}f\in\Gamma(T^{2,0}\Sigma\otimes N_f)$} and
$H_f\in\Gamma(N_f)$ for the second fundamental form and mean curvature vector of f
respectively. Further we let $A_f^\circ:=A_f-\frac 12\tildeg\otimes H_f$ be the tracefree
part of the bilinear form $A_f.$ The occasional omission of the index $t$ should here
and elsewhere be obvious. 
For a vector bundle $(\xi,D)$ with connection $D$ we
define the curvature $F\in\Gamma(\Lambda^2(TM)\otimes\mbox{End}(E))$ of
$(E,D)$ by
$$F(X,Y)\xi:=D_XD_Y\xi-D_YD_X\xi-D_{[X,Y]}\xi$$
for all $X,Y\in\Gamma(TM)$ and $\xi\in\Gamma(E).$
We then let $R$ be the curvature of $(TM,D),$ $R_f$ be the curvature of
$(f^*(TM),f^*D)$ and $R_{\Sigma\times I}$ be the curvature of $(T(\Sigma\times I),\na).$
Finally, we let $\RB{}$ and $\RB l$ denote the curvature of
$(N_f,\na),$ $(T^{0,l}\Sigma\otimes N_f,\na)$ respectively. Note, that the
restriction of
$R_{\Sigma\times I}\rvert_{(x,t)}$ to $T_x\Sigma$ equals the curvature
$R_{\Sigma_t}\rvert_x$ of $(T\Sigma,{}^{\Sigma_t}\na).$
Last, we let $\w K_t:=R_{\Sigma_t}(e_1,e_2,e_2,e_1)$ be the \it{Gaussian curvature}\rm~of
$(\Sigma,\tildeg_{\!t})$ and let $$K(T\Sigma):=R\cf(Df\c e_1,Df\c e_2,Df\c e_2,Df\c
e_1)$$
be the \it{sectional curvature of $T\Sigma$ in M,}\rm~where again $\{e_i\}_{i=1,2}$ is a
local $\tildeg_{\!t}$-orthonormal basis of $T\Sigma.$ That $\w K_t$ and $K(T\Sigma)$
is well defined, i.e. the independence with respect to the chosen local orthonormal
basis, follows from the orthogonality of the transformation map of two given
orthonormal bases and the symmetries of the curvature tensor.

\subsection*{Miscellanea.}\label{misc}
Here, we summarize the notion of various symbols and geometric quantities used in this
work.
\begin{itemize}
  \item All manifolds are understood to be second countable topological Hausdorff spaces
        together with a smooth structure. We only consider manifolds without boundary.
        For emphasis, a manifold without boundary is called \it open,\rm~if
        it is noncompact, and likewise \it closed,\rm~in case it is compact. Riemannian
        manifolds are understood to be complete with respect to the induced metric. 
  \item $\mbox{inj}(M,g)$: injectivity radius of $(M,g)$
  \item $\mu_{\tildeg}:$ induced area-measure on $(\Sigma,\w g)$
  \item Zero indexed quantities such as $A_0,~f_0,~\bet{\;\c\;}_0,~\mu_0$ etc. refer to 
        time $t=0$
  \item For a section $X\in\gg\big(T\Sigma \big)$ we occasionally abbreviate
        $\w X:=Df\c X$
  \item For a cutoff function $\g$ on a manifold N, we define the set\\
    $[\g>0]:=\{p\in N,\mbox{ such that }\g(p)>0\}$ and
    $[\g=1]$ etc. analogously to $[\g>0].$
  \item We use Einstein's convention, i.e. summation over repeated indices is used
  \item $g\cf$ is sometimes abbreviated by $\langle \;\c\;,\;\c\; \rangle$ 
  \item For $\phi\in\gg\big(N_f \big)$ we define
    $Q(A^\circ)\phi:=A^\circ(e_i,e_j)\langle A^\circ(e_i,e_j),\phi\rangle$
\end{itemize}\notiz{funtionenr\"aume definieren: $C^{k+\a,l+\b}$
und $C$}
Let $\Sigma$ be a manifold without boundary, $N$ be a manifold and $I\subset\R$ be an
interval. Then $f\in C^0(\Sigma\times I,M)$ is called \it locally proper\rm~if and only
if $$f\strichklein{$[t_1,t_2]$}\;\;\quad:\Sigma\times [t_1,t_2]\to M\mbox{ is
proper for any compact }[t_1,t_2]\subset I.$$
Analogously, $\eta\in C^0(\Sigma\times I)$ has \it locally compact support\rm~if and
only if
$$\eta\strichklein{$[t_1,t_2]$}\;\;\quad:\Sigma\times [t_1,t_2]\to\R\mbox{ has 
compact support for any compact }[t_1,t_2]\subset I.$$ Clearly, each of the conditions is
automatic if $\Sigma$ is closed. If a function is bounded and locally
proper then it obviously has locally compact support. 
\begin{definition}{\rm (Parabolic H\"older spaces, see \cite{LA01}).}\\
  Let $(x_1,t_1),~(x_2,t_2)\in\R^k\times \R.$ Then 
  $$d_m\big((x_1,t_1),(x_2,t_2) \big)
  :=\max\{\bet{x_1-x_2},\bet{t_1-t_2}^{\frac1{2m}}\}$$
  defines a metric on $\R^k\times \R.$ For $G\subset\R^k\times \R$ and a
  function $u:G\to\R^l$ we define the parabolic $\a$-H\"older coefficient
  $$[u]_{\a,m,G}
  :=\sup\Big\{\frac{\bet{u(p)-u(q)}}{d_m(p,q)^\a}:p,q\in G,~p\neq q\Big\}.$$
  Further we define $$C^{2m,1}(G)
  :=\big\{u:G\to\R^l:u,Du,D^2u,\ldots,D^{2m}u,\p_t u\in C^0(G)\big\}$$
  and $$C^{2m,1,\a}(G):=\big\{u\in C^{2m,1}(G):[D^{2m} u]_{\a,m,G}
  +[\p_t u]_{\a,m,G}<\infty\big\}.$$
\end{definition}
For a definition of H\"older spaces on manifolds, see \cite{wloka}.

\subsection*{Star- and Q-notation.}
\notiz{only real or normal valued?? - check}
Here, we orientate towards the notation used in \cite{KS02}, \cite{HAM82} and
\cite{SF02}. For an isometric immersion
$f:(\Sigma^2,\tildeg)\to(M,g)$ and $\eta,\psi\in\gg\big(T^{r,s}\Sigma\otimes (N_f)^{p,q}
\big)$ we denote by
\begin{eqnarray*}
  \eta\s \psi
\end{eqnarray*}
normal- (or real-) valued sections such that with respect to any local orthonormal frame
we have the form
\begin{eqnarray}
  \big(\eta\s \psi \big)^{\a,I}_{\b,J}
  =C^{\a,\d,\sigma,I,L,P}
    _{\b,\g,\r,J,K,R}\;\eta^{\g,K}_{\d,L}\;\psi^{\r,R}_{\sigma,P}\label{sternnotation}
\end{eqnarray}
where $\a,\b,I,J,\ldots$ are multi-indices, and for any fixed term 
the coefficients $C^{\dots}_{\dots}$ are constants depending only on $r,s,p,q$
and the dimension n.
With this definition one easily verifies that
\begin{eqnarray}
  \bet{\eta\s \psi}\le C(r,s,p,q,n)\bet{\eta}\bet{\psi}\label{starabsch}
\end{eqnarray}
and
\begin{eqnarray}
  \na(\eta\s\psi)=(\na\eta)\s\psi+\eta\s(\na\psi),\label{diffrule}
\end{eqnarray}
where $\na$ in each case stands for the respective connection. To check
the pointwise identity (\ref{diffrule})
we may for any local orthonormal frame and any fixed point $x\in\Sigma$ assume
that the respective connection one-forms are equal to zero. Then (\ref{diffrule})
follows from the product rule for $\eta^{\g,K}_{\d,L}\;\psi^{\r,R}_{\sigma,P}$.\\\\
Remarks:
\begin{enumerate}
  \item [a)] Of course, this can be generalized to more than two factors.
  \item [b)] From the way we have defined $\eta\s\psi,$ it is not only clear that 
             $\eta\s\psi=\psi\s\eta,$ but also that the indices can be interchanged
             within each multi-index.
  \item [c)] Although the star-terms do not depend on other sections than stated, we 
             sometimes omit 
             the arguments to keep overview as long as everything is clear from the
             context.\er
\end{enumerate}\notiz{c) sollte weg}

As mentioned in the introduction, it is useful to introduce a notation to encode
nonlinearities involving the Riemannian curvature tensor as in \cite{SF02} for the study
of elastic curves. It turned out that, roughly speaking, the necessary estimates are
performed in the normal bundle and therefore make the split-up $f^*(TM)=T_f\oplus N_f.$
Consider the split up 
$$P^\bot R\cf(D\eta,\w X)\phi=P^\bot R\cf(\na\eta,\w X)\phi
+P^\bot R\cf(P^\top\! D\eta,\w X)\phi.$$
If the tangent bundle has rank one as in \cite{SF02} the last summand above clearly
vanishes by symmetry properties of the curvature tensor. In our setting, where $T\Sigma$
is of rank two the notation unfortunately needs to be modified as follows:\\

Analogously to the above definition, we let for 
$\zeta,\xi\in\gg\big(T^{r,s}\Sigma\otimes\*(T^{p,q}M) \big)$ 
\begin{eqnarray}
  \zeta\star \xi\label{sternchen}
\end{eqnarray}
denote 
sections with the respective property (\ref{sternnotation}). Analogously to the above
we now have the estimate 
$\bet{\zeta\star\xi}\le C(r,s,p,q,n)\bet{\zeta}\bet{\xi}$ and
$D(\eta\star \psi)=(D\zeta)\star \xi+\zeta\star (D\xi).$
Here, $\*(TM)$ is understood to be equipped with the pull-back connection
$\*D$ induced from $(M,g,D)$ in contrast to the above definition where we consider $N_f$
to be equipped with the normal connection $\na=P^\bot D
\begin{picture}(14.3,1)\put(0,-1.8){\big\lvert}
  \put(3,-5){\footnotesize $N_f$}\end{picture}$.\\\\
{\underline{Examples:}} For $X,Y,Z\in\gg(T\Sigma),$ and
$\zeta,\xi,\varPsi\in\gg(\*(TM)),$ and if $\{e_i\}_{i=1,2}$ is a
local $\tildeg-$orthonormal frame for $T\Sigma,$ we have 
\begin{enumerate}
  \item [i)] $g\cf\big(H,(\na A)(X,Y,Z) \big)
    =(A\s \na A)(X,Y,Z)$
  \item [ii)] $P^\top R\cf(\zeta,\xi)\varPsi=R\cf(\zeta,\xi,\varPsi,\w e_i)\w e_i
      =(R\cf\star Df\star Df)(\zeta,\xi,\varPsi)$
  \item [iii)] $P^\bot R\cf (\w X,\w Y)\w Z$
\begin{eqnarray}
    &=&R\cf (\w X,\w Y)\w Z
        -R\cf (\w X,\w Y,\w Z,\w e_i)\w e_i\no\\
      &=&(R\cf \star Df\star Df\star Df
        +R\cf \star Df\star Df\star Df\star Df\star Df)(X,Y,Z).\label{rterm}
  \end{eqnarray}\er
\end{enumerate}
With $\pe kl$ we denote universal linear combinations of terms of the form
$$\na^{i_1}\!A\s \ldots\s \na^{i_l}\!A$$
when $\bet{i}:=i_1+\ldots+i_l=k.$ As usual, if the sum is empty, i.e. $l=0$ then
$\bet{i}:=0.$
With this definition one then easily verifies that
\begin{eqnarray}
  \na\pe kl=\pe{k+1}l.\label{pabsch}
\end{eqnarray}
It turns out to be useful (as (\ref{rterm}) illustrates)
to encode lower order terms involving the curvature tensor arising in the evolution
equations. To do so,
we denote by $\Q kl{(m)}$ normal- (or real-) valued universal linear combinations of
sections
\begin{eqnarray}
  \big(D^r\!R \big)\cf\star\na^{i_1\!}\!A\star \ldots\star\na^{i_\nu\!}\!A\star
  \iota_N \star\ldots\star\iota_N\star Df\star 
    \ldots\star Df\label{118}
\end{eqnarray}
when $r+\bet{i}+\nu=k+l,$ $\bet{i}\le k~(,$ and $r\le m$ in case the lower index $m$ is
given). Here, $\iota_N:(N_f,\na)\to(\*(TM),D)$ denotes the canonical injection.
Analogously, we denote by $\Q klR$ normal (or real) valued universal linear
combinations of sections
\begin{eqnarray}
  \big(D^{r_1}\!R \big)\cf\star \big(D^{r_2}\!R \big)\cf
     \star\na^{i_1\!}\!A\star \ldots\star\na^{i_\nu\!}\!A\star
     \iota_N \star\ldots\star\iota_N\star Df\star \ldots\star Df
\end{eqnarray}
when $r_1+r_2+\bet{i}+\nu=k+l$ and $\bet{i}\le k.$ For $\Q klR$ we additionally want
to assume that $\nu\geq 1$ for technical reasons.
To see what happens to the indices when differentiating covariantly
we firstly note that for a normal valued n-form $\phi$ and adapted vector fields
$X_i\in\Gamma(T\Sigma)$ we pointwise have
\begin{eqnarray}
  \hspace{-4em}(D_X\phi)(X_1,\ldots,X_n)&=&D_X\gro{\phi(X_1,\ldots,X_n)}\no\\
    &=&\na_X\gro{\phi(X_1,\ldots,X_n)}
      +\eck{D_X\gro{\phi(X_1,\ldots,X_n)},Df\c e_j}Df\c e_j\no\\
    &=&(\na_X\phi)(X_1,\ldots,X_n)-\eck{\phi(X_1,\ldots,X_n),A(X,e_i)}
      Df\c e_i\label{aufspaltungvn}
\end{eqnarray}
and thus we may write
\begin{eqnarray}
  D\phi=\na\phi+\phi\star A\star Df.\label{tangnormsplit}
\end{eqnarray}
Secondly, for $X\in\gg(T\Sigma)$ and $\phi\in\gg(N_f)$ we have
$(D\iota_N)(X,\phi)=-\eck{\iota_N(\phi),A(X,e_j)}Df\c e_j,$ i.e.
\begin{eqnarray*}
  D\iota_N=\iota_N\star A\star Df.
\end{eqnarray*}
Thus we easily get, recalling that $D(Df)=D^{2\!}f=A,$ 
\begin{eqnarray}
  \na\Q 00{}=\Q 01{}\mbox{  and  }\na\Q kl{(m)}=\Q {k+1}l{(m+1)}.\label{qdiff}
\end{eqnarray}
{\underline{Examples:}} For (\ref{rterm}) we may now write 
$P^\bot R\cf (\w X,\w Y)\w Z=\Q 00{(0)}(X,Y,Z).$ Further, we have, for example, the
simplifications $P^\bot R\cf (\na^i\!A,\w Y)\w Z=\Q i10(X,Y,Z)$ and
\begin{eqnarray*}
  P^\bot R\cf(H,Df\:\c\:)\phi&=&\big((P^\bot R\cf)\star A\star\iota_N\star Df\big)
  \s\phi\\
  &=&\big(R\cf\star A\star\iota_N\star Df+
  R\cf\star A\star\iota_N\star Df\star Df\star Df\big)
  \s\phi\\
  &=&\Q 010\s\phi.\\
\end{eqnarray*}
Recalling the definition, it is easy to see that instead of $\Q kl{}\s A$ we can write 
$\Q kl{}\star A=\Q k{l+1}{}.$\er\\
\mbox{}Summarizing, it is important to keep the estimate (\ref{starabsch}) and the
differentiation rules (\ref{pabsch}) and (\ref{qdiff}) in mind.
\section{Basic facts}
\subsection*{Fundamental identities of differential geometry.}
In this subsection, we summarize already known facts for convenience of the reader and 
derive easy facts used later on.\\\\
\begin{lemma}\label{isoembedd}
  Let $f:(\Sigma^d,\tildeg)\to(M^n,g)$ and $I:(M^n,g)\to(\R^N,\d_{eucl})$ be
  isometric immersions. Then there holds the pointwise estimates
\begin{eqnarray}
  \bet{A_{I\cf}}\le c\bet{A_I}\cf+\bet{A_f}\label{secfundabsch}
\end{eqnarray}
and
\begin{eqnarray}
    \bet{H_{I\cf}}\le c\bet{A_I}\cf+\bet{H_f},\label{meancurvabsch}
\end{eqnarray}
where $c=c(d,n)$ is a universal constant.
\end{lemma}
\bb 
Differentiating $D(I\circ f)=DI\cf\,\c\,Df$ covariantly, we get for vector fields
\mbox{$X,Y\in\Gamma(T\Sigma)$}
\begin{eqnarray}
  A_{I\cf}(X,Y)&=&D^2_{X,Y}(I\circ f)=(D^2I)\circ f(Df\c X,Df\c Y)
    +(DI)\circ f\,\c\,D^2f(X,Y)\no\\
  &=&A_I\circ f(Df\c X,Df\c Y)+(DI)\cf\c A_f(X,Y).\label{indhyp}
\end{eqnarray}
Since $I$ is isometric, the first inequality follows since
\begin{eqnarray*}
  \bet{A_{I\cf}}^2&=&\bet{A_I\cf(Df\:\c\:,Df\:\c\:)}^2+\bet{A_f}^2\\
  &\le&c\bet{A_I}^2\cf\bet{Df}^4+\bet{A_f}^2
\end{eqnarray*}
because $\bet{Df}^2=d.$ For the second, we choose a local $\tildeg-$orthonormal frame
$\{e_i\}_{1\le i\le d}$ and obtain analogously
\begin{eqnarray*}
  \bet{H_{I\cf}}^2&=&\bet{\sum_i A_I\cf(Df\c e_i,Df\c e_i)}^2+\bet{H_f}^2\\
  &\le&\sum_i\bet{A_I}^2\cf\bet{Df}^4\bet{e_i}^4+\bet{H_f}^2\\
  &\le&c\bet{A_I}^2\cf+\bet{H_f}^2.
\end{eqnarray*}\eb\\
For an arbitrary isometric immersion $f:(\Sigma,\tildeg)\to (M,g)$ and any vector
fields
$X,Y,Z,W\in\Gamma(T\Sigma),~\phi\in\Gamma(N_f)$ the following identities hold 
(see, e.g. \cite{Carmo})
\begin{eqnarray}
  R\cf(\w X,\w Y,\w Z,\w W)&=&R_\Sigma(X,Y,Z,W)\no\\ 
    &&-\eck{A(X,W),A(Y,Z)}+\eck{A(X,Z),A(Y,W)}\label{135}\\
  D_X(Df\c Y)&=&A(X,Y)+Df\c\na_X Y\label{136}\\
  &&~~~~~~~~~~~~~~~~~~~~(\mbox{Equations of Gau\ss})\no\\
  \RB{}(X,Y)\phi&=&P^\bot R\cf(\w X,\w Y)\phi+A(e_i,X)\eck{A(e_i,Y),\phi}\no\\
    &&-A(e_i,Y)\eck{A(e_i,X),\phi}~~~~\label{137}\\
  &&~~~~~~~~~~~~~~~~~~~~(\mbox{Ricci identity})\no\\
  P^\bot R\cf(\w X,\w Y)\w Z&=&(\na A)(X,Y,Z)-(\na A)(Y,X,Z)\label{138}.\\
  &&~~~~~~~~~~~~~~~~~~~~(\mbox{Equation of Mainardi-Codazzi})\no
\end{eqnarray}
The Gaussian curvature $\w K$ of $(\Sigma^2,\tildeg),$ and the sectional
curvature $K(T\Sigma)$ of $T\Sigma$ in $(M,g)$ are related as follows
\begin{eqnarray}
  \w K-K(T\Sigma)&=&\frac 12\big(\bet H^2-\bet A^2 \big)\label{139b}\\
  &=&\frac12\bet A^2-\bet{A^\circ}^2\label{140b}\\
  &=&\frac14\bet H^2-\frac12\bet{A^\circ}^2.\label{140c}
\end{eqnarray}
Namely, from the Gau\ss~equation (\ref{135}) it follows with 
$A_{ij}:=A(e_i,e_j)$ that
\begin{eqnarray*}
  \w K-K(T\Sigma)&=&R_\Sigma(e_1,e_2,e_2,e_1)
  -R\cf(\w {e_1},\w{e_2},\w{e_2},\w {e_1})\\
  &=&\eck{A_{11},A_{22}}-\bet{A_{12}}^2.\label{141}\\
  &=&\frac12\big(\eck{A_{11}+A_{22},A_{11}+A_{22}}-(\bet{A_{11}}^2
  +2\bet{A_{12}}^2+\bet{A_{22}}^2) \big)\\
  &=&(\ref{139b}).
\end{eqnarray*}
From the orthogonal decomposition $A=A^\circ\operp_{g\cf}(\frac12\tildeg\otimes H)$ in
the trace and traceless part of the bilinear form A we get
$$\bet{A}^2=\bet{A^\circ}^2+\frac 14\bet{H}^2\bet{\tildeg}^2=\bet{A^\circ}^2
+\frac12\bet{H}^2,$$
since $\bet{\tildeg}^2=\,$dim$\,\Sigma=2.$ From this, (\ref{140b}) and (\ref{140c})
follow.\\\\
Since we have to interchange second covariant derivatives of normal l-forms, we need to 
know how the curvature of $(T^{0,l}\Sigma\otimes N_f,\na)$ is related to the curvature of
$(N_f,\na)$ and $(T\Sigma,{}^\Sigma\na).$ We pointwise have for any 
$\phi\in \Gamma\left(T^{0,l}\Sigma\otimes N_f\right)$
\begin{eqnarray*}
    (\na_{X,Y}^{2}\phi)(X_1,\ldots,X_l)
    &=&\na_X\gro{(\na_Y\phi)(X_1,\ldots,X_l)}\\
    &=&\na_X\gro{ \na_Y\kle{\phi(X_1,\ldots,X_l)}}
        -\sum\limits^{l}_{k=1}\phi(X_1,\ldots,\na^2_{X,Y}X_k,\ldots,X_l)\\
    &=&\na^2_{X,Y}\gro{\phi(X_1,\ldots,X_l)}
        -\sum\limits_{k=1}^{l}\phi(X_1,\ldots,\na^2_{X,Y}X_k,\ldots,X_l)
\end{eqnarray*}
yielding
\begin{eqnarray*}
 \lefteqn{\hspace{-3em}\big(R^{\,l}_\bot(X,Y)\phi\big)(X_1,\ldots,X_l)}\\
    &=&(\na^2_{X,Y}\phi)(X_1,\ldots,X_l)
      -(\na^2_{Y,X}\phi)(X_1,\ldots,X_l)\\
    &=&\RB{}(X,Y)\phi(X_1,\ldots,X_l)
       -\sum\limits_{k=1}^{l}\phi\big(X_1,\ldots,R_\Sigma(X,Y)X_k,\ldots,X_l\big).
\end{eqnarray*}
Since
\begin{eqnarray*}
  \w K= K(T\Sigma)+\dfrac{1}{2}\big(\bet{H}^2-\bet{A}^2 \big)=\Q 00{} +A\s A
\end{eqnarray*}
and
\begin{eqnarray*}
  R_\Sigma(X,Y)Z=\w K\big(\tildeg(Y,Z)X-\tildeg(X,Z)Y\big),
\end{eqnarray*}
we get by substitution
\begin{eqnarray*}
  \phi\big(X_1,\ldots,R_\Sigma(X,Y)X_k,\ldots,X_l\big)\\
  &&\hspace{-10em}
    =\w K\big(\tildeg(Y,X_k)\phi(X_1,\ldots,X_{k-1},X,X_{k+1},\ldots,X_l)\\
  &&~~\hspace{-10em}
    -\tildeg(X,X_k)\phi(X_1,\ldots,X_{k-1},Y,X_{k+1},\ldots,X_l) \big)\\
  &&\hspace{-10em}=\big(\Q 00{}\s\phi+A\s A\s\phi\big)(X,Y,X_1,\ldots,X_l).
\end{eqnarray*}
Putting things together, we get
\begin{eqnarray}
  \RB l(\,\cdot\,,\,\cdot\,)\phi&=&A\s A\s\phi\no
     +\Q 00{}\s\phi+\RB{}(\,\cdot\,,\,\cdot\,)
     \phi(\,\cdot\,,\ldots,\,\cdot\,)\\
    &\stackrel{(\ref{137})}{=}&A\s A\s\phi+\Q 00{}\s\phi.\label{143aa}
\end{eqnarray}

\dipl{Lemma 1.5}
\begin{lemma}\label{15}
  For arbitrary $\phi\in\Gamma\big(T^{0,l}\Sigma\otimes N_f\big)$ and $l\geq 1$ we have
\begin{eqnarray*}
    (\na\na^*-\na^*\na)\phi=-\na^*T+A\s A\s\phi+\Q 00x\s\phi,
  \end{eqnarray*}
where
$T(X_0,X_1,\ldots,X_l):=(\na\phi)(X_0,X_1,\ldots,X_l)-(\na\phi)(X_1,X_0,\ldots,X_l).$
\end{lemma}
\bb

We pointwise compute for adapted vector fields
\begin{eqnarray*}
  \lefteqn{\hspace{-3em}\big((\na\na^*-\na^*\na)\phi \big)(X_1,\ldots,X_l)}\\
  &=&\na_{e_i}\gro{(\na\phi)(e_i,X_1,\ldots,X_l)}
      -\na_{X_1}\gro{(\na\phi)(e_i,e_i,X_2,\ldots,X_l)}  \\
  &=&(\na^2_{e_i,e_i}\phi)(X_1,\ldots,X_l)
    -(\na^2_{e_i,X_1}\phi)(e_i,X_2,\ldots,X_l) \\
  &&+(\na^2_{e_i,X_1}\phi)(e_i,X_2,\ldots,X_l)
    -(\na^2_{X_1,e_i}\phi)(e_i,X_2,\ldots,X_l) \\
  &=&-(\na^*T)(X_1,\ldots,X_l)+\big(\RB l(e_i,X_1)\phi \big)(e_i,X_2,\ldots,X_l). \\
  &\stackrel{(\ref{143aa})}{=}&(-\na^*T+A\s A\s\phi
    +\Q 00{}\s\phi)(X_1,\ldots,X_l).
\end{eqnarray*}\bs\\
We need two special cases of the latter lemma.
First, if we set $\phi:=A$ we get from Codazzi-Mainardi
\begin{eqnarray*}
  T(X,Y,Z)&=&(\na A)(X,Y,Z)-(\na A)(Y,X,Z)\\
  &=&P^\bot R\cf(\w X,\w Y)\w Z \\
  &=& \Q 00x(X,Y,Z).
\end{eqnarray*}

On the other hand, again using Codazzi-Mainardi we pointwise have
\begin{eqnarray}
  -(\na^*\!A)(Y)&=&(\na A)(e_i,e_i,Y)=\na_{e_i}\gro{A(e_i,Y)}
    =\na_{e_i}\gro{A(Y,e_i)}=(\na A)(e_i,Y,e_i)\no \\
  &=&(\na A)(Y,e_i,e_i)+P^{\bot}R\cf(\w e_i,\w Y)\w e_i\no \\
  &=&\na_YH+\Q 00x(Y).\label{144}
\end{eqnarray}
From this we obtain pointwise
\begin{eqnarray*}
  (\na \na ^*\!A)(X,Y)&=&\na _X\gro{(\na ^*\!A)(Y)} \\
  &\stackrel{(\ref{144})}{=}&-\na_X\kle{\na_Y H}-\na _X\gro{\Q 00x(Y)} \\
  &=&-\na ^2_{X,Y}H+\Q 01y(X,Y).\label{145}
\end{eqnarray*}
Now since $\na ^*T=1\s \na T=1\s\na\Q 00x=\Q 01y,$ \cref{15} shows that
\begin{eqnarray*}
  (\na \na ^*-\na ^*\na )A&=& A\s A\s A+\Q 01y+\Q 00x\s A,\\
  &=&\pe 03+\Q 01y
\end{eqnarray*}
since we can write $\Q 00x\s A$ as $\Q 01y.$
\notiz{rough version of simons und bem machen wegen problemen mit DR}
Eventually we get a rough version of Simons' identity
\begin{eqnarray}
  \Delta A&=& \na ^2H+P^0_3(A)+\Q 01y.\label{146}
\end{eqnarray}
\bf{Remark:}\label{simId}\rm~Following the computations thoroughly it turns out that
the difference
of $\Delta A-\na^2H$ contains an additive term with the Riemannian curvature tensor that
does \it not\rm~contain the second fundamental form. This fact results in the existence
of an extra area term when we interpolate the nonlinearities to establish $L^2$ estimates
(cf. \cref{prop44} and \cref{prop45}).\er\\
Second, we want to substitute $\phi$ by $\na \phi$ in \cref{15}:
To start with, we have
\begin{eqnarray*}
  T(X_0,\ldots,X_l)&=&(\na^2_{X_0,X_1}\phi)(X_2,\ldots,X_l)
     -(\na^2_{X_1,X_0}\phi)(X_2,\ldots,X_l)\\
  &=&\big(\RB{l-1\,}(X_0,X_1)\phi\big)(X_2,\ldots,X_l) \\
  &\stackrel{(\ref{143aa})}{=}&( A\s A\s \phi+\Q 00x\s\phi)(X_0,\ldots,X_l).
\end{eqnarray*}
Now we differentiate
\begin{eqnarray*}
  \na ^*T=1\s \na T=A\s \na A\s\phi+A\s A\s \na \phi+\Q 00x\s \na \phi+\Q 01y\s\phi
\end{eqnarray*}

so that we again get from \cref{15}
\begin{eqnarray}
  (\Delta \na -\na \Delta)\phi&=&(\na \na ^*-\na ^*\na )\na \phi\no \\
  &=&A\s \na A\s\phi+A\s A\s \na \phi+\Q 00x\s \na \phi+\Q 01y\s\phi.\label{149}
\end{eqnarray}

\dipl{Lemma 1.6}
\begin{lemma}
  Let $f\in C^{2,1}\left(\Sigma\times (t_1,t_2),M\right)$ be a variation with normal 
  velocity field $V:=\p_t f.$ Then for $\xi\in\Gamma^2(\*(TM)),~\phi\in\Gamma^2(N_f)$
  and time independent vector fields X,Y,Z we have
  \begin{align}
   && D_{\partial_{t}}^{}\kle{Df\c X}&= D_{X}^{}V&&&\label{152}\\
   &&R_f(\partial_{t}^{},X)\xi&=R\cf(V,\widetilde X)\xi&&&\label{153}\\
   && D_{\pt} P&=- D_{\partial_{t}}^{}(P^{\bot})
      =Df\c(\na V)^{*}+ \na V\c(Df)^{*}  &&&\label{154}\\
   && P D_{\partial_{t}}^{}\phi&=-Df\c (\na V)^{*}\phi  &&&\label{155}\\
   &&  \RB{}(\partial_{t}^{},X)\phi&=
      \langle A(X,e_i),\phi\rangle\na_{e_i}^{}V-\langle\na_{e_i}^{}V,\phi\rangle A(X,e_i)
      +P^{\bot} R\cf(V,\widetilde X)\phi&&&\label{156}\\
   &&(\na_{\pt}\tildeg)(X,Y)&=-2\langle A(X,Y),V\rangle   &&&\label{157}\\
   &&\partial_{t}^{}(\mm)&=-\langle H,V\rangle\mm   &&&\label{158}\\
   &&R_{\Sigma\times I}(\pt,X)Y&=\na_{\partial_{t}}^{}\kle{\na_{X}Y}&&&\no\\ 
   &&&=-\langle (\na A)(X,Y,e_i),V\rangle e_i
      +\langle A(X,Y),\na_{e_i}^{}V\rangle e_i  &&&\label{159}\\
   &&&\quad-\langle A(X,e_i),\na_{Y}^{}V\rangle e_i
     -\langle A(Y,e_i),\na_{X}^{}V\rangle e_i
     - R\cf(\widetilde Y,\widetilde {e_i},\widetilde X,V)e_i&&&\no\\
   &&(\na_{\pt}A)(X,Y)&=\na^{2}_{X,Y}V-\langle V,A(Y,e_i)\rangle A(X,e_i)
   +P^{\bot} R\cf(V,\widetilde X)\widetilde Y\label{160}\\
   &&\na_{\partial_{t}}^{}H&=\Delta V+Q(A^{o})V
      +\frac 12 H\langle H,V\rangle
      +P^{\bot} R\cf(V,\widetilde e_i)\widetilde{e_i}.&&&\label{161}
   \end{align}
\end{lemma}
\bb \\
{\bf (\ref{152}):} By properties of the connection along f, we have
\begin{eqnarray*}
  D_{\pt}\kle{Df\c X}-D_X\kle{Df\c\pt}&=&Df\c[\pt,X]\\
  &=&0
\end{eqnarray*}
since $X$ does not depend on time.\\

{\bf (\ref{153}):} Again by locality and since $[\pt,X]=0,$ we can compute
\begin{eqnarray*}
  R_f(\pt,X)\kle{\pa\cf}&=&D_{\pt}\gro{D_{X}\kle{\pa\cf}}
      -D_{X}\gro{D_{\pt}\kle{\pa\cf}}-D_{[\pt,X]}\kle{\pa\cf} \\
      &=&D_{\pt}\gro{(D\pa)\cf\c Df\c X}-D_{X}\gro{(D\pa)\cf\c V} \\
      &\stackrel{(\ref{152})}{=}&(D^2\pa)\cf(V,\w X)-(D^2\pa)\cf(\w X,V)\\
      &=&R\cf(V,\w X)(\pa\cf).
\end{eqnarray*}

{\bf (\ref{154}):}
At first, let $\psi\in\Gamma(N_f).$ Then we have for any vector field
$X\in\Gamma(T\Sigma)$
\begin{eqnarray*}
 \eck{D_{\pt}\psi,Df\c X}&=&\pt\eck{\psi,Df\c X}-\eck{\psi,D_{\pt}\kle{Df\c X}}\\
 &\stackrel{(\ref{152})}{=}&-\eck{\psi,\na_X V}\\
 &=&-\eck{\psi,\na_{e_i}V}\eck{Df\c e_i,Df\c X}
\end{eqnarray*}
and thus
\begin{eqnarray}
 PD_{\pt}\psi&=&-Df\c\eck{\psi,\na_{e_i}V}e_i.\label{162}
\end{eqnarray}
Now let $\xi\in\Gamma(f^*(TM)).$ Decomposing $\xi=\psi+Df\c Y$ in its tangential and 
normal part, we can compute for arbitrary $\phi\in\Gamma(N_f)$
\begin{eqnarray*}
 \eck{D_{\pt}\kle{P^{\bot}\xi},\phi}&=&\pt\eck{P^{\bot}\xi,\phi}
   -\eck{P^\bot\xi,D_{\pt}\phi}\\
 &=&\pt\eck{\xi,\phi}-\eck{P^\bot\xi,D_{\pt}\phi}\\
 &=&\eck{D_{\pt}\xi,\phi}+\eck{P\xi,D_{\pt}\phi}\\
 &=&\eck{D_{\pt}\xi,\phi}+\eck{P\xi,PD_{\pt}\phi}\\
 &\stackrel{(\ref{162})}{=}&\eck{D_{\pt}\xi,\phi}
   -\eck{P\xi,Df\c\eck{\phi,\na_{e_i}V}e_i }\\
 &=&\eck{D_{\pt}\xi,\phi}-\eck{\na_Y V,\phi}\\
 &=&\eck{P^\bot D_{\pt}\xi-\na_Y V,\phi}
\end{eqnarray*}
On the other hand, we get from (\ref{152})
\begin{eqnarray*}
  \eck{D_{\pt}\kle{P^\bot \xi},Df\c X}&=&\pt\eck{P^\bot\xi,Df\c X}
    -\eck{P^\bot\xi,D_{\pt}\kle{Df\c X}}\\
  &=&-\eck{P^\bot\xi,\na_X V}\\
  &=&-\eck{\xi,\na_{e_i}V}\eck{Df\c e_i,Df\c X}\\
  &=&\eck{-Df\c\eck{\xi,\na_{e_i}V}e_i,Df\c X}.
\end{eqnarray*}
Therefore
\begin{eqnarray*}
  D_{\pt}\kle{P^\bot\xi}&=&P^\bot D_{\pt}\xi-\na_{Y}V-Df\c \eck{\xi,\na_{e_i}V}e_i,
\end{eqnarray*}
and thus, using $(D_{\pt}P^\bot)(\xi)=D_{\pt}\kle{P^\bot\xi}-P^\bot D_{\pt}\xi,$ 
the second identity follows. The first follows
trivially from $P+P^\bot=Id.$\\\\
{\bf (\ref{155}):} We have
\begin{eqnarray*}
  PD_{\pt}\phi&=&-(D_{\pt} P)\phi\\
  &\stackrel{(\ref{154})}{=}&-Df\c(\na V)^*\phi-\na V\c(Df)^*\phi.\\
  &=&-Df\c(\na V)^*\phi.
\end{eqnarray*}
{\bf (\ref{156}):} Again since $[\pt,X]=0$ we can compute
\begin{eqnarray*}
  \RB{}(\pt,X)\phi&=&\na_{\pt}\kle{\na_{X}\phi}
    -\na_X\kle{\na_{\pt}\phi}-\na_{[\pt,X]}\phi\\
  &=&P^\bot D_{\pt}\kle{P^\bot D_X\phi}-P^\bot D_X\kle{P^\bot D_{\pt}\phi}\\
  &\stackrel{(\ref{155})}{=}&P^\bot(D_{\pt} P^\bot)\kle{D_X\phi}
    +P^\bot D_{\pt}\kle{D_X\phi}
    -P^\bot D_{X}\big(D_{\pt}\phi+Df\c(\na V)^*\phi \big)\\
  &\stackrel{(\ref{153}),(\ref{154})}{=}&-\na V\c(Df)^*(D_X\phi)+P^\bot R\cf(V,\w X)\phi
    -A(X,(\na V)^*\phi)\\
  &=&\eck{\phi,A(X,e_i)}\na_{e_i}V+P^\bot R\cf(V,\w X)\phi-\eck{\phi,\na_{e_i}V}A(X,e_i),
\end{eqnarray*}
where we used $(\na V)^*\phi=\eck{\phi,\na_{e_i}V}e_i$ and
\begin{eqnarray*}
  (Df)^*(D_X\phi)&\stackrel{(\ref{aufspaltungvn})}{=}&(Df)^*(\na_X\phi
    -\eck{\phi,A(X,e_i)}Df\c e_i) \\
  &=&-\eck{\phi,A(X,e_i)}e_i
\end{eqnarray*}
in the last step.\\\\
{\bf (\ref{157}):} For the evolution of the metric we get since $X$ and $Y$ do not 
depend on time
\begin{eqnarray*}
  (\na_{\pt}\tildeg)(X,Y)
  &=&\pt\gro{\tildeg(X,Y)}\\
  &=&\pt\eck{Df\c X,Df\c Y}\\
  &=&\eck{D_{\pt}\kle{Df\c X},Df\c Y}+\eck{Df\c X,D_{\pt}\kle{Df\c Y}}\\
  &\stackrel{(\ref{152})}{=}&\eck{D_X V,Df\c Y}+\eck{Df\c X,D_Y V}\\
  &=&X\eck{V,Df\c Y}-\eck{V,A(X,Y)}+Y\eck{Df\c X,V}-\eck{A(X,Y),V}\\
  &=&-2\eck{A(X,Y),V}.
\end{eqnarray*}
{\bf (\ref{158}):} Using the standard formula for the derivative of the determinant, 
i.e.\\
$\frac{d}{dt}\ln\det(\tildeg_{\!ij})=\tildeg^{ij}\pt \tildeg_{\!ij},$ we get
$$\dfrac{d}{dt}\sqrt{\det(\tildeg_{\!ij})}=\frac 12\tildeg^{ij}\pt \tildeg_{\!ij}
\sqrt{\det(\tildeg_{\!ij})}$$
and thus, since $\pt\tildeg_{\!ij}=(\na_{\pt}\tildeg)(\partial_i,\partial_j)
\stackrel{(\ref{157})}{=}-2\eck{A_{ij},V},$ finally
$$\frac d{dt}\sqrt{\det(\tildeg_{\!ij})}=-\eck{H,V}\sqrt{\det(\tildeg_{\!ij})}.$$

{\bf (\ref{159}):} Since $\na_{\pt}\kle{\na_{X}Y}$ is $C^\infty(\Sigma)-$linear in $X$
and 
$Y,$ we may assume that $\na_XY=\na_X e_i=0$ at a fixed point and at a fixed time.
\begin{eqnarray*}
  Df\c\na_{\pt}\kle{\na_XY}
  &=&D_{\pt}\kle{Df\c\na_XY}-D^2f(\pt,\na_XY) \\
  &=&D_{\pt}\gro{PD_X\kle{Df\c Y}}\\
  &=&(D_{\pt}P)\gro{D_{X}\kle{Df\c Y}}+PD_{\pt}\gro{D_{X}\kle{Df\c Y}}\\
  &\stackrel{(\ref{153})}{=}&(D_{\pt}P)\gro{A(X,Y)}
    +PD_X\gro{D_{\pt}\kle{Df\c Y}}+PR\cf(V,\w X)\w Y\\
  &=&Df\c(\na V)^*A(X,Y)+PD_{X}\kle{D_{Y}V}+R\cf(V,\w X,\w Y,\w{e_i})Df\c e_i,
\end{eqnarray*}
where we used (\ref{152}) and (\ref{154}). Since
$(\na V)^*\big(A(X,Y)\big)=\eck{A(X,Y),\na_{e_i}V}e_i$
and
\begin{eqnarray*}
  PD_X\kle{D_YV}&\stackrel{(\ref{aufspaltungvn})}{=}&
    PD_X\big(\na_YV-\eck{V,A(Y,e_i)}Df\c e_i \big)\\
  &\stackrel{(\ref{aufspaltungvn})}{=}&-\eck{\na_YV,A(X,e_i)}Df\c e_i
    -X\eck{V,A(Y,e_i)}Df\c e_i\\
  &=&-\eck{\na_YV,A(X,e_i)}Df\c e_i-\eck{\na_XV,A(Y,e_i)}Df\c e_i\\
  &&-\eck{V,(\na A)(X,Y,e_i)}Df\c e_i,
\end{eqnarray*}
the claim follows by the injectivity of $Df$ on $T\Sigma.$\\\\
{\bf (\ref{160}):} Again assuming $\na_X Y=0$ at a fixed point and a fixed time, this 
follows from
\begin{eqnarray*}
  (\na_{\pt}A)(X,Y)&=&\na_{\pt}\gro{A(X,Y)}\\
  &=&P^{\bot}D_{\partial_t}\big(D_{X}\kle{Df\c Y}-Df\c\na _XY\big) \\
  &\stackrel{(\ref{153})}{=}&P^{\bot}D_X\gro{D_{\partial_t}\kle{Df\c Y}}
    +P^{\bot}R\cf(V,\w X)\w Y \\
  &\stackrel{(\ref{152})}{=}& \na^2_{X,Y}V+P^{\bot}D_X\kle{PD_YV} 
    +P^\bot R\cf(V,\w X)\w Y
\end{eqnarray*}
using
\begin{eqnarray*}
  P^{\bot}D_X\kle{PD_YV}&=&-P^{\bot}D_X\big(\eck{V,A(Y,e_i)}Df\c e_i \big)\\
  &=&-\eck{V,A(Y,e_i)}A(X,e_i).
\end{eqnarray*}

{\bf (\ref{161}):} From (\ref{157}) we get in local coordinates
\begin{eqnarray*}
  0&=&\partial_t(g_{ij}g^{jk})=-2\langle A_{ij},V\rangle g^{jk}+g_{ij}\p_tg^{jk}\\
    \Rightarrow~~\partial_tg^{lk}&=&2\langle A_{ij},V \rangle g^{jk}g^{il} \\
    \Rightarrow~~(\p_tg^{lk})A_{lk}&=&2\langle A_{ij},V \rangle g^{jk}g^{il}A_{lk}\\
  &=&2Q(A)V.
\end{eqnarray*}

From (\ref{160}) we thus get
\begin{eqnarray*}
  \na_{\partial_t}H&=&\na_{\partial_t}(g^{lk}A_{lk}) \\
  &=&\partial_tg^{lk}A_{lk}+g^{lk}\na_{\partial_t}A_{lk} \\
  &=&2Q(A)V+\Delta V-Q(A)V+P^\bot R\cf(V,\w {e_i})\w {e_i}\\
  &=&\Delta V+Q(A^\circ)V+\frac{1}{2}\eck{H,V}H+P^\bot R\cf(V,\w{e_i})\w {e_i},
\end{eqnarray*}
since $Q(A)V=Q(A^\circ)V+\frac{1}{2}\eck{H,V}H.$ 
\eb
\setcounter{chapter}{1}
\setcounter{section}{0}
\setcounter{equation}{0}
\chapter{Lifespan Theorem}
\section{The Willmore flow}\label{wfdef}
\subsection*{Flat Situation.}
In \cite{KS02} the Willmore energy of a closed, isometrically immersed surface
$f:(\Sigma,\bar g)\to(\R^n,\d_{eucl})$
with induced area measure $d\mu$ is defined by
\begin{eqnarray}
  \W_\circ(f)=\int\limits_{\Sigma}\bet{A^\circ}^2d\mu.
\end{eqnarray}
Where $A^\circ=A-\frac 12\bar g\otimes H$ is the tracefree part of the second
fundamental form $A.$ The Gau\ss~equations and Gau\ss-Bonnet imply
\begin{eqnarray}
  \W_{\circ}(f)=\W_{A}(f)-2\pi\chi(\Sigma)
  =\W_{H}(f)-4\pi\chi(\Sigma),\label{energyidentity}
\end{eqnarray}
where 
$\W_{A}(f)=\frac 12\int_\Sigma\bet A^2d\mu$ and $\W_{H}(f)
=\frac 12\int_\Sigma\bet H^2 d\mu.$
Thus in this flat situation, the $L^2$ gradient flows of all these functionals
reduce to the same since they only differ by a topological constant. In general,
this is not true if the target manifold is curved.
\subsection*{Situation of a curved ambient manifold.}\label{curvedsit}
Now let $f:(\Sigma^2,\tildeg)\to(M^n,g)$ be an isometrically $C^2-$immersed closed
surface into an n-dimensional Riemannian Manifold. 
From (\ref{140b}) and (\ref{140c}) we have
$$\bet{A^\circ}^2=\frac12\bet A^2+K(T\Sigma)-\w K
=\frac12\bet H^2+2K(T\Sigma)-2\w K.$$
Integration and using Gau\ss-Bonnet yields
\begin{eqnarray}
  \W_{\circ,g}(f)&=&\W_{A,g}(f)+\int\limits_\Sigma
    \!K(T\Sigma)\,d\mu-2\pi\chi(\Sigma)\no\\
  &=&\W_{H,g}(f)+2\int\limits_\Sigma \!K(T\Sigma)\,d\mu-4\pi\chi(\Sigma)
    \label{nonequivfunc}
\end{eqnarray}
We decided to study $\W_{H,g}$ in this work.
\subsection*{$\bf L^2$-gradient flow for the Willmore functional.}
Let $f$ be as in \ref{curvedsit}. We define the Willmore energy (with respect to
$H$ and $g$
in f) as
\begin{eqnarray}
  \W(f):=\W_{H,g}(f):=\frac 12\int\limits_\Sigma\bet H^2d\mu.
\end{eqnarray}
In the next Lemma we compute the $L^2$ gradient of this functional for convenience of
the reader.
\dipl{Lemma 1.2}
\begin{lemma}\rm(\cite{dipl}, Lemma 1.2).\it~The
 $L^2$ gradient for the Willmore functional (in f) is given by
  \begin{eqnarray}
    \bf W\rm\it(f):=\grad_{L^2}\W(f)=\Delta_{}H+Q(A^\circ)H+P^\bot R\cf(H,\w e_i)\w e_i.
  \end{eqnarray}
\end{lemma}
\bb Let $f:\Sigma\times (-\e,\e)\to M$ be a variation of f with normal 
velocity-field $\pt f=:\phi.$ Then for the first variation of $\W$ of $f$ in direction
$\phi$ we have
\begin{eqnarray*}
  \frac d{d\e}\W(f_\e)\strichklein{$\,\e=0$}\quad\,
  &=&\frac d{d\e} \is\frac 12\bet {H_\e}^2d\mu\strichklein{$\,\e=0$}\\
  &\stackrel{(\ref{158})}{=}&
     \is\eck{H,\na_{\partial_\e}H_{\e}\strichklein{$\,\e=0$}\quad\;\;\:}d\mu
    -\frac 12\is\bet H^2\eck {H,\phi}d\mu\\
  &\stackrel{(\ref{161})}{=}&\is\eck{H,\Delta\phi +Q(A^\circ)\phi+\frac 12\eck{H,\phi}H
    +P^\bot R\cf(\phi,\w{e_i})\w {e_i}}d\mu\\
  &&-\frac 12\is\bet H^2\eck {H,\phi}d\mu.
\end{eqnarray*}
Since $\Delta$ and $Q(A^\circ)$ are self-adjoint with respect to 
$\int\eck{\;\c\;,\,\c\,}\mm,$ we further get
\begin{eqnarray*}
  &=&\is\eck{\Delta H+Q(A^\circ)H,\phi}d\mu+\is R\cf(\phi,\w{e_i},\w {e_i},H)d\mu\\
  &=&\is\eck{\Delta H+Q(A^\circ)H+P^\bot R\cf(H,\w{e_i})\w {e_i},\phi}d\mu\\
  &=&\eck{\Delta H+Q(A^\circ)H+ P^\bot R\cf(H,\w{e_i})\w {e_i},\phi}_{L^2}\\
  &=&\eck{grad_{L^2}\W(f),\phi}_{L^2},
\end{eqnarray*}
by definition of the $L^2$ gradient.\\\\
\bf{Remark:}\rm~The reason why we may restrict to normal variations is the following:
For tangential variations there always exists a flow on $\Sigma$ that generates this
tangential vector field. It is plausible that this flow only causes a reparametrization
of the surface so that the Willmore energy leaves unchanged due to its geometric
invariance of the underlying parametrization. Strictly speaking, one can show that
$\W(f\circ\psi)=\W(f)$ for all $C^2$-diffeomorphisms  $\psi:\Sigma\to\Sigma$
(see \cite{lecKuw}).\er
\begin{theorem}{ \rm(Existence and uniqueness).}~Let $f_0:(\Sigma^2,\tildeg)\to
(M^n,g)$ be an isometric $C^{4+\a}$-immersion of a
closed surface into a Riemannian manifold. Then the initial value problem 
\begin{eqnarray*}
  \pt f&=&-\mbox{\bf W} (f)\quad\mbox{on }\Sigma\times (0,T)\\
  f\strichklein{t=0}\;\,&=&f_0
\end{eqnarray*}
has a unique solution $f\in C^{4,1,\a}(\Sigma\times[0,T),M)$ on a maximal time interval
$[0,T)$ where $0<T\le\infty.$ Moreover, the restriction of f to
  $\Sigma\times (0,T)$ is smooth.\eb
\end{theorem}
The solution to the above initial value problem will be called \it{Willmore flow}~\rm for
short.  

\section{Evolution equations}
\dipl{Lemma 2.1}
\begin{lemma}\label{lemma21}
  Let $\phi\in\Gamma\big(T^{0,l-1}\Sigma\otimes N_f \big)$ be a normal (l-1)-Form along 
  a variation $f:\Sigma\times I\rightarrow M$ with normal velocity field
  $\partial_t
  f=:V.$ If $\na_{\partial_t}\phi+\Delta^2\phi=:Y,$ then for $\psi:=\na\phi$ we have
\begin{eqnarray*}
  \na_{\partial_t}\psi+\Delta^2\psi
    &=&\na Y+\sum\limits_{i+j+k=3}^{}\big(\na^iA\s \na^j A\s \na^k\phi\big)
      +A\s \na V\s \phi+V\s \na A\s\phi\\      
    &&+\Q 00x\s V\s\phi
      +\sum\limits_{\substack{j+k=1} }^{}\Delta\big(\Q 0jy\s \na^k\phi)
      +\sum\limits_{\substack{j+k=1} }^{}\Q 0jy\s\na^{k+2}\phi.\\
\end{eqnarray*}
\end{lemma}

\bb
Let $X_1,\ldots,X_l$ be time independent vector fields such that $\na X_k=0$
in a fixed point $(x_0,t_0).$ Then we pointwise have
\begin{eqnarray*}
  (\na_{\partial_t}\psi)(X_1,\ldots,X_l)
    &=&\na_{\partial_t}\gro{(\na_{X_1}\phi)(X_2,\ldots,X_l)} \\
    &=&\na_{\partial_t}\gro{\na_{X_1}\kle{\phi(X_2,\ldots,X_l)}}-\na_{\partial_t}
      \sum\limits_{k=2}^{l}\phi(X_2,\ldots,\na_{X_1}X_k,\ldots,X_l) \\
    &=&\na_{\partial_t}\gro{\na_{X_1}\kle{\phi(X_2,\ldots,X_l)}}
      -\sum\limits_{k=2}^{l}\phi(X_2,\ldots,\na_{\partial_t}\kle{
      \na_{X_1}X_k},\ldots,X_l).
\end{eqnarray*}
On the one hand we get for the first term, using (\ref{156}),
\begin{eqnarray*}
  \lefteqn{\hspace{-3em}\na_{\partial_t}\gro{\na_{X_1}\kle{\phi(X_2,\ldots,X_l)}}} \\
  &=&\na_{X_1}\gro{\na_{\partial_t}\kle{\phi(X_2,\ldots,X_l)}}
    +(A\s \na V\s\phi+\Q 00x\s V\s \phi)(X_1,\ldots,X_l) \\
  &=&\big(\na(\na_{\partial_t}\phi)\big)(X_1,\ldots,X_l)
   +(A\s \na V\s\phi+\Q 00x\s V\s \phi)(X_1,\ldots,X_l),
\end{eqnarray*}
recalling that
  $P^\bot R\cf(V,Df\:\c\:)\phi=\big((P^\bot R\cf)\star Df\star\iota_N \star\iota_N
  \big)\s V\s\phi=\Q 00{}\s V\s\phi$ and\\
  $P^\bot R\cf=R\cf+R\cf\star Df\star Df.$ For the second term, we have from (\ref{159})
\begin{eqnarray*}
  \phi(X_2,\ldots,
  \na_{\pt}\kle{\na_{X_1}X_k},\ldots,X_l)
  &=&\big(V\s \na A\s\phi+A\s \na V\s\phi\big)(X_1,\ldots,X_l).
\end{eqnarray*}
Putting things together we get
\begin{eqnarray}
  \na_{\partial_t}\psi+\Delta^2\psi-\na Y&=&\Delta^2(\na \phi)
    -\na (\Delta^2\phi)+A\s \na V\s \phi+V\s \na A\s\phi
  +\Q 00x\s V\s\phi.~~~~~~~\label{21}
\end{eqnarray}
On the other hand, we get from (\ref{149})
\begin{eqnarray}
  (\Delta \na -\na \Delta)\phi&=&A\s \na A\s\phi+A\s A\s \na \phi
    +\Q 00x\s \na \phi+\Q 01y\s\phi\label{23} \\
  \Rightarrow~~\Delta\big((\Delta \na -\na \Delta)\phi\big)
  &=&\!\!\!\sum\limits_{i+j+k=3}^{}\!\!\big(\na ^iA\s \na ^jA\s \na ^k\phi\big)
    +\!\!\!\sum\limits_{\substack{i+j+k=1} }^{}
    \!\!\!\Delta\big(\Q 0jy\s\na^k\phi\big),\label{24}
\end{eqnarray}
where we used the differentiation rule (\ref{qdiff}) for the Q-Terms. Substituting 
$\phi$ by $\Delta\phi$ in\\ (\ref{23}) implies
\begin{eqnarray}
  (\Delta \na -\na \Delta)(\Delta\phi)&=&A\s \na A\s\Delta\phi+A\s A\s \na \Delta\phi
  +\Q 00x\s \na \Delta\phi+\Q 01y\s \Delta\phi\no \\
  &=&A\s \na A\s \na ^2\phi+A\s A\s \na ^3\phi
    +\sum\limits_{j+k=1}^{}\Q 0jy\s \na ^{k+2}\phi,\label{25}
\end{eqnarray}
since terms like $\na\Delta\phi$ can be written as a contraction of 
$\tildeg^{-1}\otimes\na^3\phi.$ Adding up (\ref{24}) and (\ref{25}) yields
\begin{eqnarray*}
  \Delta^2(\na \phi)-\na (\Delta^2\phi)
  &=&\Delta\big(\Delta(\na \phi)-\na (\Delta\phi)\big)
     +(\Delta \na -\na \Delta)(\Delta\phi) \\
  &=&\4\sum\limits_{i+j+k=3}^{}\3\big(\na ^iA\s \na ^jA\s \na ^k\phi\big)
     +\2\sum\limits_{\substack{j+k=1} }^{}\2\Delta\big(\Q 0jy\s \na^k\phi)
     +\2\sum\limits_{\substack{j+k=1} }^{}\2\Q 0jy\s\na^{k+2}\phi.
\end{eqnarray*}
Substituting this into (\ref{21}) yields the claim.\eb

\dipl{Lemma 2.2}
\begin{lemma}\label{lemma22} For the Willmore flow $\partial_tf=V=-{\bf W}(f)$ 
we have the evolution equations
\begin{eqnarray}
  \lefteqn{\hspace{-2em}\na_{\partial_t}(\na^{m}A)+\Delta^2(\na^{m}A)}\label{26} \\
  &=&\pe {m+2}3+\pe m5
    +\na^m\Delta\Q 01y+\Delta\Q m1m+\na\Q {m+1}1{m+1}+\Q{m+2}1{m+1}+\Q m1R\no\\
  &=&\pe {m+2}3+\pe m5+\Q {m+2}1{}+\Q m1R\no
\end{eqnarray}
for any $m\in\mathbb N_0.$
\end{lemma}
\bb (Induction over $m\in\mathbb N_0$). Substituting $\phi$ by $\na H$ or $H$ in 
(\ref{23}) we get
\begin{eqnarray*}
  (\Delta \na-\na\Delta)(\na H)&=&A\s \na A\s \na H+A\s A\s \na ^2H
    +\sum\limits_{j+k=1}^{}\Q 0jy\s \na ^k\na H \\
  &=&P^2_3(A)+\Q  211
\end{eqnarray*}
and
\begin{eqnarray*}
  \na \big((\Delta \na -\na \Delta)H \big)&=&\na \big(A\s \na A\s H+A\s A\s \na H
   +\sum\limits_{j+k=1}^{}\Q 0jy\s \na ^kH \big)\\
  &=&\na \big(P^1_3(A)+\Q 111 \big) \\
  &=&P^2_3+\na\Q 111.
\end{eqnarray*}
Combining yields
\begin{eqnarray}
  \Delta \na ^2H-\na ^2\Delta H&=&(\Delta \na -\na \Delta)(\na H)
    +\na \big((\Delta \na -\na \Delta )H\big)\no \\
  &=&\pe 23+\na\Q 111+\Q 211,\label{27}
\end{eqnarray}
and substituting $V=-\Delta H+\pe03+\Q010$ in (\ref{160}) we obtain 
\begin{eqnarray}
  \na_{\pt}A&=&\na^2V+A\s A\s V+\Q00x\s V\no\\
  &=&-\na^2(\Delta H+\pe03+\Q010)+A\s A\s(\Delta H+\pe03+\Q010)\no\\
  &&+\Q00x\s (\Delta H+\pe03+\Q010)\label{glng1}\\
  &=&-\na^2(\Delta H)+\pe23+\na^2\Q 010+\pe05+\Q030+\Q210+\Q01R\no\\
  &\stackrel{(\ref{27})}{=}&-\Delta(\na^2H)+\pe23+\pe05+\na\Q 111+\Q211+\Q01R\no\\
  &\stackrel{(\ref{146})}{=}&-\Delta(\Delta A+\pe03+\Q01y)+\pe23+\pe05
    +\na\Q 111+\Q211+\Q01R\no\\
  &=&-\Delta^2A+\pe23+\pe05+\Delta\Q 01y+\na\Q 111+\Q 211+\Q 01R.\label{28}
\end{eqnarray}
To proceed with the induction step we use (\ref{149}) with $\phi:=\Q m1m$ and
obtain
\begin{eqnarray}
  \na\Delta\Q m1m&=&\Delta\na\Q m1m+A\s\na\!A \s\Q m1m+A\s A\s\na\Q m1m+\Q 00x\s\na\Q m1m
    +\Q 01y\s\Q m1m\no\\
  &=&\Delta\Q {m+1}1{m+1}+\Q {m+1}3m+\Q {m+1}3{m+1}+\Q {m+1}1R+\Q m2R\no\\
  &=&\Delta\Q {m+1}1{m+1}+\Q {m+3}1{m+1}+\Q {m+1}1R.\label{commutate1}
\end{eqnarray}
Letting $\psi:=\na\phi$ for $\phi:=\na^mA,$ \cref{lemma21} now yields  
\begin{eqnarray*}
  \lefteqn{\hspace{0em}\na_{\pt}(\na^{m+1}A)+\Delta^2(\na^{m+1}A)}\\
   &=&\na \big(\na_{\pt}(\na^mA)+\Delta^2(\na^mA) \big)
     +\!\!\!\!\sum\limits_{i+j+k=3}^{}\!\!\!\!\big(\na^iA\s \na^j A\s\na^{k+m}A\big)\\
   &&+A\s (\na\Delta H+\pe13+\Q111)\s \na^m A+(\Delta H
     +\pe03+\Q010)\s \na A\s\na^m A\\
   &&+\Q00x\s(\Delta H+\pe03+\Q010)\s\na^m A
     +\3\sum\limits_{\substack{j+k=1} }^{}\Delta\big(\Q 0jy\s \na^{k+m}\!A)
      +\3\sum\limits_{\substack{j+k=1} }^{}\Q 0jy\s\na^{k+m+2}\!A\\
   &=&\pe{m+3}3+\pe{m+1}5+\na^{m+1}\Delta\Q 01y
     +\na\Delta\Q m1m+\na^2\Q {m+1}1{m+1}+\na\Q {m+2}1{m+1}+\na\Q m1R\\
   &&+\Q {m+1}31+\Q {m+1}30
     +\Q {m+2}20+\Q m40+\Q m2R+\Delta\Q{m+1}11+\Q{m+3}11\\
   &\stackrel{(\ref{commutate1})}{=}&\pe{m+3}3+\pe{m+1}5
     +\na^{m+1}\Delta\Q 01y+\Delta\Q {m+1}1{m+1}+\na\Q {m+2}1{m+2}+\Q {m+3}1{m+2}
     +\Q {m+1}1R,
\end{eqnarray*}
which proves the induction step and thus the lemma.\eb
\section{Energy- and integral estimates}
\dipl{Lemma 3.1}
\begin{lemma}\label{lemma31}
  Let $f:\Sigma\times I\rightarrow M^n$ be a variation with normal velocity field 
  $\partial_tf=V$ and let $\phi\in\Gamma\big(T^{0,l}\Sigma\otimes N_f \big)$ with
  $\na_{\partial_t}\phi+\Delta^2\phi=Y.$ Then for arbitrary
  $\eta\in C^2(\Sigma\times I)$ with locally compact support
  \begin{eqnarray}
    \lefteqn{\hspace{-4em}\dfrac{d}{dt}\is\frac{1}{2}\lvert\phi\rvert^{2}\eta\,\mm
      +\is\langle\Delta\phi,\Delta(\eta\phi)\rangle\mm-
      \is\langle Y,\eta\phi\rangle\mm}\no\\
    ~~~~~~~~~~~&=&\is\eta\sum_{k=1}^{l}\langle 
    A(e_{i_{k}},e_{j}),V\rangle\langle \phi(e_{i_1},\ldots,e_{i_k},\ldots,e_{i_l}),
       \phi(e_{i_1},\ldots,e_{i_{k-1}},e_j,e_{i_{k+1}},\ldots,e_{i_l})\rangle\mm\no\\
    &&-\is\frac{1}{2}\lvert\phi\rvert^{2}\langle H,V\rangle\eta\,\mm
      +\is\frac{1}{2}\lvert\phi\rvert^{2}\partial_{t}\eta\,\mm.\label{31}
  \end{eqnarray}
\end{lemma}
Here $\w g(e_i,e_j)=\d_{ij}$ locally on $\Sigma\times I$ and summation over 
$j,i_\nu\in\{1,2 \},~1\le\nu\le l,$ is used.\\\\  
\bb 
With $\Phi_k(X,Y):=\eck{\phi(e_{i_1},\ldots,e_{i_{k-1}},X,e_{i_{k+1}},\ldots,e_{i_l}),
\phi(e_{i_1},\ldots,e_{i_{k-1}},Y,e_{i_{k+1}},\ldots,e_{i_l})}$ we have
\begin{eqnarray*}
\lefteqn{\hspace{-5em}\dfrac{d}{dt}\is\frac{1}{2}\eta\eck{\phi(e_{i_1},\ldots,e_{i_l}),
  \phi(e_{i_1},\ldots,e_{i_l})}\ds}\\
  &=&-\is\eck{\Delta^2\phi,\eta\phi} \ds+\is\eta\eck{ Y,\phi}\ds
   +\is\eta\sum_{k=1}^l\Phi_k(e_i, \na_{\partial_t}^{}e_i)\ds\\
  &&-\is\frac{1}{2}\eta\eck{ H,V}\lvert\phi\rvert^2\ds
   +\is\frac{1}{2}\partial_t\eta\lvert\phi\rvert^2\ds.
\end{eqnarray*}
Since $\tildeg(e_i,e_j)=\delta_{ij}$ we compute, using (\ref{157}),
$$\tildeg( \na_{\partial_t}^{}e_{i_k},e_j)
    +\tildeg(e_{i_k},\na_{\partial_t}^{}e_j)
  =-(\na_{\partial_t}^{}\tildeg)(e_{i_k},e_j)
  =2\eck{A(e_{i_k},e_j),V},$$ which implies $\na_{\partial_t}e_{i_k}
  =2\langle A(e_{i_k},e_j),V \rangle e_j
    -\tildeg(e_{i_k},\na_{\partial_t}e_j)e_j.$
From this we get
\begin{eqnarray*}
  \sum\limits_{k=1}^{l}\Phi_k\big(e_i,\na_{\partial_t}e_i \big)
  &=&\sum\limits_{k=1}^{l}\Big(2\langle A(e_i,e_j),V \rangle \Phi_k(e_i,e_j)
    -\Phi_k\big(\tildeg(e_i,\na_{\partial_t}e_j)e_j,e_i \big) \Big) \\
  &=&\sum\limits_{k=1}^{l}\Big(2\langle A(e_i,e_j),V \rangle \Phi_k(e_i,e_j)
    -\Phi_k\big(e_j,\tildeg(e_i,\na_{\partial_t}e_j)e_i \big) \Big) \\
  &=&\sum\limits_{k=1}^{l}\langle A(e_i,e_j),V \rangle \Phi_k(e_i,e_j).
\end{eqnarray*}
Since $\Delta$ is self-adjoint, the claim follows.\eb

\dipl{Lemma 3.3}
\begin{lemma}
  Under the assumptions of the previous Lemma let $\eta:=\g^s,$ where\\ 
  \mbox{$\g\in C^2(\Sigma\times I)$} is locally proper with
  $0\le\g\le 1$ and $s\geq 4.$ Then  for $c=c(n,s)$
  \begin{eqnarray}
    \lefteqn{\hspace{-2em}\dfrac{d}{dt}\is\bet{\phi}^2\g^s\mm
       +\is\bet{\na^2\phi}^2\g^s\mm
       -\is2\langle Y,\phi \rangle \g^s\mm}\label{32}\\
    &\le&\is\langle A\s \phi \s \phi ,V \rangle \g^s\mm
       +\is\bet{\phi}^2s\g^{s-1}\partial_t\g\mm
      +c\is\bet{\phi}^2\g^{s-4}\big(\bet{\na\g}^4
      +\g^2\bet{\na^2\g}^2 \big)\mm\no \\
    &&+c\is\bet{\phi}^2\big(\bet{\na A}^2+\bet{A}^4 \big)\g^s\mm
      +c\big(\norm{R}{2}{\infty}{}+\norm{DR}{4/3}{\infty}{} \big)
      \is\bet{\phi}^2\g^{s-2}\mm.\no
  \end{eqnarray}
\end{lemma}
\bb Because of \cref{lemma31} it suffices to estimate
\begin{eqnarray*}
  \is\bet{\na ^2\phi}^2\gamma^s\ds&\le& c\is\bet{\phi}^2\gamma^{s-4}\big(\bet{\na\g}^4
    +\g^2\bet{\na^2\g}^2 \big)\mm
    +c\is\bet{\phi}^2\big(\bet{\na A}^2+\bet{A}^4 \big)\g^s\mm \\
  &&+2\is\langle \Delta\phi,\Delta(\g^s\phi) \rangle \mm+c\big(\norm{R}{2}{\infty}{}
    +\norm{DR}{4/3}{\infty}{}\big)
    \is\bet{\phi}^2\g^{s-2}\mm.
\end{eqnarray*}
Since (\ref{32}) is scale invariant, we may assume (in the non-flat case) that
$\norm{R}{2}{\infty}{}+\norm{DR}{4/3}{\infty}{}=1.$
Differentiating
\begin{eqnarray}
  \na^2(\g^s\phi)&=&s(s-1)\g^{s-2}(\na\g)\otimes(\na\g)\otimes\phi\no\\ 
  &&+s\g^{s-1}(\na^2\g)\otimes\phi 
  +2s\g^{s-1}(\na\g)\otimes\na\phi+\g^s\na^2\phi\label{na2phi}
\end{eqnarray}
we get 
\begin{eqnarray}
  \is\lvert \na^2\phi\rvert^2\gamma^s\ds
  &\le&\is\eck{ \na^2\phi, \na^2(\gamma^s\phi)}\ds
    +c\is\lvert \na^2\phi\rvert\gamma^{s-1}\lvert \na\gamma\rvert \bet{\na\phi}\ds\no\\
  &&+c\is \bet{\na^2\phi}\gamma^{s-2}\big(\bet{\na\g}^2
    +\gamma \bet{\na^2\g}\big)\bet{\phi}\ds\no\\
  &\le&\is\eck{ \na^2\phi, \na^2(\gamma^s\phi)} \ds
     +\frac{1}{4}\is\bet{\na^2\phi}^2\gamma^s\ds
     +c\is\bet{\na\phi}^2\gamma^{s-2}\bet{\na\g}^2\ds\no\\
  &&+c\is\bet{\phi}^2\gamma^{s-4}\big(\bet{\na\g}^4
     +\gamma^2\bet{\na^2\g}^2\big)\ds,\label{33}
\end{eqnarray}
where we have used Young's inequality in the last estimate. 
Now we use partial integration for the third term
\begin{eqnarray}
  \lefteqn{\hspace{-2em}\is\bet{\na\phi}^2\g^{s-2}\bet{\na\g}^2\mm}\no \\
  ~~~~&=&\is\langle \phi,\na^*\big(\g^{s-2}\bet{\na\g}^2\na\phi \big) \rangle\mm\no \\
  &=&-\is\langle\phi,\na_{e_i}\big(\g^{s-2}\bet{\na\g}^2\big)\na_{e_i}\phi\rangle\mm
    -\is\langle \phi,\g^{s-2}\bet{\na\g}^2\Delta\phi \rangle \mm\no \\
  &\le&-\is\langle \phi,\Delta\phi \rangle \g^{s-2}\bet{\na\g}^2\mm
    +c\is\bet{\phi}\bet{\na\phi}\g^{s-3}\bet{\na\g}^3\mm
    +c\is\bet{\phi}\bet{\na\phi }\g^{s-2}\bet{\na\g}\bet{\na^2\g}\mm\no \\
  &\le&\e\is\bet{\na^2\phi}^2\g^s\mm+c(\e)\is\bet{\phi }^2\g^{s-4}\bet{\na\g}^4\mm\no\\
  &&+\frac{1}{2}\is\bet{\na\phi }^2\g^{s-2}\bet{\na\g}^2\mm
    +c\is\bet{\phi }^2\g^{s-2}\bet{\na^2\g}^2\mm,\label{34}
\end{eqnarray} 
where we have used  $\bet{\na_{e_i}(\g^{s-2}\bet{\na\g}^2)\na_{e_i}\phi}
  \le c\bet{\na\phi }\big(\g^{s-3}\bet{\na\g}^3+\g^{s-2}\bet{\na^2\g}\bet{\na\g}\big)$
  together with Young's inequality. Now choosing $\e$ small enough, adding up (\ref{33}) 
  and (\ref{34}), and absorbing yields
\begin{eqnarray}
  \lefteqn{\is\bet{\na^2\phi}^2\g^s\mm+\is\bet{\na\phi}^2\g^{s-2}\bet{\na\g}^2\mm}\no\\
  &\le&\dfrac{3}{2}\is\langle \na^2\phi ,\na^2(\g^s\phi) \rangle \mm
    +c\is\bet{\phi}^2\g^{s-4}\big(\bet{\na\g}^4+\g^2\bet{\na^2\g}^2\big)\mm.\label{35}
\end{eqnarray}

Again with partial integration, (\ref{23}) and 
$\na(\g^s\phi)=s\g^{s-1}(\na\g)\otimes\phi+\g^s\na\phi$ we compute
\begin{eqnarray}
  \lefteqn{\hspace{0em}\is\langle \na^2\phi ,\na^2(\g^s\phi ) \rangle \mm }\label{36} \\
  &=&\is\langle \na^*\na^2\phi ,\na(\g^s\phi ) \rangle \mm=\is\langle 
    -\Delta \na\phi ,\na(\g^s\phi )\rangle\mm\no \\
  &=&\is\langle -\na\Delta\phi ,\na(\g^s\phi ) \rangle \mm
    +\is\langle (\na\Delta-\Delta \na)\phi ,\na(\g^s\phi ) \rangle \mm\no
\end{eqnarray}
\begin{eqnarray}
  &\stackrel{(\ref{23})}\le&\is\langle \Delta\phi ,\Delta (\g^s\phi ) \rangle \mm
    +\is\langle A\s A\s \na\phi ,\na\phi  \rangle \g^s\mm
    +\is\langle A\s \na A\s \phi ,\na\phi  \rangle \g^s\mm\label{37} \\
  &&+c\is\bet{A}^2\bet{\phi}\bet{\na\phi}\g^{s-1}\bet{\na\g}\mm
    +s\is\bet{A}\bet{\na A}\bet{\phi}^2\g^{s-1}\bet{\na\g}\mm\label{38} \\
  &&+c\is(1+\bet{A})\bet{\phi}\big(\bet{\na\phi}\g^s
    +\bet{\phi}\bet{\na\g}\g^{s-1} \big)\mm\label{39}\\
  &&+c\is\bet{\na\phi}\big(\bet{\na\phi}\g^s
    +\bet{\phi}\bet{\na\g}\g^{s-1} \big)\mm\label{39neu},
\end{eqnarray}

where in (\ref{39}) and (\ref{39neu}) we estimated
\begin{eqnarray*}
  \bet{\Q00x\s \na\phi }
    &\le&c\bet{\Q00x}\bet{\na\phi}\le c\norm{R}{}{\infty}{}\bet{\na\phi}\\
    &\le&c\bet{\na\phi}
\end{eqnarray*}
and
\begin{eqnarray*}
  \bet{\Q01y\s \phi }&\le&c\bet{\Q01y}\bet{\phi}
    \le c\big(\norm{DR}{}{\infty}{}+\norm{R}{}{\infty}{}\bet{A} \big)\bet{\phi}\\
  &\le&c\big(1+\bet{A} \big)\bet{\phi}.
\end{eqnarray*}
The two integrals in (\ref{38}) can be estimated with Young
\begin{eqnarray}
  \is\bet{A}^2\bet{\phi}\bet{\na \phi}\g^{s-1}\bet{\na\g}\mm
  &\le&\e\is\bet{\na \phi}^2\g^{s-2}\bet{\na\g}^2\mm
    +c(\e)\is\bet{\phi}^2\bet{A}^4\g^s\mm,\label{311}\\
    \is\bet{A}\bet{\na A}\bet{\phi}^2\g^{s-1}\bet{\na\g}\mm
  &\le&\is\bet{\phi}^2\g^{s-2}\bet{\na\g}^2\bet{A}^2
    +\is\bet{\phi}^2\bet{\na A}^2\g^s\no \\
  &\le&\is\bet{\phi}^2\g^{s-4}\bet{\na\g}^4\mm
    +\is\bet{\phi}^2\big(\bet{\na A}^2+\bet{A}^4 \big)\g^s\mm.~~~~
\end{eqnarray}
With partial integration we estimate the second term in (\ref{37})
\begin{eqnarray*}
  \lefteqn{\hspace{-1em}\is\bet{A}^2\bet{\na \phi}^2\g^s\mm} \\
  &\le&-\is\bet{A}^2\langle \phi ,\Delta \phi  \rangle \g^s\mm
    +\is A\s \na  A\s \phi \s \na \phi\g^s\mm
    +c\is\bet{A}^2\bet{\phi}\bet{\na \phi}\g^{s-1}\bet{\na\g}\mm \\
  &\le&\e\is\bet{\na ^2\phi}^2\g^s\mm+c(\e)\is\bet{\phi}^2\bet{A}^4\g^s\mm\\
  &&+\is A\s\na A\s\phi\s\na\phi\g^s\mm+\e\is\bet{\na\phi}^2\g^{s-2}\bet{\na\g}^2\mm 
\end{eqnarray*}
and further
\begin{eqnarray*}
  \is A\s\na A\s\phi\s\na\phi\g^s\mm\le\dfrac{1}{2}\is\bet{A}^2\bet{\na\phi}^2\g^s\mm
    +c\is\bet{\phi}^2\bet{\na A}^2\g^s\mm.
\end{eqnarray*}

Adding up the last two estimates and absorbing yields
\begin{eqnarray}
  \lefteqn{\is\bet{A}^2\bet{\na \phi}^2\g^s\mm
    +\is A\s \na A\s \phi \s \na \phi \g^s\mm}\label{313}\\
  &\le&\e\is\bet{\na ^2\phi}^2\g^s\mm+\e\is\bet{\na \phi}^2\g^{s-2}\bet{\na\g}^2\mm
    +c(\e)\is\bet{\phi}^2\big(\bet{\na A}^2+\bet{A}^4 \big)\g^s\mm.\no
\end{eqnarray}

Combining (\ref{35}),(\ref{36}), and (\ref{311}) to (\ref{313}) yields
\begin{eqnarray*}
  \lefteqn{\is\bet{\na ^2\phi}^2\g^s\mm
    +\is\bet{\na \phi}^2\g^{s-2}\bet{\na\g}^2\mm } \\
  &\le&\dfrac{3}{2}\is\langle \Delta \phi ,\Delta (\g^s\phi ) \rangle \mm
    +\e\is\bet{\na ^2\phi}^2\g^s\mm
    +\e\is\bet{\na \phi}^2\g^{s-2}\bet{\na\g}^2\mm\\
  &&+c\is\bet{\phi}^2\g^{s-4}\big(\bet{\na\g}^4+\g^2\bet{\na^2\g}^2 \big)\mm
    +c(\e)\is\bet{\phi}^2\big(\bet{\na A}^2+\bet{A}^4 \big)\g^s\mm \\
  &&+\mbox{ the integrals in (\ref{39}) and (\ref{39neu}), coming from the
  estimated Q-terms.}
\end{eqnarray*}

So we only have to estimate (\ref{39}) and (\ref{39neu}) term by term
\begin{eqnarray*}
  \bullet~~(1+\bet{A})\bet{\phi}\bet{\na \phi}\g^s&\le&\bet{\phi}^2\g^s
      +\bet{\na \phi}^2\g^s+\bet{A}^2\bet{\phi}^2\g^s \\
    &\le&2\bet{\phi}^2\g^{s}+\bet{\na \phi}^2\g^s+\bet{A}^4\bet{\phi}^2\g^s,
      ~~~~~~~~~~~~~~~~~~~~~~~~~~~~~~~~~~~~~~~~~~~~~~~~
\end{eqnarray*}
where we further estimate $\bet{\na \phi}^2\g^s$ in the third point.
\begin{eqnarray*}
  \bullet~~(1+\bet{A})\bet{\phi}^2\bet{\na\g}\g^{s-1}
    &\le&(1+\bet{A}^2)\bet{\phi}^2\big(\g^s+\bet{\na\g}^2\g^{s-2} \big)\\
    &\le&2\bet{\phi}^2\g^s+2\bet{A}^4\bet{\phi}^2\g^s+2\bet{\phi}^2\bet{\na\g}^4\g^{s-4}.
      ~~~~~~~~~~~~~~~~~~~~~~~~~~~~~~~~~~~~
\end{eqnarray*}
From \cref{korollar62} with $p=2$ we get
\begin{eqnarray*}
  \bullet~~\is\bet{\na \phi}^2\g^s\mm&\le&\e\is\bet{\na ^2\phi}^2\g^s\mm
  +c(\e)\is\bet{\phi}^2\g^{s-2}\mm.~~~~~~~~~~~~~~~~~~~~~~~~~~~~~~~~~~~~~~~~~~~~~~~~~~~~
\end{eqnarray*}
\begin{eqnarray*}
  \bullet~~\bet{\na \phi}\bet{\phi}\bet{\na\g}\g^{s-1}&\le&\bet{\na \phi}^2\g^s
    +\bet{\phi}^2\bet{\na\g}^2\g^{s-2} \\
  &\le&\bet{\na \phi}^2\g^s+\bet{\phi}^2\g^s+\bet{\phi}^2\bet{\na\g}^4\g^{s-4}.
     ~~~~~~~~~~~~~~~~~~~~~~~~~~~~~~~~~~~~~~~~~~~~~~~~~~
\end{eqnarray*}
Choosing $\e>0$ appropriately the claim follows.\eb
For $\w\g\in C_c^\infty(M)$ with $0\le\w\g\le 1$ let $\g:=\w\g\circ f.$ This implies
\begin{align}
  \bet{\na\g}&=\bet{D\w\g\cf\c Df}\le c\norm{D\w\g}{}{\infty}{}
  &\5\mbox{and}\qquad
  \bet{\na^2\g}&=\bet{\na(D\w\g\cf\c Df)}\label{314}\\
  &&&=\bet{D^2\w\g\cf(Df\,\c\,,Df\,\c\,)+D\w\g\cf\c A}\no \\
  &&&\le c\big(\norm{D^2\w\g}{}{\infty}{}+\norm{D\w\g}{}{\infty}{}\bet{A} \big).\no
\end{align}
Furthermore, we specialize to the Willmore flow, i.e.
\begin{eqnarray}
  V=\partial_tf=-\Delta H+\pe03+\Q010.\label{315}
\end{eqnarray}
For the time derivative we thus get
\begin{eqnarray}
  \partial_t\g&=&D\w\g\cf\c\partial_tf\no \\
  &=&D\w\g\cf\c\big(-\Delta H+\pe03+\Q010 \big).\label{316}
\end{eqnarray} 

\dipl{Lemma 3.4}
\begin{lemma}\label{lemma34}
  Let $f:\Sigma^2\times I\rightarrow M^n$ be a smooth locally proper Willmore flow. Then
  for 
  $\phi:=\na^mA$ with $m\in\N_0,$ $\g=\w \g\cf$ satisfying (\ref{314}) and $s\geq2m+4$
  we have
  \begin{eqnarray}
    \lefteqn{\hspace{-0em}\dfrac{d}{dt}\is\bet{\phi}^2\g^s\mm
      +\dfrac{3}{4}\is\bet{\na^2\phi}^2\g^s\mm}\label{hiergibtsQ}\\
    &\le&\is\big(\pe{m+2}3+\pe m5+\Q{m+2}1{m+1}+\Q m1R\,\big)\s\phi\,\g^s\mm\no\\
    &&+\is\big<\na^m\Delta\Q 01y+\Delta\Q m1m+\na\Q {m+1}1{m+1}
      ,\phi \big>\,\g^s\mm\label{x317}
      +c\,C_{scal}^{2m+4}\iz\bet{A}^2\g^{s-4-2m}\mm,\no
  \end{eqnarray}
where $C_{scal}=\sum\limits_{i=1}^{2}\norm{D^i\w\g}{1/i}{\infty}{}
  +\sum\limits_{i=0}^{1}\norm{D^{i\!}R}{\frac{1}{i+2}}{\infty}{}$ and $c=c(n,s).$
\end{lemma}

\bb After scaling we may assume that
$C_{scal}=1.$
We estimate the terms in (\ref{32}). From (\ref{26}) we know that 
$Y=\pe{m+2}3+\pe m5+\na^m\Delta\Q 01y+\Delta\Q m1m+\na\Q {m+1}1{m+1}+\Q{m+2}1{m+1}
  +\Q m1R,$ 
and substituting (\ref{315}) we obtain
\begin{eqnarray*}
  \lefteqn{\hspace{-3em}\is2\langle Y,\phi\rangle\g^s\mm
    +\is\langle A\s\phi\s\phi,V \rangle\g^s\mm
    +c\is\bet{\phi}^2\big(\bet{\na A}^2+\bet{A}^4 \big)\g^s\mm} \\
  &=&\is\big(\pe{m+2}3+\pe m5+\Q{m+2}1{m+1}+\Q m1R\big)\s\phi\,\g^s\mm\\
  &&+\is\big<\na^m\Delta\Q 01y+\Delta\Q m1m+\na\Q {m+1}1{m+1}
      ,\phi \big>\,\g^s\mm,
\end{eqnarray*}
where we used that 
$\langle A\s\phi\s\phi,\Q 010\rangle=\Q m30\s\phi=\Q {m+2}1{m+1}\s\phi.$
So it suffices to estimate each of
\begin{eqnarray}
  \is\bet{\phi}^2\g^{s-1}\partial_t\g\mm,~~\is\bet{\phi}^2\g^{s-4}\big(\bet{\na\g}^4
  +\g^2\bet{\na^2\g}^2 \big)\mm,
  ~~\is\bet{\phi}^2\g^{s-2}\mm\label{318}
\end{eqnarray}
by
\begin{eqnarray*}
  \is\big(\pe{m+2}3+\pe m5 \big)\s \phi\,\g^s\mm+\e\is\bet{\na ^2\phi}\g^s\mm
    +c(\e)\iz\bet{A}^2\g^{s-4-2m}\mm.
\end{eqnarray*}
At first, we get using \cref{korollar63} \info{Korollar 6.3} with 
$\phi:=A,~k:=m,~p=2$ and substituting $s$ by $s-4$ 
\begin{eqnarray*}
  \is\bet{\phi}^2\g^{s-4}\mm\le\is\bet{\na\phi}^2\g^{s-2}\mm
    +c\iz\bet{A}^2\g^{s-4-2m}\mm.
\end{eqnarray*}
Analogously, it follows from the same corollary with $\phi :=A,~k:=m+1,~p:=2$
and substituting $s$ by $s-2$
\begin{eqnarray}
  \is\bet{\na\phi}^2\g^{s-2}\mm\le\e\is\bet{\na^2\phi}^2\g^s\mm
    +c(\e)\iz\bet{A}^2\g^{s-4-2m}\mm.\label{318a}
\end{eqnarray}

Combining the last two estimates yields
\begin{eqnarray}
  \is\bet{\phi}^2\g^{s-4}\mm+\is\bet{\na \phi}^2\g^{s-2}\mm
    \le\e\is\bet{\na ^2\phi}^2\g^s\mm
    +c(\e)\iz\bet{A}^2\g^{s-4-2m}\mm,\label{319}
\end{eqnarray}
which treats the last term in (\ref{318}). From the evolution equation (\ref{316}) we
have
\begin{eqnarray}
  \is\bet{\phi}^2\g^{s-1}\partial_t\g\mm
    =\is\bet{\phi}^2\g^{s-1}D\w\g\cf\c\big(-\Delta H+\pe03+\Q010\big)\mm.\label{320}
\end{eqnarray}
We estimate the second summand on the right-hand side of (\ref{320}) with
Young\\\mbox{$(p=4,~q=4/3)$}
\begin{eqnarray}
  \is\bet{\phi}^2\g^{s-1}D\w\g\cf\c\pe03\mm
  &\le&c\is\bet{\phi}^2\g^{\frac{s}{4}-1}\bet{A}^3\g^{\frac{3}{4}s}\mm\no \\
  &\le&c\is\bet{\phi}^2\bet{A}^4\g^s\mm+c\is\bet{\phi}^2\g^{s-4}\mm\label{321} \\
  &\le& c\is\pe m5\s \phi \g^s\mm+\e\is\bet{\na^2\phi}^2\g^s\mm\no\\
  &&+c(\e)\iz\bet{A}^2\g^{s-4-2m}\mm,\no
\end{eqnarray}

where we used the interpolation inequality (\ref{319}) in the last step. 
The last summand in (\ref{320}) is estimated as follows:
\begin{eqnarray*}
  \is\bet{\phi}^2\g^{s-1}D\w\g\cf\c\Q010\mm
  &\le&c\is\bet{\phi}^2\bet{A}\g^{s-1}\mm\\
  &\le&c\is\bet{\phi}^2\g^{s-2}\mm+c\is\bet{\phi}^2\bet{A}^4\g^s\mm.
\end{eqnarray*}

The first summand on the right-hand side above can again be treated by (\ref{319}). 
For the second, note that
$\bet{\phi}^2\bet{A}^4$ can be written as $\pe m5\s \phi.$
We estimate the remaining first summand in (\ref{320}) as in
\cite{KS02}. Namely, let
$\{e_i\}_{i=1,2}$ be an adapted local orthonormal frame. Then we pointwise have, using
Einstein's convention
\begin{eqnarray*}
  -\na^*\big((D\w\g)\cf\c\na H \big)&=&\p_{e_i}\big((D\w\g)\cf\c \na_{e_i}H \big)\\
  &=&(D^2\w\g)\cf( Df\!\cdot\! e_i,\na_{e_i}H)+(D\w\g)\cf\c\Delta H,
\end{eqnarray*}
\begin{eqnarray}
  \Rightarrow~~(D\w\g)\cf\c\Delta H&=&-\na^*\big((D\w\g)\cf\c \na H \big)
    -(D^2\w\g)\cf( Df\!\cdot\! e_i,\na_{e_i}H).\label{322}
\end{eqnarray}
Since
\begin{eqnarray*}
  \is\bet{\phi}^2\g^{s-1}\na^*(D\w\g\cf\c\na H)\mm
  &=&\is {\w g}^{\,*}\big(\na(\bet{\phi}^2\g^{s-1}),D\w\g\cf\c\na H \big)\mm\\
  &=&\is\na_{e_i}(\bet{\phi}^2\g^{s-1})D\w\g\cf\c\na_{e_i}H\mm
\end{eqnarray*}
we get using partial integration
\begin{eqnarray*}
  \lefteqn{-\is\bet{\phi}^2\g^{s-1}D\w\g\cf\c\Delta H\mm } \\
  &=&\is\na_{e_i}(\bet{\phi}^2\g^{s-1})D\w\g\cf\c \na_{e_i}H\mm
    +\is\bet{\phi}^2\g^{s-1}D^2\w\g\cf( Df\!\cdot\! e_i,\na_{e_i}H)\mm \\
  &\le& c\is\big(\bet{\phi}\bet{\na\phi}\g^{s-1}
    +\bet{\phi}^2\g^{s-2}\bet{\na\g} \big)\bet{\na A}\bet{D\w\g\cf}\mm
    +c\is\bet{\phi}^2\g^{s-1}\bet{D^2\w\g\cf}\bet{\na A}\mm \\
  &\stackrel{(\ref{314})}{\le}& c\is\bet{\phi}\bet{\na\phi}\bet{\na A}\g^{s-1}\mm
    +c\is\bet{\phi}^2\bet{\na A}\g^{s-2}\mm+c\is\bet{\phi}^2\bet{\na A}\g^{s-1}\mm
\end{eqnarray*}
and further with Young
\begin{eqnarray*}
  &\le& c\is\bet{\na \phi}^2\g^{s-2}\mm+c\is\bet{\phi}^2\bet{\na A}^2\g^s\mm
    +c\is\bet{\phi}^2\g^{s-4}\mm.\\
  &\stackrel{(\ref{319})}{\le}&\is\pe{m+2}3\s \phi \g^s\mm
    +\e\is\bet{\na ^2\phi}^2\g^s\mm
    +c(\e)\iz\bet{A}^2\g^{s-4-2m}\mm.
\end{eqnarray*}
Altogether, we have estimated (\ref{320}), remaining the second term in (\ref{318}). But
for this term we estimate
\begin{eqnarray*}
  \lefteqn{\hspace{-7em}\is\bet{\phi}^2\g^{s-4}\big(\bet{\na\g}^4
    +\g^2\bet{\na^2\g}^2\big)\mm} \\
  &\stackrel{(\ref{314})}{\le}& c\is\bet{\phi}^2\g^{s-4}\mm
    +c\is\bet{\phi}^2\big(1+\bet{A}^2 \big)\g^{s-2}\mm \\
  ~~~~~~~~~~~~~~~~~~~~~~~~&\le& c\is\bet{\phi}^2\g^{s-4}\mm
    +c\is\bet{\phi}^2\bet{A}^4\g^s\mm \\
  &\stackrel{(\ref{319})}{\le}&\e\is\bet{\na ^2\phi}^2\g^s\mm
    +c(\e)\iz\bet{A}^2\g^{s-4-2m}\mm
    +c\is\pe m5\s \phi \g^s\mm,
\end{eqnarray*}
completing the proof.\eb
\section{Sobolev inequalities for Riemannian Manifolds}
\dipl{Theorem 4.1}\notiz{vl noch ein beispiel dass A statt H stehen muss...}
\begin{theorem}{\rm (Michael-Simon Sobolev inequality for Riemannian
    manifolds I).}\label{mssi}~\\
  Let $f:(\Sigma^d,\tildeg\,)\to(M^n,g)$ be an isometric $C^2$-immersion where
  $(\Sigma,\tildeg\,)$ and $(M,g)$ are (open or closed) complete Riemannian
  manifolds of dimension $d\geq2$ and $n=d+m$ respectively. If
  $\Lambda:=\norm{ricci_{(M,g)}}{}{\infty}{}~^{\hspace{-1.1em}\nf 12}\hspace{0.0em}
  +\inj(M,g)^{-1}$ is bounded, then for any
  $u\in C^1_c(\Sigma)$ we have
  \begin{equation}
    \bigg(\is \bet{u}^{\frac d{d-1}} d\mu \bigg)^{\frac{d-1}{d}}\le
    c\is\bet{\na u}d\mu+c\is\bet{A_f}\bet{u}d\mu+c\Lambda\!\!
    \is\bet ud\mu,\label{mssieqn} 
  \end{equation}
  where $c=c(d,n)$ is a universal constant.
\end{theorem}
\begin{theorem}{\rm (Michael-Simon Sobolev inequality for Riemannian  
  manifolds II).}\label{mssic}~\\
  Let $f:(\Sigma^d,\tildeg\,)\to(M^n,g)$ be an isometric $C^2$-immersion of a
  complete Riemannian manifold of dimension $d\geq 2$ into a closed Riemannian manifold
  $(M,g)$ of  dimension $n=d+m.$ Assume that there exists an isometric immersion
  $I:(M,g)\to\R^N$ with the property that 
  $\Vert A_I\Vert_{L^\infty(M,g)}<\infty.$ Then for any $u\in C^1_c(\Sigma)$ we have
  \begin{equation}
    \bigg(\is \bet{u}^{\frac d{d-1}} d\mu \bigg)^{\frac{d-1}{d}}\le
    c\is\bet{\na u}d\mu+c\is\bet{H_f}\bet{u}d\mu
      +c\Vert A_I\Vert_\infty\is\bet ud\mu,\label{mssieqnc} 
  \end{equation}
  where $c=c(d,n).$
\end{theorem}  
\bf{Remark:}\rm~By Nash's embedding theorem the existence of such an immersion is clearly
automatic in case $M$ is compact. To see that (\ref{mssieqn}) and (\ref{mssieqnc}) cannot
hold in general without an extra term on the right-hand side, just take the standard
embedding of the Sphere
$\mathbb S^d\subset\mathbb S^n\subset(\R^{n+1},\d_{eucl}),$ where
$\mathbb S^d=\mathbb S^n\cap E$ and $E=\R e_1\oplus\ldots\oplus\R e_{d+1}$ is considered
as a linear subspace of $\R^{n+1}.$ Note that $\mathbb S^d$ is a totally geodesic
submanifold of $\mathbb S^n$ and thus the second fundamental form vanishes
identically.\notiz{achtung *****************************}
If we then define $u:\mathbb S^d\to\R$ by $u\equiv1,$ (\ref{mssieqn}) and 
(\ref{mssieqnc}) would obviously imply that $\mathbb S^d$ had zero volume with
respect to the induced metric.\er\\
\bf{Proof of \cref{mssic}:}\rm~Let
$I:(M,g)\hookrightarrow(\R^N,\d_{eucl})$ be an immersion as in the statement of
the Theorem. Our starting point is the euclidean inequality
\begin{equation*}
    \bigg(\is \bet{u}^{\frac d{d-1}} d\mu \bigg)^{\frac{d-1}{d}}\le
    c\is\bet{\na u}d\mu+c\is\bet{H_{I\cf}}\bet{u}d\mu
\end{equation*}
(see proof of \cref{mssi}). The claim follows immediately by substitution of
(\ref{meancurvabsch}).\eb\\
\bf{Proof of \cref{mssi}:}\rm~Assume 
that $f:(\Sigma,\bar g)\hookrightarrow(\R^n,\d_{eucl})$ is an isometric
embedding and we have given $h\in C^1_c(U)$ for $U\subset\R^n$ open and $h\geq 0.$
Applying \cite{MS}, Theorem 18.6, we get
\begin{eqnarray*}
  \Big(\is h^pd\mu_{\bar g} \Big)^{1/p}\le c\is \big(\bet{{}^\Sigma\na h}
    +\bet{H_f}h \big)d\mu_{\bar g},
\end{eqnarray*}
where $c=c(d)$ and $p:=\frac d{d-1}.$ Now if $\w u\in C_c^1(\Sigma),$ we extend $\w u$
from supp$\,\w u\subset\!\subset\Sigma,$ e.g. using slice coordinates, to an appropriate
open subset $U\subset\R^n$ such that the extension $u$ has compact support in $U$ and
obtain
\begin{eqnarray}
   \Big(\is \bet{u}^pd\mu_{\bar g} \Big)^{1/p}\le c\is \big(\bet{{}^\Sigma\na u}
    +\bet{H_f}\bet{u} \big)d\mu_{\bar g},\label{msvorl}
\end{eqnarray}
since $\bet{u}$ is at least Lipschitz in $U.$ In case f is merely an isometric
immersion, one can choose an isometric embedding $E:(\Sigma,\bar g)\to\R^N,$ applying
(\ref{msvorl}) to $f_\e:=(f,\e E):(\Sigma,(1+\e^2)\bar g)\to\R^{N+n}$ and then letting 
$\e\searrow0$ (see \cite{KS02}).\\\\
We now want to prove the case $(M,g)\neq(\R^n,\d_{eucl}).$ Define
$r_0:=c(n)\Lambda^{-1},$
where $\Lambda,c$ are as in \cref{harmcoord1}. Let $\{B_{r_0(p_i)}
\}_{i\in\N}$
be a uniformly locally finite covering of M and $\{\widetilde\eta_i \}_{i\in\N}$ a
partition of unity subordinate to $\{B_{r_0(p_i)}\}_{i\in\N}$ as in \cref{456}. 
Further, let $\{\psi_i:B_{r_0}(p_i)\rightarrow V_i\subset\R^n\}_{i\in\N},$
$\psi_i=\{y^\a_i \}_{1\le\a\le n}$ be a countable atlas of harmonic coordinates as in
\cref{harmcoord1}, and $\{x^j\}_{1\le j\le d}$ arbitrary coordinates on $\Sigma.$ 
Since we want to localize in the target, we may fix some arbitrary $i\in\N$ and
omit it almost every time. Note that we have
\begin{eqnarray}
  \frac 1c\d\le \big(g_{\a\b} \big)\le c\d\label{maequiv}
\end{eqnarray}
and
\begin{eqnarray}
  \bet{{}^M\Gamma_{\b\g}^\a}\le c\Lambda.\label{christoffel}
\end{eqnarray}
We want to introduce some notation. Let
$\Sigma_{r_0,i}:=f^{-1}\big(B_{r_0}(p_i) \big)\subset\Sigma.$ On $V\subset\R^n$ we
consider the Riemannian metrics $g$ and $\d,$ where $\d$ is the standard metric and,
by slight abuse of notation, g stands for the coordinate representation 
$(g_{\a\b})_{1\le\a,\b\le n}$  of g with respect to the harmonic chart $\psi.$ Also,
define isometric immersions
$$\bar f:=\big(\Sigma_{r_0},\bar g \big)\to(V,\d)$$
and
$$\w f:=\big(\Sigma_{r_0},\tildeg\,\big)\to(V,g)$$
such that $\bar f=\w f=\psi\circ f\strichklein{$\Sigma_{r_0}$}\;\;\,$ as maps between
manifolds. From $c^{-1} \d\le g\le c\d$ it follows that
\begin{eqnarray}
  \frac 1c\bar g\le\tildeg\le c\bar g\label{sigmaaequiv}
\end{eqnarray}
by properties of the pullback metric. From this, we get an estimate for the Gram
determinants
\begin{eqnarray}
  \frac 1{c^d}\det(\bar g)\le\det(\tildeg)\le c^d\det(\bar g)\label{detbil}
\end{eqnarray}
as follows: Since this is a pointwise estimate, we can choose for any $x\in\Sigma_{r_0}$
a local chart $\{x^k\}_{k=1,2}$ around $x$ with $\bar g_{ij}(x)=\d_{ij}$ and  
obtain from $c^{-1}\d\le\big(\widetilde g_{ij}(x)\big)\le c\d$ the desired
equivalence. We remark that although the determinant is not defined for
bilinear forms, inequality (\ref{detbil}) does make sense by invariance with respect to
coordinate transformations.\\


We now want to apply the flat-case inequality (\ref{msvorl}) locally to 
$\Sigma^*:=\Sigma_{r_0}$ and $u^*:=\eta u\in C^1_c(\Sigma^*),$ where we define
$\eta:=\w{\eta}\cf.$
If we denote by $\bar\mu,\w\mu,\barGamma,\w\Gamma,...$ etc. the 
locally defined geometric quantities on $\big(\Sigma_{r_0},\bar g \big)$ or 
$\big(\Sigma_{r_0},\tildeg\, \big)$ respectively, we get for any
$u\in C^1_c(\Sigma)$ using dominated convergence
\begin{eqnarray}
  \norm{u}{}{L^p(\w\mu)}{}&\le&\sum^\infty_{i=1}\norm{\eta_iu}{}{L^p(\w\mu)}{}\no \\
  &=&\sum\limits^\infty_{i=1}\bigg(\int\limits_{\Sigma_{r_0,i}}\bet{\eta_iu}^p
    d\w\mu_i \bigg)^{1/p}\no \\
  &\stackrel{(\ref{detbil})}{\le}&
    c \sum\limits^\infty_{i=1}\bigg(\int\limits_{\Sigma_{r_0,i}}\bet{\eta_iu}^p d\bar
    \mu_i\bigg)^{1/p}\no \\
  &\le& c\sum\limits^{\infty}_{i=1}\bigg(\int\limits_{\Sigma_{r_0,i}}
    \lvert\nabla(\eta_iu )\rvert_{\bar g_i}d\bar\mu_i
      +\int\limits_{\Sigma_{r_0,i}}\bet{H_{\bar f}}_{\d}
      \bet{u}\eta_id\bar\mu_i \bigg).\label{47a}
\end{eqnarray}

To estimate the first integral in (\ref{47a}), we compute (for fixed i) 
$$\bet{\na(\eta u)}_{\bar g}\le c\bet{\na(\eta u)}_{\tildeg}
  \le c\left(\bet u\bet{\na\eta}_{\tildeg}+\eta\bet{\na u}_{\tildeg} \right)$$
and further for $N=N(n)$ as in \cref{456}
$$\bet{\na\eta}_{\tildeg}=\bet{\na(\widetilde\eta\circ f)}_{\tildeg}
=\bet{D\widetilde\eta\cf\cdot Df}_{\tildeg}
  \le c\bet{D\widetilde\eta}_{g}\cf\bet{Df}_{\w f}\le c(N)r_0^{-1}= c\Lambda$$
by construction of $\widetilde\eta$ in \cref{456}. Combining, we get with 
(\ref{detbil}), i.e. $d\bar\mu\le c^{d/2}\,d\w\mu,$ and dominated convergence
\begin{eqnarray*}
  \sum\limits_{i=1}^{\infty}\int\limits_{\Sigma_{r_0,i}}
    \bet{\na(\eta_iu)}_{\bar g_i}d\bar\mu_i
  &\le&c\Lambda\sum\limits^{\infty}_{i=1}
    \int\limits_{\Sigma_{r_0,i}}\bet ud\w\mu_i+c\is\bet{\na u}_{\tildeg}d\w\mu \\
  &\le& c\Lambda N(n)\is\bet ud\w\mu+c\is\bet{\na u}_{\tildeg}d\w\mu.
\end{eqnarray*}
The last inequality holds, since $\w\mu_i=\w\mu\llcorner\chiunten{\Sigma_{r_0,i}}$ and
for any $x\in\Sigma$ we have\\ 
  $\sum\limits_{i=1}^\infty\chi_{\Sigma_{r_0,i}}(x)
  \le\sum\limits_{i=1}^\infty\chi_{B_{r_0}(p_i)}\circ f(x)\le N(n)$ by construction of 
the $\{\widetilde\eta_i\}$ in \cref{456}.
To estimate the second integral in (\ref{47a}), we have to estimate
\begin{eqnarray}
  \bet{\bar H}_\d^2\le c\bet{\w A}^2_{\w f}+c\Lambda^2.\label{haveto}
\end{eqnarray}
on $\Sigma_{r_0}.$
Since this is a pointwise estimate, we now fix for arbitrary 
$x_0\in\Sigma_{r_0}$ Riemannian normal coordinates $\{x^j \}_{j=1,2}$ with respect to 
$\tildeg$ around $x_0$. We want to remark that the second fundamental form $A$ is a
natural map in the following sense: For an isometry $\varphi$ it is easy to see that 
$A^\a_{\varphi^*f}=A^\a_{f\circ\varphi}=\varphi^*(A^\a_f)$ and thus 
$\bet{A_{\varphi^*f}}_{\varphi^*f}=\bet{A_f}_f
{\circ\varphi}.$ Therefore
(\ref{haveto})
is not affected when choosing coordinates (clearly, the geometry of f should not depend
on the parametrization).
For the next calculation we want to define the Matrix valued functions
\begin{eqnarray*}
  F_{jk}^\a:=\p^2_{jk}f^\a,~~G_{jk}^\a:={}^\Sigma\Gamma_{jk}^l\p_l f^\a\mbox{~~and~~}
  C_{jk}^\a:={}^M\Gamma_{\b\g}^\a\cf\p_j f^\b\p_k f^\g.
\end{eqnarray*} 
Then the Gauss formula reads $A=F-G+C$ in general, and in our case (in $x_0$) with the
obvious notation $\bar A=F-\barG\mbox{~~and~~}\w A=F+C.$ Contracting, this implies for
the mean curvature of $\bar f$
\begin{eqnarray*}
  \bar H=\bar g^{jk}F_{jk}-\bar g^{jk}\barG_{jk}\mbox{~~and~~}
    \w A=F+C.
\end{eqnarray*}
Because $c^{-1}\d\le \big(\bar g_{jk}\big)\le c\d$ we get using polarization
$\bet{\bar g_{jk}},\bet{\bar g^{jk}}\le c.$ From this and by orthogonality of
$\bar H\perp_\d \bar g^{jk}\barG_{jk}$ we can now estimate
\begin{eqnarray*}
  \bet{\bar H}_\d^2&=&\bet{\bar g^{jk}F_{jk}}_\d^2-\bet{\bar g^{jk}\barG_{jk}}^2_\d
  \le c\bet{F}_\d^2
  \le c\bet{\w A}_\d^2+c\bet{C}_\d^2.
\end{eqnarray*}
Taking $$\bet{\w A}^2_\d=\sum_{j,k=1}^2\bet{\w A_{jk}}^2_\d
\le c\sum_{j,k=1}^2\bet{\w A_{jk}}^2_g=c\bet{\w A}^2_{\w f}(x_0)$$
into account, the Theorem is now proved since we
can estimate
$$\bet{C}^2_\d=\sum_{jk\a}\bet{{}^M\Gamma_{\b\g}^\a\cf\p_j f^\b\p_k f^\g}^2
\stackrel{(\ref{christoffel})}{\le} c\Lambda^2\bet{\p f}^4_\d\le
c\Lambda^2\bet{Df}^4_{\w f}(x_0)\le c\Lambda^2.$$\eb

\section{Interpolation of lower-order terms}

\dipl{Lemma 4.2}
\begin{lemma}\label{lemma42}
  Let $(M^n,g)$ be a Riemannian manifold, 
  $\norm{ricci_{(M,g)}}{1/2}{\infty}{}+\inj(M,g)^{-1}=:\Lambda<\infty$ and
  $f:(\Sigma^2,\w g)\rightarrow (M,g)$ an isometric $C^4-$immersion. Then for any 
  $\g\in C^1_c(\Sigma)$ with $0\le\g\le 1$ we have
  \begin{eqnarray}
    \lefteqn{\hspace{-2em}\is\bet A^6\g^4\mm+\is\bet A^2
      \bet{\na A}^2\g^4\mm}\label{413}\\
    &\le& c\iz\bet A^2\mm\is\left(\bet{\na ^2\!A}^2\g^4+\bet{A}^6\g^4 \right)\mm
      +c\big(\norm{\na\g}{4}{\infty}{}
      +\Lambda^4\big)\bigg(\iz\bet A^2\mm \bigg)^2,\no
  \end{eqnarray}
where c=c(n).
\end{lemma}
\bb
Again, since the above inequality is scale invariant, we may assume that 
$\norm{\na\g}{4}{\infty}{}+\Lambda^4=1.$
We want to apply \cref{mssi} to $u:=\bet A\bet{\na A}\g^2$. If $u$ should not be 
differentiable, we take $u_\e:=\sqrt{\bet A^2\bet{\na A}^2+\e^2}\g^2$ instead, so
that since 
$\bet{\na u_\e}\le\bet{\na A}^2\g^2+\bet A\bet{\na^2\!A}\g^2+c\bet A\bet{\na
A}\g+c\e\g$
we get from (\ref{mssieqn}) after letting $\e\searrow 0$ 
\begin{eqnarray*}
  \lefteqn{\is\bet A^2\bet{\na A}^2\g^4d\mu}~~~~~~~~~~~\\
  &\le& c\bigg( \is\bet A\bet{\na ^2\!A}\g^2d\mu+\is\bet{\na A}^2\g^2d\mu
    +\is\bet A\bet{\na A}\g d\mu\\
  &&+\is\bet A^2\bet{\na A}\g^2d\mu+\is\bet A\bet{\na A}\g^2d\mu\bigg)^2\\
\end{eqnarray*}
and further with Young's inequality and Cauchy-Schwarz
\begin{eqnarray}\label{415}
  \le c\iz\bet A^2\mm\is(\bet{\na ^2\!A}^2\g^4+\bet A^6\g^4)\mm
    +c\bigg( \iz\bet A^2\mm \bigg)^2+c\bigg( \is\bet{\na A}^2\g^2\mm \bigg)^2
\end{eqnarray}
since
\begin{eqnarray}\label{416}
  \bigg( \is\bet A^4\g^2\mm \bigg)^2\le\iz\bet A^2\mm\is\bet A^6\g^4\mm.
\end{eqnarray}

Using partial integration for the last term we get
\begin{eqnarray}\label{417}
  \is\bet{\na A}^2\g^2\mm&=&\is\langle \na ^*(\g^2\na A),A\rangle\mm\nonumber\\
  &\le&c\is\bet A\bet{\na ^2 A}\g^2\mm+c\is\bet A\bet{\na A}\g\mm\nonumber\\
  &\le&c\bigg( \is\bet A^2\mm\is\bet{\na ^2\!A}^2\g^4\mm \bigg)^{1/2}
    +c\iz\bet A^2\mm+\frac{1}{2}\is\bet{\na A}^2\g^2\mm.~~~~~~~~
\end{eqnarray}\notiz{approxargument muss weg - siehe unten}
Absorbing and substituting in (\ref{415}) we have an estimate for the second summand in 
(\ref{413}). With the same approximation argument as above, apply \cref{mssi} to 
  $u_\e:=\sqrt{\bet A^6+\e^2 }\g^2\le (\bet A^3+\e)\g^2$
and on account of $\bet{u_\e}
\le 3\bet A^2\bet{\na A}\g^2+2\bet A^3\g\bet{\na\g}+2\e\g\bet{\na\g}$
we get after letting $\e\searrow 0$
\begin{eqnarray}
  \is\bet A^6\g^4\mm&\le&c\bigg(\is\bet{\na A}^2\g^2 \mm\bigg)^2
    +c\bigg(\is\bet{A}^4\g^2\mm \bigg)^2+c\bigg(\iz\bet A^2\mm  \bigg)^2.
\end{eqnarray}
Combining with (\ref{416}) and (\ref{417}), eventually yields the estimate for the first 
summand in (\ref{413}). The Lemma then follows after rescaling.\eb\\
We define for an arbitrary section $\phi\in\Gamma^0\big(T^{0,l}\Sigma\otimes N_f
\big)$
and any measurable $U\subset\Sigma$
\begin{eqnarray*}
  \norm{\phi }{}{p,U}{}:=\bigg(\int_U\bet{\phi}^p\mm \bigg)^{\frac{1}{p}}
  \qquad\mbox{and}\qquad\norm{\phi }{}{\infty,U}{}:=
  \sup_{U}\bet{\phi}.
\end{eqnarray*}
We abbreviate $\norm{\phi}{}{\infty,\Sigma}{}=:\norm{\phi}{}{\infty}{}$ and when there is
no ambiguity we also write
\begin{eqnarray}
  \norm{S }{}{\infty}{}:=\sup_{M}\bet{S}
\end{eqnarray}
for sections $S\in\Gamma^0\big(\*(T^{k,l\!}M) \big).$

\dipl{Lemma 4.3}
\begin{lemma}\label{lemma43}
  Let $(M^n,g)$ be a Riemannian manifold with the property that
  \mbox{$\norm{ricci_{(M,g)}}{1/2}{\infty}{}$} $+\inj(M,g)^{-1}<\infty$,
  and $f:(\Sigma^2,\tildeg)\rightarrow (M,g)$ be a proper isometric immersion.
  Then for arbitrary $\phi \in\Gamma\big(T^{0,l}\Sigma\otimes N_f \big)$ and
  $\g=\w\g\cf$ as in (\ref{314}) we have
  \begin{eqnarray}
  \norm{\phi}{4}{\infty}{1}\!\!\!\!&\le&\!\!\!\!c\norm{\phi}{2}{2}{0}
        \left(\norm{\na^2\phi}{2}{2}{0}
      +\norm{\bet{\phi}^2\bet{A}^4}{}{1}{0}
      +c_{scal}^4\norm{\phi}{2}{2}{0} \right).\label{421}
  \end{eqnarray}
  Moreover, for $\phi:=A$ and provided that
  \begin{eqnarray}\label{422}
    \norm{A}{2}{2}{0}\le\e_0
  \end{eqnarray}
  for $\e_0$ small enough depending only on the dimension n we get
  \begin{eqnarray}
    \norm{A}{4}{\infty}{1}\le c\norm{A}{2}{2}{0}\left(\norm{\na^2\!A}{2}{2}{0}
    +c_{scal}^4\norm{A}{2}{2}{0} \right),\label{423}
  \end{eqnarray}
where $c_{scal}=\norm{D\w\g}{}{\infty}{}+\norm{D^2\w\g}{1/2}{\infty}{}
+\norm{ricci_{(M,g)}}{1/2}{\infty}{}+\inj(M,g)^{-1}$ and $c=c(n)$ is a universal
constant.
\end{lemma}
\bb
Since the above estimates are scale invariant we may assume that $c_{scal}=1.$ Define
$\psi:=\g^2\phi.$ With $m=2,~p=4$ and $\a=2/3$ we infer from the multiplicative 
Sobolev-inequality (\ref{66}) and Kato's inequality
\begin{eqnarray}\label{424}
  \norm{\psi}{}{\infty}{}\le c\norm{\psi}{1/3}{2}{}\left(\norm{\na\psi}{}{4}{}
    +\norm{\bet{\psi}A}{}{4}{}+\norm{\psi}{}{4}{} \right)^{2/3}.
\end{eqnarray}
Using partial integration we can interpolate
\begin{eqnarray*}
  \norm{\na\psi}{}{4}{}\le\norm{\psi}{1/2}{\infty}{}\norm{\na^2\psi}{1/2}{2}{}.
\end{eqnarray*}
We estimate the second and third summand by $\norm{\bet{\psi}A}{}{4}{}
  \le\norm{\psi}{1/2}{\infty}{}\norm{\bet{\psi}^{1/2}A}{}{4}{}$ and\\
  $\norm{\psi}{}{4}{}\le\norm{\psi}{1/2}{\infty}{}\norm{\psi}{1/2}{2}{}$ respectively.
Substituting these estimates in (\ref{424}) we obtain using\\ $(a+b)^r\le c(r)(a^r+b^r)$
($a,b\geq 0,~r>0$)
\begin{eqnarray}\label{426}
  \norm{\psi}{4}{\infty}{}\le c\norm{\psi}{2}{2}{}\left(\norm{\na^2\psi}{2}{2}{}
   +\norm{\bet{\psi}^2\bet{A}^4}{}{1}{}+\norm{\psi}{2}{2}{} \right). 
\end{eqnarray}
Now since trivially $\norm{\phi}{4}{\infty}{1}\le\norm{\psi}{4}{\infty}{}$ and
$\norm{\psi}{2}{2}{}\le\norm{\phi}{2}{2}{0}$ we have after substituting in (\ref{426})
\begin{eqnarray}\label{427}
\norm{\phi}{4}{\infty}{1}\le c\norm{\phi}{2}{2}{0}\left(\norm{\na^2\psi}{2}{2}{}
  +\norm{\bet{\psi}^2\bet{A}^4}{}{1}{}+\norm{\phi}{2}{2}{0}\right).
\end{eqnarray}
Taking (\ref{314}) into account, we compute for the first summand
\begin{eqnarray*}
  \bet{\na^2(\g^2\phi)}&\le&2\bet{\na\g}^2\bet\phi+2\g\bet{\na^2\g}\bet\phi
    +4\g\bet{\na\g}\bet{\na \phi}+\g^2\bet{\na ^2\phi}\\
  &\le&\bet{\na ^2\phi}\g^2+c\g\bet{\na \phi}+c\bet\phi\chi_{[\g>0]}+c\g\bet A\bet\phi
\end{eqnarray*}
so that
\begin{eqnarray*}
  \norm{\na ^2\psi}{2}{2}{}&\le& c\is\bet{\na ^2\phi}^2\g^4\mm
    +c\is\bet{\na \phi}^2\g^2\mm
    +c\is\bet\phi^2\left(\chi_{[\g>0]}+\g^2\bet A^2\right)\mm\\
  &\le&c\iz\left(\bet{\na ^2\phi}^2+\bet\phi^2\right)\mm+c\is\bet\phi^2\bet A^4\g^4\mm
\end{eqnarray*}
since from \cref{korollar62} (with $p:=s:=2$) we have after renaming $\e$
\begin{eqnarray*}
  c\is\bet{\na \phi}^2\g^2\mm\le\e\is\bet{\na ^2\phi}^2\g^4\mm+c(\e)\iz\bet\phi^2\mm.
\end{eqnarray*}
Putting the obtained estimates in (\ref{427}) yields (\ref{421}) as
$\norm{\bet\psi^2\bet A^4}{}{1}{}\le\norm{\bet\phi^2\bet A^4}{}{1}{0}.$
Furthermore, assuming (\ref{422}) we get from \cref{lemma42} after absorbing
\begin{eqnarray*}
  \norm{\bet{\psi}^2\bet{A}^4}{}{1}{}=
  \is\bet A^6\g^4\mm&\le&\frac{1}{2}\is\bet{\na^2\!A}^2\g^4\mm+c\norm{A}{4}{2}{0}\\
                    &\le&\frac{1}{2}\norm{\na^2\!A}{2}{2}{0}+\norm{A}{2}{2}{0},
\end{eqnarray*}
which finishes the proof after rescaling.
\eb

\dipl{Proposition 4.4}
\begin{prop}\label{prop44}
  Let $f:\Sigma\times[0,T]\rightarrow M$ be a Willmore flow on a closed surface,
  $\g=\w\g\cf$ as in (\ref{314}) 
  and assume
  \begin{eqnarray}
    \sup\limits_{0< t< T}\norm{A}{2}{2}{0}\le\e^*\label{430}
  \end{eqnarray}
  for $\e^*$ small enough, depending only on the dimension n. Then for any $t\in[0,T]$
  \begin{eqnarray}
    \lefteqn{\hspace{-0em}\is\bet A^2\g^4 d\mu
      +\frac{1}{2}\int\limits_{0}^{t}\is
        \left(  \bet{\na^2 A}^2+\bet A^2\bet{\na A}^2
      +\bet A^6  \right)\g^4d\mu d\t}\label{431}\\
     &\le& \is\bet{A_{0}}^2_{0}\g^4_0 d\mu_{0}
      +c\,C_{scal}^4\int\limits_0^t\int\limits_{[\g>0]}\bet{A}^2\mm d\tau
        +c\norm{DR}{2}{\infty}{}\sup\limits_{0\le t
        \le T}\mu\big([\gamma>0]\big)t,\no
  \end{eqnarray}
  where $C_{scal}=\sum_{i=1}^{2}\norm{D^i\w\g}{1/i}{\infty}{}
  +\sum_{i=0}^{1}\norm{D^{i\!}R}{\frac{1}{i+2}}{\infty}{}+\inj(M,g)^{-1},$ $c=c(n)$ is a
  universal constant and the zero-indexed quantities refer to time t=0.
\end{prop}
\bb Since (\ref{431}) is scale invariant, we may assume that $C_{scal}=1.$
From \cref{lemma34} we know that
\begin{eqnarray}
  \lefteqn{\hspace{0em}\dfrac{d}{dt}\is\bet{ A}^2\g^4\mm
      +\dfrac{3}{4}\is\bet{\na^2 A}^2\g^4\mm}\label{432}\\
   &\le&\is\big(\pe 23+\pe 05+\Q{2}1{1}+\Q 01R\,\big)\s A\,\g^4
      +\big<\Delta\Q 01y+\na\Q 111
      ,A \big>\,\g^4\mm+c\iz\bet{A}^2\mm.
\end{eqnarray}
To begin, we estimate 
$\bet{(\pe23+\pe05)\s A}\le c(\bet{A}^3\bet{\na^2\!A}+\bet{A}^2\bet{\na A}^2+\bet{A}^6).$
Using partial integration we obtain
\begin{eqnarray*}
  \is\big<\Delta\Q 01y+\na\Q 111,\g^4\phi\big>\mm
    \le \is\bet{\Q 01y}\bet{\na^2(\g^4 A)}\mm+\is\bet{\Q 111}\bet{\na(\g^4 A)}\mm.
\end{eqnarray*}
Recalling the definition of the Q-Notation we may further estimate
\begin{itemize}
  \item $\bet{\Q 01y}\le c\bet{DR}\cf+c\bet{A}\bet{R}\cf$
  \item $\bet{\Q 111}\le c\bet{R}\cf\big(\bet{\na A}+\bet{A}^2\big)       
        +c\bet{DR}\cf\bet{A}$
  \item $\bet{\Q 211}\le c\bet{R}\cf\big(\bet{\na^2\!A}+\bet{A}\bet{\na\!A}
        +\bet{A}^3\big)+c\bet{DR}\cf\big(\bet{\na\!A}+\bet{A}^2 \big)$
  \item $\bet{\Q 01R}\le c\bet{R}^2\cf\bet{A}$
\end{itemize}
and moreover
\begin{itemize}
  \item $\bet{\na(\g^4 A)}\le\bet{4\g^3(\na\g)\otimes A+\g^4\na\! A}
        \le c\g^3\bet{A}+\g^4\bet{\na\! A}$
  \item $\bet{\na^2(\g^4 A)}\stackrel{(\ref{na2phi})}{=}
        c\g^2\bet{A}+c\g^3(1+\bet{A})\bet{A}+c\g^3\bet{\na\! A}+\g^4\bet{\na^2\!A}$
\end{itemize}
where we used (\ref{314}) in the last step.
Thus, to estimate the right-hand side of (\ref{432}) we may estimate each of the
following terms
\begin{eqnarray*}
\lefteqn{\hspace{-0em}\is\bet{A}^3\bet{\na^2\!A}\g^4\mm
      \nr{-40.0}{-10}{9}{-42.2}{-11.6}{a}
    +\is\bet{A}^2\bet{\na\!A}^2\g^4\mm
      \nr{-40}{-10}{9}{-42.2}{-12.5}{b}
    +\is\bet{A}^6\g^4\mm\nr{-27}{-10}{9}{-29.0}{-11.8}{c}}\\
  &+&\4\is\bet{DR}\cf\big(
    \bet{A}\g^2\nr{-14}{-10}{9}{-16.4}{-12.4}{d}
    +\bet{A}^2\g^3\nr{-14}{-10}{9}{-15.8}{-11.8}{e}
    +\bet{\na\!A}\g^3\nr{-14}{-10}{9}{-15.5}{-12.3}{f}
    +\bet{\na^2\!A}\g^4\big)\mm\nr{-34}{-10}{9}{-35.9}{-11.4}{g}
    +\2\is\bet{A}^2\g^2\nr{-14}{-10}{9}{-16.1}{-12.3}{h}s
    +\bet{A}^3\g^3\nr{-14}{-10}{9}{-15.2}{-12.4}{i}
    +\bet{A}\bet{\na\!A}\nr{-14}{-10}{9}{-15.2}{-12.0}{j}\g^3\mm\\
  &+&\4\is\bet{A}\bet{\na^2\!A}\nr{-14}{-10}{9}{-15.8}{-12.4}{k}\g^4
    +\bet{\na\!A}^2\g^4\nr{-14}{-10}{9}{-15.2}{-12.6}{l}
    +\bet{A}^2\bet{\na\!A}\nr{-14}{-10}{9}{-17.3}{-11.7}{m}\g^4
    +\bet{A}^4\g^4\nr{-14}{-10}{9}{-16.3}{-11.7}{n}\mm
\end{eqnarray*}
by
$$(\tau+\e^*)\is\bet{\na^2\!A}^2\g^4\mm+c_\tau\norm{A}{2}{2}{0}
  +c_\tau\norm{DR}{2}{\infty}{}\mu\big([\g>0]).$$
At first, using absorption, we get from \cref{lemma42}
\begin{eqnarray}
  \is\bet{A}^6\g^4\mm+\is\bet{A}^2\bet{\na\!A}^2\g^4\mm
    \le c\norm{A}{2}{2}{0}+\e^*\is\bet{\na^2\!A}^2\g^4\mm,\label{lemma211*}
\end{eqnarray}
which estimates integral \it b)\rm\: and \it c)\rm\: already. We have\\\\
\tabcolsep0.3em
\begin{tabular}{rcl}
  $\mbox{\it{a)}}$
    &$\le
    $
    &$\tau\int_\Sigma\bet{\na^2\!A}^2\g^4\mm+c_\tau\int_\Sigma\bet{A}^6\g^4\mm$\\\\
    
$\mbox{\it{d),\;f),\;g)}}$
    &$\le$
    &$c_\tau\norm{DR}{2}{\infty}{}\mu\big([\g>0]\big)+c_\tau\norm{A}{2}{2}{0}
      +\tau\int_\Sigma\bet{\na\!A}^2\g^6\mm+\tau\int_\Sigma\bet{\na^2\!A}^2\g^4\mm$\\\\
\end{tabular}\\
and from \cref{korollar63} we get, letting $\phi=A,~k=1,~p=2$ and $s=4,$
\begin{eqnarray}
  \is\bet{\na\!A}^2\g^4\mm\le\is\bet{\na^2\!A}^2\g^6\mm
    +c\norm{A}{2}{2}{0}\label{na1absch1},
\end{eqnarray}
\tabcolsep0.1em
\begin{tabular}{rcl}
  $~~~~\hspace{0.05em}\mbox{\it e),\;h)\rm}$
    &$\le$
    &$c\norm{A}{2}{2}{0},$\\\\

  $\mbox{\it i)\rm}$
    &$\le$
    &$c\norm{A}{2}{2}{0}+\is\bet{A}^4\g^6\mm
      \le c\norm{A}{2}{2}{0}+\is\bet{A}^{6}\g^{4}\mm$\\
    &$\stackrel{(\ref{lemma211*})}{\le}$
    &$c\norm{A}{2}{2}{0}+\e^*\int_\Sigma\bet{\na^2\!A}^2\g^4\mm,$\\\\

  $\mbox{\it j)\rm}$
    &$\le$
    &$c_\tau\norm{A}{2}{2}{0}+\tau\int_\Sigma\bet{\na\!A}^2\g^4\mm,
    \mbox{$\qquad$(again using $(\ref{na1absch1})$)}$\\\\

  $\mbox{\it k)\rm}$
    &$\le$
    &$c_\tau\norm{A}{2}{2}{0}+\tau\int_\Sigma\bet{\na^2\!A}^2\g^4\mm,$\\\\

  $\mbox{\it l)\rm}$
    &&\7\3(use again \cref{korollar63})\\\\

  $\mbox{\it m),\;n)\rm}$
    &$\le$
    &$c\norm{A}{2}{2}{0}
      +\int_\Sigma\bet{A}^6\g^4\mm+\int_\Sigma\bet{A}^2\bet{\na\!A}^2\g^4\mm$\\
    &$\stackrel{(\ref{lemma211*})}{\le}$
    &$c\norm{A}{2}{2}{0}+\e^*\is\bet{\na^2\!A}^2\g^4\mm.$\\
\end{tabular}\\
~\\\\Putting things together we have now shown that
\begin{eqnarray}
  \lefteqn{\hspace{-2em}\dfrac{d}{dt}\is\bet{ A}^2\g^4\mm
      +\dfrac{3}{4}\is\bet{\na^2 A}^2\g^4\mm}\no\\
  &\le&c_\tau\norm{DR}{2}{\infty}{}\sup\limits_{0\le t\le T}\mu\big([\g>0] \big)
    +(\tau+\e^*)\is\bet{\na^2\!A}^2\g^4\mm+c_\tau\norm{A}{2}{2}{0},\label{dtabsch}
\end{eqnarray}
recalling that $\mu\big([\g>0] \big)\le\mu(\Sigma)
\le \sqrt{2T}\mathcal W(f_0)+\mue{f_0}\,\,(\Sigma)<\infty.$
Adding up (\ref{lemma211*}) and (\ref{dtabsch}), using small curvature concentration
(\ref{430}) and absorbing finally yields
\begin{eqnarray*}
  \lefteqn{\hspace{-5em}\dfrac{d}{dt}\is\bet{ A}^2\g^4\mm
      +\dfrac{3}{4}\is\bet{\na^2 A}^2\g^4\mm
      +\is\bet{A}^6\g^4\mm+\is\bet{A}^2\bet{\na\!A}^2\g^4\mm}\\
  &\le&c\norm{DR}{2}{\infty}{}\sup\limits_{0\le t\le T}\mu\big([\g>0] \big)
    +c\norm{A}{2}{2}{0}.
\end{eqnarray*}
Integrating over $[0,T]$ and rescaling finally yields the claim.\eb
\begin{assumption}\label{skalannahme2}
Assume that $(M,g)$ is of bounded Geometry of order K and that $\varrho$ is chosen
with
$$0<\varrho<\min\Big\{\frac\pi{2\kappa},\inj_{(M,g)}\Big\}.$$
Let $\w\g\in C^\infty(M)$ be a cutoff function satisfying
\begin{eqnarray}
  \chi_{B_{\varrho/2}(x_0)}\le\w\g
 \le\chi_{B_{\varrho}(x_0)}\mbox{~~~and~~~}\sum\limits_{i=1,2}\varrho^i
 \norm{D^i\w\g}{1/i}{L^\infty(M)}{}<c(n),\label{cgamma}
\end{eqnarray}
where $c(n)$ is a universal constant. Then choose $1\le R_K<\infty$ with
\begin{eqnarray}
\varrho\sum\limits_{i=1}^{K}\norm{D^i\!R}{1/(i+2)}{L^\infty(M,g)}{}<R_K.
\label{riemannbound}
\end{eqnarray}
\end{assumption}
\cref{452} shows the existence of such a cutoff function for any $x_0\in M.$\er
\dipl{Proposition 4.5}
\begin{prop}\label{prop45}
  Let $f:\Sigma^2\times [0,T]\rightarrow (M^n,g)$ be a proper Willmore flow. Then
for $\phi=\na^m\!A,$ $m\in\N_0$ and $\g=\w\g\cf$ with $\w\g$ as in (\ref{cgamma}), we
have for all $s\geq 2m+4$
  \begin{eqnarray}
    \lefteqn{\hspace{-0em}\frac{d}{dt}\is\bet{\phi}^2\g^s\mm
    +\frac{1}{2}\is\bet{\na^2\phi}^2\g^s\mm }\label{434} \\
    &\le&c(R_{m+1})\big(\varrho^{-4}+\norm{A}{4}{\infty}{0}\big)
      \is\bet{\phi}^2\g^s\mm\no\\
      &&+c(R_{m+1})\big(1+\norm{\varrho A}{\max\{4,2m\}}{\infty}{0}\big)
      \norm{A}{2}{2}{0}\varrho^{-2m-4}
    +c\norm{D^{m+1}\!R}{2}{L^\infty(M,g)}{}\,\mu\big([\gamma>0]\big),\no
  \end{eqnarray}
  where $c$ only depends on $n,m$ and the constant in (\ref{cgamma}).
\end{prop}
\bb Since (\ref{434}) is scale invariant, we may assume that 
$\varrho= R_{m+1}.$ We estimate the terms in \cref{lemma34}, i.e.
$$\is\big(\pe{m+2}3+\pe m5+\Q{m+2}1{m+1}+\Q m1R\,\big)\s\phi\,\g^s
      +\big<\na^m\Delta\Q 01y+\Delta\Q m1m+\na\Q {m+1}1{m+1}
      ,\phi \big>\,\g^s\mm.$$
Analogously to Kuwert and Sch\"atzle (\cite{KS02}, Proposition 4.5, (4.15))
we have
\begin{eqnarray*}
  \lefteqn{\hspace{-2em}\int_{\Sigma} \big(\pe{m+2}3+\pe m5 \big)\s \phi\,\g^s\mm}\\
  &\le&\dfrac{1}{16}\is\bet{\na^2\phi}^2\g^s\mm
    +c\norm{A}{4}{\infty}{0}\is\bet{\phi}^2\g^s\mm
    +c\big(1+\norm{A}{4}{\infty}{0}\big)\norm{A}{2}{2}{0}.
\end{eqnarray*}
Similar to the proof of \cref{prop44} we use partial integration obtaining
\begin{eqnarray}
  \lefteqn{\hspace{-0em}\is\big<\na^m\Delta\Q 01y+\Delta\Q m1m+\na\Q {m+1}1{m+1}
      ,\g^s\phi \big>+\big(\Q{m+2}1{m+1}+\Q m1R\,\big)\s\phi\,\g^s\mm}\label{qtermab10}\\
  &\le&c\is\bet{\Q m1x}\bet{\na^2(\g^s\phi)}\mm
    +\is\bet{\Q {m+1}1{m+1}}\bet{\na(\g^s\phi)}\mm
    +\is\big(\bet{\Q {m+2}1{m+1}}+\bet{\Q m1R}\big)\bet{\phi}\g^s\mm.\no
\end{eqnarray}
It will turn out that problems arise when trying to estimate the first integral on the
right-hand side of (\ref{qtermab10}) in case the homogeneity $\nu$ of the second
fundamental form is too large in $\Q m1{}.$ More precisely, this is the case when
$\nu>m/2+1.$
Fortunately, we may assume that $\nu\le m/2+1.$ This can be seen as follows:
With the obvious notation, we have the decomposition
$\Q m1{}=\Q m1{\nu\le m/2+1}+\Q m1{\nu>m/2+1}.$ Since we can write
$\Q m1{\nu> m/2+1}=\Q m1m$ and because $\na\Q m1m=\Q {m+1}1{m+1},$ we now estimate
\begin{eqnarray*}
  \is\big<\na^m\Delta\Q 01y,\g^s\phi \big>\mm
  &\le&\is\bet{\Q m1{\nu\le m/2+1}}\bet{\na^2(\g^s\phi)}\mm
    +\is\bet{\na\Q m1{\nu>m/2+1}}\bet{\na(\g^s\phi)}\mm\\
  &\le&\is\bet{\Q m1{\nu\le m/2+1}}\bet{\na^2(\g^s\phi)}\mm
    +\is\bet{\Q {m+1}1{m+1}}\bet{\na(\g^s\phi)}\mm,
\end{eqnarray*}
which justifies the above assumption.
Thus, by definition of the Q-Notation, we may estimate
(recall that $\mu=i_1+\ldots+i_\nu,$ i.e. $\nu=0$ implies $\mu=0$):
\begin{eqnarray*}
  \bet{\Q m1{\nu\le m/2+1}}
    &\le\7&\sum\limits_{\substack{r+\mu+\nu=m+1\\\mu\le m\\\nu\le m/2+1}}
    \6\bet{D^{r\!}R}\cf\bet{\na^{i_1}\!A}\c\ldots\c\bet{\na^{i_\nu}\!A}
  \le\3\sum\limits_{\substack{1\le\mu+\nu\le m+1\\\nu\le m/2+1}}
    \3\bet{\na^{i_1}\!A}\c\ldots\c\bet{\na^{i_\nu}\!A}
    +\bet{D^{m+1\!}R}\cf
\end{eqnarray*}
and analogously
\begin{eqnarray*}
  \bet{\Q {m+1}1{m+1}}&\le&\7\sum\limits_{\substack{r+\mu+\nu=m+2\\ r,\mu\le m+1}}
    \6\bet{D^{r\!}R}\cf\bet{\na^{i_1}\!A}\c\ldots\c\bet{\na^{i_\nu}\!A}
  \le\3\sum\limits_{\substack{1\le\mu+\nu\le m+2\\\mu,\nu\le m+1}}
    \3\bet{\na^{i_1}\!A}\c\ldots\c\bet{\na^{i_\nu}\!A}
    +\bet{A}^{m+2},\\
  \bet{\Q {m+2}1{m+1}}
    &\le&\7\sum\limits_{\substack{r+\mu+\nu=m+3\\\mu\le m+2\\r\le m+1\\\mbox{~}}}
    \6\bet{D^{r\!}R}\cf\bet{\na^{i_1}\!A}\c\ldots\c\bet{\na^{i_\nu}\!A}\\
  &\le&\6\sum\limits_{\substack{1\le\mu+\nu\le m+3\\\mu,\nu\le m+1}}
    \3\bet{\na^{i_1}\!A}\c\ldots\c\bet{\na^{i_\nu}\!A}
    +\bet{\na^2\phi}+\bet{A}^{m+3}+\bet{A}^{m+2}+\bet{A}^{m+1}\bet{\na\!A},\\
  \bet{\Q {m}1R}&\le&\7\4\sum\limits_{\substack{r_1+r_2
    +\mu+\nu=m+1\\\mu\le m\\\nu\geq 1}}
    \7\bet{D^{r_1\!}R}\cf\bet{D^{r_2\!}R}\cf
    \bet{\na^{i_1}\!A}\c\ldots\c\bet{\na^{i_\nu}\!A}
  \le\5\sum\limits_{\substack{1\le\mu+\nu\le m+1}}
    \3\bet{\na^{i_1}\!A}\c\ldots\c\bet{\na^{i_\nu}\!A}.
\end{eqnarray*}
Using that $\bet{\na(\g^s\phi)}=\bet{s\g^{s-1}(\na\g)\otimes\phi+\g^s\na\phi}
\le c\bet{\phi}\g^{s-1}+c\bet{\na\phi}\g^s$
and
\begin{eqnarray*}
  \bet{\na^2(\g^s\phi)}&\stackrel{(\ref{na2phi})}{\le}&
    c\g^{s-2}\bet{\phi}+c\g^{s-1}(1+\bet{A})\bet{\phi}+c\g^{s-1}\bet{\na\phi}
    +c\g^s\bet{\na^2\phi}\\
  &\le&c\sum\limits_{p=0}^1\g^{s+p-2}\bet{\na^p\phi}+c\g^{s-2}\bet{A}\bet{\phi}
    +c\bet{\na^2\phi}\g^s
\end{eqnarray*}
we arrive, collecting terms and rearranging, at $(\ref{qtermab10})\le$
\begin{eqnarray*}
    &c&\7\4\sum\limits_{\substack{1\le \mu+\nu\le m+3\\\mu,\nu\le m+1}}\sum_{p=0}^1
     \is\bet{\na^{i_1}\!A}\c\ldots\c\bet{\na^{i_\nu}\!A}\bet{\na^p\phi}\g^{s+p-2}\mm\\
    &+&c\6\sum\limits_{\substack{1\le\mu+\nu\le m+1\\\nu\le m/2+1}}
      \is\bet{\na^{i_1}\!A}\c\ldots\c\bet{\na^{i_\nu}\!A}\bet{\na^2\phi}\g^s\mm
  +c\norm{D^{m+1\!}R}{}{L^\infty(M)}{}\sum_{p=0}^2\is\bet{\na^p\phi}\g^{s+p-2}\mm\\
    &+&c\sum\limits_{p=0}^1\bet{A}^{m+2}\bet{\na^p\phi}\g^{s+p-2}\mm
      +c\is\bet{A}^{m+3}\bet{\phi}\g^s\mm
    +c\is\bet{A}^{m+1}\bet{\na\!A}\bet{\phi}\g^s\mm.
\end{eqnarray*}
We now estimate each of the integrals above:\\
\begin{eqnarray}
  \mbox{1)}&\le&\4\sum\limits_{\substack{1\le \mu+\nu\le m+3\\\mu,\nu\le m+1}}
    \sum_{p=0}^1
    \Bigg(\is\bet{\na^{i_1}\!A}^2\c\ldots\c\bet{\na^{i_\nu}\!A}^2
    \g^{s-2}\mm\Bigg)^{1/2}\Bigg(\is\bet{\na^p\phi}^2\g^{s+2p-2}\mm
    \Bigg)^{1/2}\!\!\!\!.\label{firstint}~~~~
\end{eqnarray}
For the first integral in (\ref{firstint}) we use \cref{korollar66} with 
$k:=\mu\le m+1,~\tilde s=s-2\geq 2m+2\geq 2k,$ $r=2\nu\geq 2$
obtaining 
\begin{eqnarray*}
  \is\bet{\na^{i_1}\!A}^2\c\ldots\c\bet{\na^{i_\nu}\!A}^2\g^{s-2}\mm
  &\le&c\norm A{2\nu-2}\infty0 
    \Bigg(\is\bet{\na^\mu\!A}^2\g^{s-2}\mm+\norm{A}{2}{2}{0}\Bigg)\\
  &\le&c\norm A{2\nu-2}\infty0\Bigg(\tau
  \is\bet{\na^{2\!}\phi}^2\g^{s}\mm+(1+\tau^{-1})\norm{A}{2}{2}{0}\Bigg),
\end{eqnarray*}
where we used \cref{korollar64} with $k:=\mu,$ $l:=m+2$ and $\tilde s:=s-2
\geq 2(l-1).$ For the
second integral
in (\ref{firstint}) we want to employ \cref{korollar64} with $k:=m+p,$ $l:=m+2,$  
$\phi:=A$ and $\tilde s:=s+2p-2\geq 2(l-1)$ so that further
\begin{eqnarray*}
  (\ref{firstint})&\le&c\sum\limits_{\nu=1}^{m+1}
    \norm{A}{\nu-1}{\infty}{0}
    \Bigg(\tau\is\bet{\na^{2\!}\phi}^2\g^{s}\mm
    +(1+\tau^{-1})\norm{A}{2}{2}{0}\Bigg)\\
  &\le&\e\is\bet{\na^{2\!}\phi}^2\g^{s}\mm
    +c_\e(1+\norm{A}{2m}{\infty}{0})\norm{A}{2}{2}{0}.\hspace{17em}
\end{eqnarray*}
\begin{eqnarray*}
  \mbox{2)}
    &\le&\e\is\bet{\na^2\phi}^2\g^s\mm
      +c_{\e}\6\sum\limits_{\substack{1\le\mu+\nu\le m+1\\\nu\le m/2+1}}
      \is\bet{\na^{i_1}\!A}^2\c\ldots\c\bet{\na^{i_\nu}\!A}^2\g^s\mm\\
    &\stackrel{(\ref{65})}{\le}&\e\is\bet{\na^2\phi}^2\g^s\mm
      +c_\e\2
      \5\sum\limits_{\substack{0\le k+\nu\le m\\0\le\nu\le m/2}}\5
      \norm{A}{2\nu}{\infty}{0}\Bigg(
      \is\bet{\na^{k}\!A}^2\g^s\mm+\norm{A}{2}{2}{0}\Bigg)\hspace{2em}\\
    &\stackrel{(\ref{63})}{\le}&
      \e\is\bet{\na^2\phi}^2\g^s\mm+c_\e(1+\norm{A}{m}{\infty}{0})
      \Bigg(\tau\is\bet{\na^2\phi}^2\g^s\mm+(1+\tau^{-1})\norm{A}{2}{2}{0}\Bigg)\\
    &\le&
      \e\is\bet{\na^2\phi}^2\g^s\mm+c_\e(1+\norm{A}{2m}{\infty}{0})
      \norm{A}{2}{2}{0}.  
\end{eqnarray*}
\begin{eqnarray*}
  \mbox{3)}
    &\le&\e\is\bet{\na^2\phi}^2\g^s\mm+c_\e\norm{D^{m+1\!}R}{2}{L^\infty(M)}{}
      \mu\big([\g>0]\big)+c_\e\is\bet{\phi}^2\g^s\mm+c_\e\norm{A}{2}{2}{0},
\end{eqnarray*}
where we used $\norm{D^{m+1\!}R}{}{L^\infty(M)}{}\bet{\na^p\phi}\g^{s+p-2}
\le \sigma\bet{\na^p\phi}^2\g^s
+c_\sigma \norm{D^{m+1\!}R}{2}{L^\infty(M)}{}\chi_{[\g>0]}$
and $$\is\bet{\na\phi}^2\g^s\mm\stackrel{(\ref{63})}{\le}
\tau\is\bet{\na^2\phi}^2\g^s\mm+c_\tau\norm{A}{2}{2}{0}.$$
\begin{eqnarray*}
  \mbox{4)}
    &\le&c\iz\bet{A}^{2m+2}\mm
      +c\sum\limits_{p=0}^1\is\bet{A}^2\bet{\na^p\phi}^2\g^s\mm\\
    &\stackrel{(\ref{63})}{\le}&
    \norm{A}{2m}{\infty}{0}\norm{A}{2}{2}{0}+c\norm{A}{2}{\infty}{0}
      \Bigg(\e\is\bet{\na^2\phi}^2\g^s\mm+(1+\e^{-1})\norm{A}{2}{2}{0} \Bigg)\\
    &\le&\e\is\bet{\na^2\phi}^2\g^s\mm
      +c_\e\Big(1+\norm{A}{\max\{4,2m\}}{\infty}{0}\Big)\norm{A}{2}{2}{0}.\\
  \mbox{5)}
    &\le&c\norm{A}{4}{\infty}{0}\is\bet{\phi}^2\g^s\mm
      +c\norm{A}{2m}{\infty}{0}\norm{A}{2}{2}{0}\\
    &\le&c\norm{A}{4}{\infty}{0}\is\bet{\phi}^2\g^s\mm
      +c\Big(1+\norm{A}{\max\{4,2m\}}{\infty}{0} \Big)\norm{A}{2}{2}{0}.
\end{eqnarray*}
\begin{eqnarray*}
  \mbox{6)}
    &\le&c\is\bet{A}^{2m+2}\g^s\mm
      +\is\bet{\na\!A}^2\bet{\phi}^2\g^s\mm\\
    &\stackrel{(\ref{65})}{\le}&
      c\norm{A}{2m}{\infty}{0}\norm{A}{2}{2}{0}
      +c\norm{A}{2}{\infty}{0}\Bigg(
      \is\bet{\na\phi}^2\g^s\mm+\norm{A}{2}{2}{0}\Bigg)
\end{eqnarray*}
\begin{eqnarray*}
    &\stackrel{(\ref{63})}{\le}&
    c\norm{A}{2m}{\infty}{0}\norm{A}{2}{2}{0}
      +c\norm{A}{2}{\infty}{0}\Bigg(
      \tau\is\bet{\na^2\phi}^2\g^s\mm+(1+\tau^{-1})\norm{A}{2}{2}{0}\Bigg)\\
    &\le&\e\is\bet{\na^2\phi}^2\g^s\mm
      +c_\e\Big(1+\norm{A}{\max\{4,2m\}}{\infty}{0} \Big)\norm{A}{2}{2}{0}.
\end{eqnarray*}
Putting things together and choosing $\e$
small enough, the desired estimate (\ref{434}) follows after rescaling.\eb

The next lemma proves a H\"older estimate of the volume function. It was taken from
Kuwert (\cite{lecKuw}, Lemma 3.5) and carried over to the Riemannian setting. 
\dipl{Lemma 4.6}
\begin{lemma}\rm(\cite{lecKuw}).\it\label{lemma46}~Let
  $f:\Sigma\times [0,T]\to M^n$ be a Willmore flow on a closed surface. Then the
  volume function   $vol_\Sigma(t):=\int_{\Sigma}d\mu_{\tildeg(t)}$ for the
  volume of $(\Sigma,\tildeg(t))$ at time t satisfies
  \begin{eqnarray}
    \bet{vol_\Sigma(t_1)-vol_\Sigma(t_2)}
    \le\sqrt{2}\bet{t_1-t_2}^{1/2}\W(f_0)\label{443}
  \end{eqnarray}
for $t_1,t_2\in[0,T],$ where
$\W(f_0):=\frac{1}{2}\int_{\Sigma}\bet{H}^2\mm\strichklein{t=0}\quad.$
\end{lemma}

\bb
By definition of the Willmore flow, we have with ${\bf W}:={\bf W}(f)=grad_{L^2}\W$
\begin{eqnarray}
  \dfrac{d}{dt}\dfrac{1}{2}\is\bet{H}^2\mm=-\is\bet{{\bf W}}^2\mm,\label{444}
\end{eqnarray} 
and therefore
\begin{eqnarray}
  \dfrac{1}{2}\is\bet{H}^2\mm
    \le\dfrac{1}{2}\is\bet{H}^2\mm\strich{t=0}\quad=\W(f_0)<\infty.\label{445}
\end{eqnarray}
We further obtain from (\ref{444})
\begin{eqnarray}
  \int\limits^{T}_{0}\is\bet{{\bf W}}^2\mm\,dt
  &=&\frac{1}{2}\is\bet{H}^2\mm\strich{t=0}\quad
    -\frac{1}{2}\is\bet{H}^2\mm\strich{t=T}\no \\
  &\le&\dfrac{1}{2}\is\bet{H}^2\mm\strich{t=0}=\W(f_0).\label{446}
\end{eqnarray}

Now for the derivative of the volume function we have
\begin{eqnarray}
  \bet{vol'(t)}^2&=&\bigg\lvert\partial_t\is\mm\bigg\rvert^2\no \\
  &\stackrel{(\ref{158})} {=}&
    \bigg\lvert\is\langle H,{\bf W} \rangle \mm\bigg\rvert^2\no \\
  &\le&\is\bet{H}^2\mm\is\bet{{\bf W}}^2\mm\no \\
  &\stackrel{(\ref{445})}{\le}&2\W(f_0)\is\bet{{\bf W}}^2\mm,\label{447}
\end{eqnarray}
so that we obtain for arbitrary $t_1,t_2\in [0,T]$
\begin{eqnarray*}
  \bet{vol(t_2)-vol(t_1)}&=&\bigg\lvert\int\limits_{t_1}^{t_2}vol'(t)dt\bigg\rvert \\
  &\le& \bet{t_2-t_1}^{\frac{1}{2}}\bigg\lvert\int\limits^{t_2}_{t_1}vol'(t)^2dt 
    \bigg\rvert^{\frac{1}{2}} \\
  &\stackrel{(\ref{447})}{\le}&\sqrt{2\W(f_0)}\bet{t_2-t_1}^{\frac{1}{2}}
    \bigg(\int\limits^{T}_{0}\is\bet{{\bf W}}^2\mm\,dt \bigg)^{\frac{1}{2}}\\
  &\stackrel{(\ref{446})}{\le}& \sqrt{2}\W(f_0)\bet{t_2-t_1}^{\frac{1}{2}}.
\end{eqnarray*}
The claim follows.\eb\\
Now in addition to \cref{skalannahme2} assume that $\Lambda_K\geq R_K$ is chosen such
that
\begin{eqnarray}
  \varrho\big(\mue{f_0}\:(\Sigma)+\varrho^2\W(f_0)\big)\sum\limits_{i=1}^{K}
    \norm{D^i\!R}{\frac{3}{i+2}}{L^\infty(M,g)}{}\le\lam K<\infty.\label{anfangsbed}
\end{eqnarray}
\dipl{Proposition 4.7}
\addcontentsline{toc}{subsection}{Interior estimates I}
\begin{prop}{\rm (Interior estimates I).}\label{prop47}~
  Let $f:\Sigma\times [0,T]\to M^n$ be a Willmore flow of a closed surface, and
  $\g=\w\g\cf$ as in (\ref{cgamma}). Assume that $\varrho^{-4}T\le\widehat T$ holds. 
  If $$\sup\limits_{0\le t\le T}\iz\bet{A}^2\mm<\e^*,$$
  where $\e^*>0$ is small enough depending only on n, then for any $m\in\N_0$
  \begin{eqnarray*}
    \norm{\na^m\!\!A}{}{\infty}{1}
      \le c\big(n,m,\widehat T,\lam{m+3},\a_0(m+2)\big)\varrho^{-(m+1)},
  \end{eqnarray*}
where $\sum_{j=0}^{l}\varrho^{\,j}\norm{\na^j\!A}{}{2}{0}\strich{t=0}\;\;\;\le\a_0(l),$
$\lam{m+3}$ is as in (\ref{anfangsbed}), and c also depends on the constant in
(\ref{cgamma}).
\end{prop}
\bb By scale invariance, we may assume that $\varrho=1.$ In the sequel, the
constants
$c$ may additionally depend on the constant in (\ref{cgamma}). Now for
$0\le\sigma<\t\le 1$ we use cutoff functions $\psi_{\sigma,\t}\in C^\infty(\R)$
with $0\le\psi_{\sigma,\t}\le 1$ satisfying $\psi_{\sigma,\tau}(s)=0$ for $s\le \sigma$ 
and $\psi_{\sigma,\t}(s)=1$ for $s\geq\t.$ Define
$\g_{\sigma,\t}:=\psi_{\sigma,\tau}\circ\g.$ With $\sigma:=0,~\t:=\frac{1}{2}$ and wlog
$\e^*\le 1$ it follows from (\ref{431}) that
\begin{eqnarray*}
  \int\limits_{0}^{T}\abint {\frac{1}{2}}\bet{\na^2\!A}^2+\bet{A}^6\mm\,dt
     &\le&  \e^*
       +c\,\e^* C_{scal}^{\,4}T+c\norm{DR}{2}{\infty}{}\sup\limits_{0\le t
       \le T}\bet{\Sigma}T\\
     &\le& c(T,\lam 1)
\end{eqnarray*}
where we used
\begin{eqnarray}
  \norm{DR}{2}{\infty}{}\sup\limits_{0\le t\le T}\bet{\Sigma}
    \stackrel{(\ref{443})}{\le}
    c(T,R_1)\norm{DR}{}{\infty}{}
    \big(\W(f_0)+\mue{f_0}\:(\Sigma) \big)
    \le c(T,\lam 1).\label{riemvolabsch}
\end{eqnarray}
From (\ref{423}) with
$\sigma=\frac{1}{2}$ and $\t=\frac{3}{4}$ we get, integrating over $[0,T],$
\begin{eqnarray*}
  \int_0^T \norm{A}{4}{\infty}{\frac{3}{4}}dt
    &\le& c\e^*\big(c(T,\lam 1)
      +c\e^* \big) \\
    &\le& c(T,\lam 1).
\end{eqnarray*}\notiz{das mit G noch besser machen.alles allgemein, dann spezialisieren!}
In order to use the integral estimate in \cref{prop45} we restrict to 
$0\le m\le 2$ at first, and let $\sigma=3/4$ and $\t=7/8.$ Since we have to go
through the following estimates for different values for
$\sigma,\t$ later on when assuming $m>2,$ we do not substitute $\sigma,\t$ in the first
instance. From (\ref{434}) with
$\phi =\na^m\!A,~s=2m+4$ and $0\le t\le T$ we obtain
\begin{eqnarray}
  \lefteqn{\hspace{-5em}\dfrac{d}{dt}\is\bet{\phi}^2\g^s_{\sigma,\t}\mm
    +\dfrac{1}{2}\abint{\t}\bet{\na^2\phi}^2\mm }\label{vonhiernochmal} \\
  &\le& c( R_{m+1})\big(1+\norm{A}{4}{\infty}{\sigma}
\big)\bigg(\is\bet{\phi}^2\g^s_{\sigma,\t}\mm
    +\dfrac{1}{2}\int\limits_{0}^{t}\abint\t\bet{\na^2\phi}^2\mm\,dt' \bigg)\no \\
  &&+c(R_{m+1})\big(1+\norm{A}{\max\{4,2m\}}{\infty}{\sigma} \big)\e^*+c(T,\lam{m+1}),\no
\end{eqnarray} 
where we have estimated the last summand in (\ref{434}) similar to (\ref{riemvolabsch}).
Now from Gronwall's Lemma we obtain
\begin{eqnarray}
  \lefteqn{\hspace{-6em}\sup\limits_{0\le t\le T}\abint\t\bet{\na^{m}\!\!A}^2\mm
    +\dfrac{1}{2}\int\limits_0^T\abint\t\bet{\na^{m+2}\!\!A}^2\mm\,dt}\no\\
  &\le&\exp\bigg(c( R_{m+1})\int\limits_0^T
    \big(1+\norm{A}{4}{\infty}{\sigma}\big)dt\bigg)
    \cdot\bigg[\int\limits_{[\g_{_0}>0]}\bet{\na^{m}\!\!A_{_0}}^2_{_0}d\mu_{_0}\no\\
  &&+c( R_{m+1})\int\limits_0^T\big(1+\norm{A}{\max\{4,2m\}}{\infty}{\sigma} \big)dt
    +c(T,\lam{m+1})\bigg]\label{gronlemma}\\
  &\le& c(n,m,T,\lam{m+1},\a_0(m)).\no
\end{eqnarray}
From (\ref{423}) we get for $\g:=\g_{\sigma,\t},~\sigma=\frac{7}{8},~\t=\frac{15}{16}$
using absorption
\begin{eqnarray*}
  \norm{A}{4}{\infty}{\frac{15}{16}}&\le& c\e^*
    \Big(c\big(n,T,\lam 3,a_0(2)\big)
    +c\e^*\Big) \\
  &\le& c\big(n,T,\lam 3,\a_0(2)\big),
\end{eqnarray*}
so that we now have for all $\nu\in\N_0$
\begin{eqnarray}
  \int\limits_{0}^{T}\norm{A}{\nu}{\infty}{\frac{15}{16}}dt
    \le c\big(n,T,\lam 3,\a_0(2),\nu\big).\label{449}
\end{eqnarray}
Now if $m>2,$ we let $\sigma=15/16$ and
$\t=31/32$ and obtain, again using (\ref{gronlemma}),
\begin{eqnarray}
  \sup\limits_{0\le t\le T}\abint{\frac{31}{32}}\bet{\na^{m}\!\!A}^2\mm
  +\dfrac{1}{2}\int\limits_0^T\abint{\frac{31}{32}}\bet{\na^{m+2}\!\!A}^2\mm\,dt
    &\le&c\big(n,m,T,\lam 3,\lam{m+1},a_0(2),a_0(m)\big)\no\\
  &=:&\w c(m),\label{gronwallfinal}
\end{eqnarray}
where we used 
\begin{eqnarray}
  \int\limits_0^T\big(1+\norm{A}{2m}{\infty}{\frac{15}{16}} \big)dt
  \stackrel{(\ref{449})}{\le}c\big(n,m,T,\lam 3,\a_0(2)\big).
\end{eqnarray} 
From the above, the bound (\ref{gronwallfinal}) now holds for any $m\in\N_0.$
Finally, the claim follows from \cref{lemma43}, if we set 
$\phi:=\na^m\!A,~\g:=\g_{\sigma,\t},~\sigma=31/32$ and $\t=1:$
\begin{eqnarray*}
  \norm{\na^m\!A}{4}{\infty}{1}&\le&\w c(m)\big(\w c(m+2)
    +\norm{A}{4}{\infty}{\frac{31}{32}}\w c(m)+\w c(m) \big) \\
  &\le& c(n,m,T,\lam{m+3},\a_0(m+2)).
\end{eqnarray*}
Rescaling yields the claim.\eb

\section{Estimating the lifespan}\label{lifespanEst}
If $g$ is the given background metric on $M$, consider for
$\r>0$ the scaled metric 
$\rg:=\r^{-2}g$ and define an isometric immersion
\begin{eqnarray}
  f_\r:\big( \Sigma,{}^\r\tildeg \big)\times\big[0,\frac{1}{\varrho^{4}}T \big)
      \rightarrow(M,\rg)\label{51}
\end{eqnarray}
by
\begin{eqnarray}
  f_\r(x,t):=f(x,\r^4t),
\end{eqnarray}
i.e. we rescale parabolically.
By definition we have ${}^\r\tildeg:=f_{\r,t}^*(\rg)=\r^{-2}\tildeg.$ As one easily 
verifies from the scaling properties of the various geometric quantities (see e.g.
\cite{dipl}, Lemma 3.2), we have that $f$ is a (maximal) Willmore flow, if and only
if $f_\r$ is a (maximal) Willmore flow.\\\\
We now want to prove \cref{lifespan} (cf. the respective euclidean proof in
\cite{KS02}). Since the statement of \cref{lifespan} is scale invariant, we
may assume that for $\r$ therein given we have $\r=1.$
We first want show that
$t^+(1)=:t^+<T.$ By definition, we have
$t^+\le T.$ To prove that $t^+<T,$ we will lead the case
\begin{eqnarray}
  t^+=T<\infty,\label{notpossible}
\end{eqnarray}
to a contradiction. More precisely, in this case we can extend the flow,
contradicting the assumed maximality of the lifespan $T.$\\\\
Choose $\Lambda_K,R_K\geq 1$ with
\begin{eqnarray}
    \sum\limits_{i=0}^{K}\norm{D^{i\!}R}{\frac{1}{i+2}}{L^\infty(M,g)}{}
    \le R_K\le\Lambda_K\mbox{\quad and \quad}
    \big(\mue{f_0}\;(\Sigma)+\W(f_0)\big)\sum\limits_{i=3}^{K}
    \norm{D^{i\!}R}{\frac{3}{i+2}}{L^\infty(M,g)}{}
    \le\lam K\label{annahmecurv}
\end{eqnarray}
for $K\in\N_0$ to be determined. From \cref{prop47} we then get, letting
$\e_0^2:=\e^*$ and $\widehat T:=t^+,$
\begin{eqnarray}
  \sum_{i=0}^{m}\norm{\na^i\!A}{}{\infty}{}\le A_m\le
    c(m,n,t^+,\lam{m+3},\a_0(m+2),(\ref{cgamma})),\label{aabsch}
\end{eqnarray}
for some $A_m\in\R,$ allowing $c$ also to depend on the constant in (\ref{cgamma}).
Now by definition of the Q-terms
\begin{eqnarray*}
  \norm{\Q {k}ls}{}{\infty}{}
    &\le& c(k,l,n)\!\!\!\!\!\sum\limits_{\substack{r+\mu+\nu= k+l\\
      \mu\le k\\r\le s } }^{}\6\big(\norm{D^r\!R}{}{\infty}{}
      \norm{\na^{i_1}\!\!A}{}{\infty}{}\c\ldots\c
      \norm{\na^{i_\nu}\!\!A }{}{\infty}{} \big)\\
    &\le&c\big(k,l,n, R_{s},A_k\big)
\end{eqnarray*} 
and analogously
\begin{eqnarray*}
  \norm{\Q {k}lR}{}{\infty}{}
  &\le& c(k,l,n)\7\sum\limits_{\substack{r_1+r_2+\mu+\nu= k+l\\
      \mu\le k\\\nu\geq 1}}^{}\7\big(
      \norm{D^{r_1}\!R}{}{\infty}{}\norm{D^{r_2}\!R}{}{\infty}{}
      \norm{\na^{i_1}\!\!A}{}{\infty}{}\c\ldots\c
      \norm{\na^{i_\nu}\!\!A }{}{\infty}{} \big)\\
  &\le& c(k,l,n, R_{k+l-1},A_k).
\end{eqnarray*}
Now since
\begin{eqnarray*}
  V=\partial_t f&\stackrel{(\ref{315})}{=}&1\s \na^2\!A+\pe 03+\Q01{0}
\end{eqnarray*}
we get from
\begin{eqnarray}
  h:=\na_{\partial_t}\w g\stackrel{(\ref{157})}{=}-2\langle A,\partial_t f \rangle
    =A\s\na^2 A+\pe 04+\Q02{0}\label{metricevol4}
\end{eqnarray}
for the $m$-th covariant derivative 
\begin{eqnarray}
  \norm{\na^mh}{}{\infty}{}&\le& c(m,n,A_{m+2}, R_{m}).\label{MteVonH}
\end{eqnarray}
From (\ref{159}) we obtain for the evolution of the Christoffel symbols
induced by $\w g$
\begin{eqnarray*}
  \na_{\partial_t}\!{}^\Sigma\na&=&V\s\na A+A\s\na V+\Q 00{}\s V \\
  &=&\pe 32+\pe14+\Q 121+\Q 210+\Q030+\Q01R,
\end{eqnarray*} 
so that
\begin{eqnarray}
  \norm{\na^m(\na_{\partial_t}\!{}^\Sigma\na)}{}{\infty}{}
    &\le& c(m,n,A_{m+3}, R_{m+1}).\label{510}
\end{eqnarray}
Now to extend the flow, we take (\ref{aabsch}) into account
for the flow in a neighbourhood of T.
Using \cref{harmcoord1}, we can now fix, for $K$ as in (\ref{anfangsbed}) to be
determined later on, a harmonic radius $r_0=c(n,K, R_K)>0$ and a countable atlas
$$\{\psi_{i,r_0}:B^g_{r_0}(p_i)\to V_i\subset B^n_{2r_0} \}_{i\in\N}$$ of harmonic
coordinates for $M,$ such that in these coordinates 
\begin{itemize}
  \item [i)] $\frac{1}{2}\d\le(g_{\a\b})\le 2\d$ on $B^g_{r_0}(p_i)$
  \item [ii)] $\sup_{B^g_{r_0}(p_i)}\bet{\p^\g g_{\a\b}}\le c(n,K, R_K)$ for all
    $1\le\bet{\g}\le K+1$
  \item [iii)] $\sup_{B^g_{r_0}(p_i)}\bet{\p^\g \Gamma^\a_{\!\b\tau}}
    \le c(n,K, R_K)$ for all $0\le\bet{\g}\le K,$
\end{itemize}
and such that $\bigcup_{i\in\N}B_{r_0}(p_i)$ is a uniformly locally finite cover as
described in \cref{hilfslemma51}.
We now show that there exists a well-defined coordinate representation for $f$ in a 
neighbourhood of $T$. Since 
$\lVert\p_t f \rVert_{\infty}\le c(n,A_2,R_0)$ we clearly have for any $t_1,t_2\in[0,T)$
and $x\in\Sigma$
\begin{eqnarray*}
  d^g\big(f(x,t_1),f(x,t_2) \big)\le\Big\lvert\int_{t_1}^{t_2}\bet{\p_t f}_{g}(x,\tau)
  d\tau\Big\rvert\le c(n,A_2,R_0)\bet{t_2-t_1}.
\end{eqnarray*}
Thus for $0<\bet{T-\tau}\le\frac {r_0}4c(n,A_2,R_0)$ we get,
because $f\big(B^{\tildeg(\tau)}_{r_0/2}(x)\big)\subset B_{r_0/2}^g(f(x,\tau))\subset
B^g_{r_0}(p_{i_0})$ for some $i_0\in\N$ by
construction, that for all $t\in[\tau,T)$ and $x\in\Sigma$
\begin{eqnarray*}
  f\big(B_{r_0/4}^{\tildeg(t)}(x),t \big)\subset B_{r_0}^g(p_{i_0}).
\end{eqnarray*}
Clearly, using topological data instead of geometric, Lebesgue's number lemma ensures the
existence of such a radius 
(not in terms of $r_0$) as above $r_0/4,$  which would suffice. To lead
(\ref{notpossible}) to a contradiction we may assume that $\tau=0.$\\
Now fix a local chart $\{\varphi:U\to V\subset\R^2\}$ with 
$c^{-1}\d\le(\tildeg_{\!ij}(0))\le c\d$ and diam$_{\tildeg(0)}(U)\le r_0/4.$
Lemma 14.2 in \cite{HAM82} shows that the metrics $\w g(t)$ on $0\le t<T$ are
equivalent, i.e. we have 
$c_{\tildeg}^{-1}\tildeg\!(0)\le\tildeg\!(t)\le c_{\tildeg}\tildeg\!(0)$
and thus
\begin{eqnarray}
  c_{\tildeg}^{-1}\d\le\big(\w g_{ij}(t)\big),\big(\w g^{ij}(t)\big)
    \le c_{\tildeg}\d,\label{512}
\end{eqnarray}
where $c_{\tildeg}=c(n,A_2, R_0,T).$\\
\notiz{das oben brauchen wir jetzt nicht mehr??}
Let $\w \Gamma$ be the associated Christoffel symbols and denote the coordinate
derivative by $\p.$ For any tensor $T\in\Gamma^m(T^{r,s}\Sigma)$ we have the formula
\begin{eqnarray}
  \na^{m}T=\p^m T+\sum\limits^{m}_{l=1}
    \sum\limits_{\;k_1+...+k_l+k\le m-1}\p^{k_1}\w\gg\cdot\ldots\cdot
    \p^{k_l}\w\gg\cdot\p^k T,\label{514}
\end{eqnarray}
where here and in what follows a product as in (\ref{514}) comprises of universal
linear combinations of coordinate functions as $\p^{k_j}\w\gg,~\p^kT$ etc.
The above formula is immediate for $m=1$ and follows then by induction. 
Therefore letting $\w\gg_m:=\bet{\w\gg}+\ldots+\bet{\p^m\w\gg}$ we have
$$\bet{\p^m T}\le c(n,m,\wgm{m-1})
\Big( \bet{\na^{m}T}+\bet{\p^{m-1}T}+\ldots+\bet T \Big)$$ and hence by induction
\begin{eqnarray}
  \bet{\p^{m} T}&\le& c(n,m,\wgm{m-1})\Big( \bet{\na^{m}T}+\bet{\na^{m-1}T}
     +\ldots+\bet T \Big)\no\\
  &\stackrel{(\ref{512})}{\le}&c(n,m,\wgm{m-1},c_{\tildeg})
    \Big(\norm{\na^{m}T}{}\infty{}+\ldots+\norm{T}{}\infty{} \Big).\label{515}
\end{eqnarray}
\begin{lemma}\label{gammadotlemma}
  We have
  \begin{eqnarray}
  \bbet{\p^{m}(\partial_{t}\w\gg)},~\bbet{\p^m\w\gg}
    \le c(n,m,A_{m+3}, R_{m+1},c_{\tildeg},
      \wgm{m}\strichklein{t=0}\quad,T).\label{516}
\end{eqnarray}
\end{lemma}
\bb (Induction over $m\in\N_0$). For $m=0$ we get for the coordinate functions
\begin{eqnarray*}
  \bbet{\partial_t\w\gg}=\bbet{\partial_t\na}
    &\stackrel{(\ref{512})}{\le}& c(c_{\tildeg})
      \norm{\partial_t\na}{}{\infty}{}\\
    &\stackrel{(\ref{510})}{\le}&c(n,A_3, R_1,c_{\tildeg})
\end{eqnarray*}
and hence by integration 
$\bbet{\w\gg}\le c(n,A_3, R_1,c_{\tildeg},\wgm{}\strichklein{t=0}\quad,T).$
For the induction step we obtain from (\ref{515}) and (\ref{510})
\begin{eqnarray*}
  \bbet{\p^{m+1}(\partial_{t}\w\gg)}&\le& c(n,m+1,\wgm{m},c_{\tildeg})
      \big( \norm {\na^{m+1}(\partial_{t}\na)}{}\infty{}
      +\ldots+\norm {\partial_{t}\na}{}\infty{} \big)\\
   &\le&c(n,m+1,\wgm{m},c_{\tildeg},A_{m+4}, R_{m+2})\\
   &\le&c(n,m+1,A_{m+4}, R_{m+2},c_{\tildeg},\wgm{m}\strichklein{t=0}\quad,T).
\end{eqnarray*}
Integrating over $[0,T)$ yields the claim.\eb\\
Summarizing, we now have
\begin{enumerate}
\item [i)] $\frac12\delta\le (g_{\a\b})\le 2\delta\qquad$and
   $\qquad c_{\tildeg}^{-1}\d\le(\tildeg_{\!ij})\strichklein{t=0}\quad
   \le c_{\tildeg}\d$
\item [ii)] $\bet{\p^\gamma \Gamma}\le\Gamma_K:=c(n,K, R_K)\quad$for $\g\le K,$
   where $K$ is as in (\ref{riemannbound})
\item [iii)] $\sum_{l=0}^{L}
  \norm{\widetilde\na^l\!\widetilde A}{}
       {L^\infty(\widetilde G(t),U)}{}
  \le A_L:=c(L,n,t^+,\lam{L+3},\a_0(L+2),(\ref{cgamma}))$
\item [iv)] $\bet{\partial_t\widetilde G}_{\w G(t)}\le P
  :=c(n,A_2, R_0)\quad\forall t\in [0,T)$
\item [v)] $\bet{\partial^\gamma\widetilde\Gamma}\le\wgm N
  :=c(n,N,A_{N+3}, R_{N+1},c_{\tildeg},\wgm{N}\strichklein{t=0}\quad,T)\quad$for 
  $\bet\g\le N$ and all $t\in[0,T).$
\end{enumerate}
Since $\norm{\p_tf}{}{\infty}{}\le c(n,A_2,R_0)
\le c(n,t^+,\Lambda_5,\a_0(4),(\ref{cgamma}))$ we may assume that
$f(\Sigma,t)\subset\!\subset B_R$ for all $t\in[0,T]$ and then that $(M,g)$ is of bounded
geometry of infinite order.
Expressing the $L^2-$gradient of the Willmore functional $\bf W\rm$ in local
coordinates, it is easy to see that $\bf W\rm$ is a universal linear combination of
elements in $\{\na^2\!A,\tildeg,\tildeg^{-1},\p f,g\cf,R\cf \}.$ From
\cref{zeitableitungen} with $m=4,~s=0$ it follows for $p\le8,~q\le 5,~1\le l\le2$
and fixing $K=14$ that
\begin{eqnarray*}
  \bet{\p^{p+1}\!f}\le c(\wgm{7},A_7,\Gamma_7)\mbox{~~~~and~~~~}
  \bet{\p^q\p_t^lf}\le c(\wgm{8},A_{11},\Gamma_{8},R_{9},c_{\tildeg},t^+).
\end{eqnarray*}
Substituting the bounds for the second fundamental form finally yields
\begin{eqnarray}
  \bet{\p^{p+1}\!f},~\bet{\p^q\p_t^lf}\le c(n,t^+,\lam{14},\a_0(13),
  \wgm 8\strichklein{t=0}\quad),\label{coordbound2}
\end{eqnarray}
where c may also depend on the constants in (\ref{cgamma}).\\\\
Now since $\p^{p+1}\!f$ is uniformly Lipschitz continuous in the space variable
for $p\le7,$ we get that the one-parameter family
$\{\p^{p+1}\!f \}_{t\in[t^+-\d,t^+)}$ is equicontinuous and thus it follows from the
Arzel\`a-Ascoli that
$\lim_{t\nearrow t^+}\p^{p+1}\!f=\p^{p+1}\!f_{t^+}$ uniformly for
$p\le 7.$ With the same argument and additionally using that
$\p^q\p_t f$ is also uniformly Lipschitz continuous in the time-variable $(q\le 4)$,
we moreover obtain the uniform convergence
$\lim_{t\nearrow t^+}\p^q\p_t f=\p^q\p_t f_{t^+}$ for $q\le 4.$
\notiz{nochmal checken mit l=2 benoetigt? und wenn ja wo?}
Thus we can define $$\widehat f(\,\cdot\,,t):=
\left\{\begin{array}{l@{\quad for \quad}l}f(\,\cdot\,,t) & 0
\le t<t^+ \\f_{t^+} &t=t^+\end{array}\right.$$
where $\widehat f\in C^{8,0}\cap C^{4+1,1+1}$ on $\Sigma\times [t^+-\d,t^+].$ Because
$f_{t^+}\in C^{4+1}(\Sigma)\subset C^{4+\a}(\Sigma)$ we apply short-time existence
obtaining a Willmore flow $h\in C^{4,1,\a}(\Sigma\times [t^+,\tau))
\cap C^\infty(\Sigma\times (t^+,\tau)).$ We now define
$$E(\,\cdot\,,t):=
\left\{\begin{array}{l@{\quad for \quad}l}f(\,\cdot\,,t) & 0
\le t<t^+ \\f_{t^+} &t=t^+\\ h(\,\cdot\,,t)&t^+<t\le t^++\d<\tau.\end{array}\right.$$
From the discussion above it is clear that $E\in C^{4,0}(\Sigma\times[0,t^++\d],M).$
To see that furthermore $E\in C^{4,1},$ it is enough to check time-differentiability in
$t^+.$ Taking $\widehat f\in C^{8,0}$ into account, we see that
$$\lim_{t\nearrow t^+}\p^q\p_t E=\lim_{t\nearrow t^+}\p^q\p_t\widehat f
=\lim_{t\nearrow t^+}\p^q\bf W\it(\widehat f)
=\p^q\bf W\it(f_{t^+})=\lim_{t\searrow t^+}\p^q\p_th=\lim_{t\searrow t^+}\p^q\p_tE,$$
where we used that $h$ is a Willmore flow defined on $[t^+,\tau).$ It remains to check
that $E$ is parabolically H\"older continuous. We remark that a function is
parabolically $\a$-H\"older continuous, if it is $\a$-H\"older continuous in space and
$\frac \a4$-H\"older continuous in time. Thus $\widehat f$ is parabolically $\a$-
H\"older continuous and, by the regularity properties obtained from short-time
existence, this also holds true for h restricted to $[t^+,t^++\d].$ Using the definition
of parabolic H\"older continuity together with a simple triangle argument, it is not hard
to show that also $E\in C^{4,1,\a}(\Sigma\times [0,t^++\d],M).$ Since we have also shown
that E is a Willmore flow, in particular in $t^+,$ and recalling $t^+=T,$ this obviously
contradicts the maximality of T and hence assumption (\ref{notpossible}) is wrong. \\

We now want to show the inequality on the right-hand side of (\ref{lifespan2}). Therefore
assume that $t^+<T\le\infty.$ Then 
\begin{eqnarray}
  \chi(1,t^+)\geq\e^2_0,\label{nachunten}
\end{eqnarray}
since when assuming the contrary, there exists by maximality of $t^+$ a sequence
$t_i\searrow t^+$ with $\chi(1,t_i)\geq \e_0^2,$ contradicting 
upper semicontinuity of the concentration function $\chi(1,\:\c\:).$
From \cref{452} we know that for any $p\in M$
we have a cutoff function
$\w\g\in C^2(M)$ with \mbox{$\Norm{\w\g}_{C^2(M)}\le C(n)$} and 
$\chi_{{B_{1/2}(p)}}\le\w\g\le\chi_{{B_{1}(p)}}.$
As $\overline{B_1(p)}$ can be covered by a number $\Gamma=\Gamma(n)$
of balls with radius $1/2$ by \cref{hilfslemma51}, we obtain from \cref{prop44} with
$\e_0^2:=\e^*$
for $0\le t\le t^+$
\begin{eqnarray*}
  \Gamma^{-1}\chi(1,t)&\le&\chi(1/2,t)\le\chi(1,0)
    +c\int\limits_0^t\chi(1,s)ds+c\norm{DR}{2}{L^\infty(M)}{}
    \sup\limits_{0\le t\le t^+}\bet{\Sigma}t\\
  &\le&\chi(1,0)+c\int\limits_0^t\chi(1,s)ds
    +\w ct, 
\end{eqnarray*}
where we have used  \cref{lemma46} (note that 
$C_{scal}\stackrel{(\ref{bg2small})}{\le} c(n))$, Cauchy's
inequality, and abbreviated
$\w c:=\norm{DR}{2}{L^\infty(M)}{}(\bet{\Sigma}_{t=0}+\sqrt{2t^+}\W(f_0)).$ From
Gronwall's inequality we infer that
\begin{eqnarray*}
  \chi(1,t)\le(\Gamma\chi(1,0)+\w c\,)\exp(c\Gamma t).
\end{eqnarray*}
Using this inequality for $t:=t^+$ together with (\ref{nachunten}) and estimating
$$\bet{\Sigma}_{t=0}+\sqrt{2t^+}\W(f_0)
\le 2\exp(ct^+)(\bet{\Sigma}_{t=0}+\W(f_0))$$
we obtain $$\e_0^2\le \big(\Gamma\chi(1,0)+
\norm{DR}{2}{L^\infty(M)}{}(\bet{\Sigma}_{t=0}+\W(f_0))\big)\exp\big(c(\Gamma)t^+\big)$$
and arrive at
$$t^+\geq c(\Gamma)\log\bigg(\frac{\e_0^2}
{\Gamma\chi(1,0)+
\norm{DR}{2}{L^\infty(M)}{}(\bet{\Sigma}_{t=0}+\W(f_0))}\bigg).$$
\cref{lifespan} now follows after rescaling.
\setcounter{section}{0}
\setcounter{equation}{0}
\chapter{Blow-up of singularities}
\notiz{check einleitung}
In \cite{Blatt} it was shown that singularities may occur either in finite time or at
infinity if the
ambient manifold is $(\R^n,\d_{eucl}).$ To study assumed singularities in the
aforementioned cases if the ambient manifold is more generally
a (possibly noncompact) complete Riemannian manifold $(M^n,g)$ of bounded geometry we
want to perform a blow-up procedure as in \cite{KS01}. Namely, we show show that the
blow-up is an immersed (time independent) Willmore surface that is either an embedded
round sphere, or contains at least one component which is
a nonumbilic (compact or noncompact) Willmore surface. To do so, it is essential to
have a mass-density estimate and interior estimates to ensure compactness, i.e.
$C^{k}-$subconvergence for the blow-up sequence using a reparametrization.
The mass-density estimate and the interior estimates will be provided in the forthcoming
sections.
\section{Monotonicity formulas for Riemannian Manifolds}
\label{blow-up}
In the following lemma, it is again important to have ``good'' coordinates for $M$. The
proof of the Michael-Simon Sobolev inequality (\ref{mssieqn}) suggests to make use of
harmonic coordinates since the bounds on the curvature of $M$ are rather week (cf.
\cref{harmcoord1}).
Unfortunately, we need to have Christoffel symbols growing at most linearly.
Although there exist estimates (cf. \cite{jostkarcher}) for the Christoffel
symbols in harmonic coordinates such that they are bounded by the squared
radius of the geodesic ball they are defined on, we only have these bounds for a
fixed chart. Also, at some point we make use of Gauss's Lemma so that we choose
Riemannian normal coordinates this time, under the drawback that the geometry of $(M,g)$
has to be bounded of order one. \notiz{$\Sigma_R$ definieren}
\begin{lemma}{\rm(Monotonicity formula for Riemannian
  manifolds).}\label{kleineradien}~Let
  $(\Sigma^2,\tildeg)$ and $(M^n,g)$ be complete (possibly noncompact) Riemannian
  manifolds 
  and $f:(\Sigma,\tildeg)\rightarrow(M,g)$ be a proper isometric $C^{2}-$immersion. If
  $r_0>0$ can be chosen such that $r_0\Lambda<c(n),$ where 
  $\Lambda:=\norm{R}{1/2}{\infty}{}+\norm{DR}{1/3}{\infty}{}+\inj(M,g)^{-1}$ and
  $c=c(n)$  is a small universal
  constant, then
  \begin{eqnarray}
  \dfrac{\lvert\Sigma_{\sigma}\rvert}{\sigma^2}
    \le C\left(\dfrac{\lvert\Sigma_{\varrho}\rvert}{\varrho^2}
    +\W(\Sigma_{\varrho})\right)\label{monformula}
  \end{eqnarray}
  for all $0<\sigma\le\varrho\le r_0,$ where $C=C(n)$ only depends on the dimension n. 
\end{lemma}~\\
\notiz{checke, ob wir die ganzen nummern brauchen...}
\bf{Remark:}\rm~To see that in general (\ref{monformula}) cannot hold for all
$0<\sigma\le\varrho<\infty,$ at least in case $\Sigma$ and $M$ are closed, again take the
standard sphere 
$\mathbb S^2\subset\mathbb S^n\subset\R^{n+1}$ as described in the remark of
\cref{mssi}. Since $\W(\mathbb S^2\subset\mathbb S^n)=0,$ letting 
$\varrho\nearrow\infty$ would imply that $\bet{\Sigma_\sigma(p)}=0$ for any
$\sigma$ and any $p\in M.$\er\\ 
\bb
  The proof imitates the proof of a monotonicity formula in case
  $(M,g)=(\R^n,\d)$ 
  (see \cite{SIM93}) using Riemannian normal coordinates locally. One has to 
  find an appropriate test vector field in the first variation formula 
  $$\is div_\Sigma \phi\,d\mu=-\is g\cf(\phi\cf, H) d\mu,$$ where $div_\Sigma\phi
    =g\cf(D\phi\cf\cdot Df\cdot \tau_i,Df\cdot \tau_i),$
  $\phi\in C^{0,1}_c(M,TM)$ is a Lipschitz vector field with compact support, $H$ is
the
  mean curvature of f, $\{\tau_i\}_{i=1,2}$ is a locally defined 
  $\tildeg-$orthonormal basis and summation over repeated indices is used.
  
  \quad For this, let $p\in M$ be arbitrary but fixed, $\sigma,~\varrho$ as in the
statement
  and fix Riemannian normal coordinates 
  $\{y^\a \}=\varphi_p:B_{r_0}(p)\rightarrow B_{r_0}^n\subset\mathbb R^n,$
  with respect to $g$  centred at p. From \cref{riemcoord1} we know that
  \begin{equation*}
    \begin{array}{rlll}
      i)&(1-C r^2\Lambda^2)\delta
        \le (g_{\a\b})\le(1+C r^2\Lambda^2)\delta&&\mbox{on $B_r^g(p)$ for all } 
        0\le r<r_0\label{metabs}\\
      ii)&\sup_{B_r^g(p)}\bet{\Gamma_{\b\d}^\a }\le C r\Lambda^2&&
        \mbox{for all } 0\le r<r_0.
    \end{array}
  \end{equation*}
  For $x\in B_{r_0}^n$ define the vector field 
  $\bar\phi\in C^{0,1}_c(B_{r_0}^n,\R^n)$ by 
  $$\bar\phi (x):=\left(\dfrac{1}{\max(\lvert x\rvert,\sigma)^2}
    -\dfrac{1}{\varrho^2}\right)_+ x$$ and $$\phi:=(\varphi_p)^*\bar \phi.$$

  To simplify notation we may by locality assume that 
  $\Sigma\hookrightarrow M$ is embedded. We compute for $\{\tau_i\}$ as above
  $$D\phi\c\tau_i=\tau_i^\a D\phi\c\p_\a=\tau_i^\a D(\phi^\beta\pb)\c\p_\a
  =\tau_i^\a(\p_\a \phi^\beta)\pb+\tau_i^\a\phi^\beta\Gamma^\gamma_{\alpha\beta}\pg.$$
  Let $P=P^\top_g$ be the $g$-orthonormal projection onto the tangent bundle $T\Sigma,$ 
  in coordinates
\begin{eqnarray*}\label{defproj}
    P_\alpha^\beta&=&g(P(\pa),\pg)g^{\beta\gamma}\\
                  &=&g(\tau_i,\pa)g(\tau_i,\pg)g^{\beta\gamma},
\end{eqnarray*}
  and define $$\sigsig:=\{q\in\Sigma\cap B_{r_0}:
  \bet{q}<\sigma\}$$ and the annulus
  $$\sigsigrho:=\{q\in\Sigma\cap B_{r_0}:0<\sigma\le\bet{q}<\varrho\},$$
  where here and in what follows we identify $(B_{r_0},g)$ with 
  $(B_{r_0}^n,(g_{\a\b}))$ via the coordinates $\varphi_p$.\\
  First, assume that $q\in\sigsig.$ We get since 
    $\phi(q)=(\sigma^{-2}-\varrho^{-2})q^\a\pa$\\
  \begin{eqnarray*}
    D\phi\c\tau_i(q)
    &=&(\sigma^{-2}-\varrho^{-2})
       \left[  g^{\beta\gamma}  g(\tau_i,\pb)\delta^\alpha_\gamma\pa
       +  g^{\beta\g}  g(\tau_i,\pb)q^\delta\Gamma^\alpha_{\g\d}\pa\right]\\
    &=&(\sigma^{-2}-\varrho^{-2})\left[\tau_i
       +\tau_i^\gamma q^\delta\Gamma_{\gamma\delta}^\alpha\pa\right]
  \end{eqnarray*}
  which yields, using $\tau_i^\g=g(\tau_i,\p_\varrho) g^{\varrho\g},$
  $$g(D\phi\c\tau_i,\tau_i)=(\sigma^{-2}
  -\varrho^{-2})(2+P^\gamma_\alpha q^\delta\Gamma^\alpha_{\gamma\delta}).$$
  Therefore, abbreviating 
  $f^\alpha\Gamma^\gamma_{\alpha\beta}\cf P^\beta_\gamma\cf=:f\,\Gamma P$
  we generally have for an immersion f
  \begin{eqnarray}\label{divergenzeins}
      div_\Sigma\phi=(\sigma^{-2}-\varrho^{-2})(2+f\,\Gamma P).
  \end{eqnarray}
  Now assume that $q\in\Sigma_{\sigma,\varrho}.$
  Since $\phi(q)=(\lvert q\rvert^{-2}-\varrho^{-2})q^\alpha\pa$ and
  $\pg \bet{q}^{-2}=-2 q^\rho\d_{\rho\g}\betrag q {-4}$ it follows with
  $X(q):=q^\alpha\pa$
  $$D\phi\c\tau_i=\tau_i^\gamma\left(\betrag q {-2}
    -\varrho^{-2}\right)\pg-2q^\varrho\d_{\varrho\g}\betrag q {-4}\tau_i^\gamma X
    +\tau_i^\gamma(\betrag q {-2}-\varrho^{-2})q^\alpha\Gamma_{\gamma\alpha}^\delta\pd.$$
  Due to the existence of two Riemannian metrics $g$ and $\d,$ we need to be careful,
  when lowering indices as in $q^\varrho\d_{\varrho\g}$.
  Since $X$ is a radial vector field, we infer from Gauss's
  Lemma that $$\d(P^\top_{g}\!X,X)=g(P^\top_g\!X,X)=\bet{P^\top_g\!X}^2_g.$$
  Thus, using $q^\r\d_{\r\g}\tau_i^\g g(X,\tau_i)=\d(X,\tau_i)g(X,\tau_i)
    =\d(P^\top_{g}\!X,X),$ we now have 
  $$g(D\phi\c\tau_i,\tau_i)\strichklein{$\,q$}
    =(\betrag q {-2}-\varrho^{-2})(2+q^\alpha\Gamma_{\alpha\gamma}^\delta P^\gamma_\d)
    -2\betrag q {-4}\bet{P^\top_g\!X}^2_g,$$
  or generally,
  \begin{eqnarray}\label{divergenzzwei}
     div_\Sigma\phi
      =(\betrag f {-2}-\varrho^{-2})(2+f\,\Gamma P)-2
        \betrag{P_g^\top\! X}{2}_{g}\cf\betrag f {-4}.
  \end{eqnarray}
  In the sequel, we do not distinguish between $X$ and  $X\cf.$
  By (\ref{divergenzeins}) and (\ref{divergenzzwei}) we get
  \begin{eqnarray}
    \frac 12\;\isr div_\Sigma \phi
      &=&\sigma^{-2}\betrag\sigsig{}-\varrho^{-2}\betrag\sigrho{}
        +\frac{1}{2\sigma^2}\iss f\,\Gamma P
        -\frac{1}{2\varrho^2}\isr f\,\Gamma P\no\\
      &&+\issr\betrag f {-2}-\issr\dfrac{\betrag {P_g^\top\! X}{2}_{ g }}{\betrag f 4}
        +\frac{1}{2}\issr\betrag f{-2} f\,\Gamma P.\label{mittelrekruemmungeins}
  \end{eqnarray}
  Since $\phi(q)=\left(\max(\bet{q},\sigma)^{-2}-\varrho^{-2}\right)_+X(q)$ it follows,
  abbreviating $g\cf(Y,Z)$ with $Y\c Z,$
  \begin{eqnarray}\label{mittlerekruemmungzwei}
    -\isr\phi\c H&=&-\iss\sigma^{-2}X\c H-\issr\betrag f{-2} X\c H
    +\isr\varrho^{-2}X\c H. 
  \end{eqnarray}
  Completing the square yields 
  \begin{eqnarray*}
    \sigma^{-2}\betrag\sigsig{}&+&\issr\left\lvert\frac 1 4 H
       +\dfrac{P^\bot_g\!X}{\betrag{X}{2}_\delta}\right\rvert^2_{ g }\\
    &\le&\varrho^{-2}\betrag{\sigrho}{}+1/8\,\W (\sigrho)+\frac 1 {2\varrho^2}\isr X\c H
       -\frac 1 {2\sigma^2}\iss X\c H-\frac 1 {2\sigma^2} \iss f\,\Gamma P\\
    && -\frac 1 2 \issr\betrag X {-2}_\delta f\,\Gamma P
       +\frac 1 {2\varrho^2}\isr f\,\Gamma P
       +\issr\betrag X {-2}_\delta\left(\dfrac{\betrag X2_g}{\betrag X
         2_\delta}-1\right).
  \end{eqnarray*}
  Because $\betrag {X\c H}{}\le c\bet{X}_\delta\betrag{H}{}_g\le cr\betrag H{}_g$ on
  $\Sigma_r$ for any $r\le r_0,$ we get after absorption and dropping the square terms on
  the left-hand side
  \begin{eqnarray*}
    \dfrac{\betrag\sigsig{}}{\sigma^2}
    &\le&c(n)\bigg(\dfrac{\betrag{\sigrho}{2}}{\varrho^2}
      +\W (\sigrho)\bigg)-\frac 12\,(\sigma^{-2}\kreis 1-\varrho^{-2}\kreis 2)
        \iss f\,\Gamma P\\
    &&+\issr\betrag f {-2}\bigg(\dfrac{ g _{\alpha\beta}\cf f^\alpha f^\beta}
      {\betrag f 2} 
      \put(-11,22.6){\circle{8}}\put(-13,20){\scriptsize{3}} -1\bigg)
      -\frac 12\issr\big(\betrag f{-2}\kreis 4-\varrho^{-2}\kreis 5\big)f\,\Gamma P.
  \end{eqnarray*}
  In what follows, we estimate the integrals on the right-hand side.\\\\
\underline{Ad 1:} To begin with, we trivially have $\lvert f^\a\rvert\le\sigma$ on 
$\sigsig$ by definition of $\sigsig.$ The projection $P=P^\top_g$ is, as expected, also
estimated:
$\sum_\b\bet{P^\b_\a}^2=\bet{P\p_a}_\d^2\le c\bet{P\p_a}_g^2\le c\bet{\p_\a}_g^2
  =c g_{\a\a}\le c(n).$
Since $\bet{\Gamma}_\d\le c\sigma\Lambda^2$ on $B_\sigma(p)$ we get
$$\betrag{f\,\Gamma P}{} \le c(n)\sigma^2\Lambda^{2}\le c(n)r_0^2\Lambda^2\le 1/2,$$
choosing $r_0$ smaller if necessary. Therefore 
\begin{eqnarray}
  \bigg\lvert\dfrac{1}{2\sigma^2}\iss f\,\Gamma P\bigg\rvert
    &\le&\dfrac{\betrag{\sigsig}{}}{4\sigma^2}.
\end{eqnarray}\\
\underline{Ad 2:} Since $\varrho^{-2}\le\sigma^{-2}$ this integral is treated as the
latter.\\\\
\underline{Ad 3:} If one again uses Gauss's Lemma, it is immediate that this integral
vanishes. Namely, it is $\bet{X}^2_g=\bet{X}^2_\d$ because $\exp$ is a radial
isometry. If one had other coordinates, only using that the eigenvalues of 
$(g_{\a\b})$ grow at most quadratically (as it is the case here), one could
alternatively argue as follows without using Gauss's Lemma: Taking
\begin{eqnarray}\label{ew1}
  -c(n)\Lambda^{2} r^{2}\d\le ( g _{\a\b}-\d_{\a\b})\le c(n)\Lambda^{2} r^{2}\d
\end{eqnarray}
on $B_r(p)$ into account, we get for any $\sigma\le s\le r\le\varrho$
\begin{eqnarray*}
  \bigg\lvert\int\limits_{\Sigma_{s,r}}\betrag f {-2}
      \Big(\dfrac{ g_{\a\b}f^\a f^\b}{\betrag f 2}-1\Big)\bigg\rvert
    &\le&c(n)\dfrac{r^2}{s^2}\betrag{\Sigma_{s,r}}{}\Lambda^{2}.
\end{eqnarray*}
Restricting to special annuli, i.e. annuli with ratio $r/s\le2$, we can further estimate
\begin{eqnarray}
  &\le&4c(n)\betrag{\Sigma^{}_{s,r}}{}\Lambda^{2}\no\\
  &\le&\dfrac{\betrag{\Sigma^{}_{s,r}}{}}{r_0^2}\label{annulusabschaetzung}
\end{eqnarray}
since $r_0\Lambda\le(4c(n))^{-1/2}$ choosing $r_0$ smaller if
necessary. Now choose a geometric decomposition of 
$\sigsigrho$ in annuli $\Sigma_{\sigma_i,\sigma_{i+1}},$ such that
$\dfrac{\sigma_{i+1}}{\sigma_i}<2.$ From (\ref{annulusabschaetzung}) it   follows with
$N=N(\sigma,\varrho)<\infty$ that
\begin{eqnarray*}
  \bigg\lvert\issr\betrag f {-2}
    \Big(\dfrac{ g_{\a\b}f^\a f^\b}{\betrag f 2}-1\Big)\bigg\rvert
  &\le&\dfrac{1}{r_0^2}\sum^{N}_{i=1}
    \betrag{\Sigma_{\sigma_i,\sigma_{i+1}}}{}   
    =\dfrac{\betrag{\sigsigrho}{}}{\varrho^2}
    \le\dfrac{\betrag{\sigrho}{}}{\varrho^2}.
\end{eqnarray*}
\underline{Ad 4:} Analogously to the third integral, again using a decomposition, we can
estimate
\begin{eqnarray*}
  \int\limits_{\Sigma_{\sigma_i,\sigma_{i+1}}}
      \betrag f {-2}\betrag{f\,\Gamma P}{}
  &\le& c(n)\frac{1}{\sigma_i^2}
      \betrag{\Sigma_{\sigma_i,\sigma_{i+1}}}{}\sigma^2_{i+1}\Lambda^2
  \le \dfrac{\lvert\Sigma_{\sigma_i,\sigma_{i+1}}\rvert}{\varrho^2}.
\end{eqnarray*}
\underline{Ad 5:} Similarly, 
\begin{eqnarray*}
  \dfrac{1}{\varrho^2}\issr\betrag{f\,\Gamma P}{}
  &\le&c(n)\Lambda^2\betrag{\sigsigrho}{}\le \dfrac{\betrag{\sigrho}{}}{\varrho^2}.
\end{eqnarray*}
Combining, we finally arrive at 
$\sigma^{-2}\bet{\Sigma_{\sigma}}\le c(n)\big(\r^{-2}\bet{\Sigma_{\r}}
  +\W(\Sigma_{\r}) \big).$\eb\\
Unfortunately, the monotonicity formula of \cref{kleineradien} only holds for small
radii, with smallness depending on the geometry of $(M,g).$ For large radii we are
forced to make additional assumptions to prove the following Corollary.
\begin{korollar}{\rm(Mass-density estimate in the large by means of a mass
bound).}\label{massbound3}\\
  Let $f:\Sigma\times[0,T)\to(M,g)$ be a Willmore flow of a closed surface into a
  Riemannian manifold of bounded geometry. Assume that 
  $\mathfrak M_{f}:=\sup_{t\in[0,T)}\mue{f_t}\:(\Sigma)<\infty.$
  Then for all $R>0$
  \begin{eqnarray}
    \frac{\bet{\Sigma_R}}{R^2}\le c\,\mathfrak M_f
    \big(\norm{R}{}{\infty}{}+\norm{DR}{2/3}{\infty}{}+\inj(M,g)^{-2}\big)+c\W(f_{t=0}).
    \label{massdensity1}
  \end{eqnarray}
  In case $T<\infty,$ we have 
  \begin{eqnarray}
    \mathfrak M_f\le\sqrt{2T}\W(f_0)+\mue{f_0}\:(\Sigma)<\infty.\label{TFiniteBound}
  \end{eqnarray}
  If the sectional
  curvature $K^M$ of $(M,g)$ is uniformly negative, i.e. if 
  $K^M\le-\hat\kappa^2<0,$ then it is
  \begin{eqnarray}
    \mathfrak M_f\le\frac{\W(f_0)-4\pi\chi(\Sigma)}{2\hat\kappa^2}.\label{kappabound}
  \end{eqnarray}
  (\ref{TFiniteBound}) and (\ref{kappabound}) hold for all smooth Riemannian manifolds.
\end{korollar}
\bb
We clearly have $\W(\Sigma_\varrho)\le \W(f_t)\le\W(f_0)$ by definition of the
$L^2-$gradient flow
(see the proof of \cref{lemma46}). The claim now follows easily from
\cref{kleineradien},~\cref{lemma46} and (\ref{nonequivfunc}).\eb\\
\section{Interior estimates  II}
Although we already have proven interior estimates for the second fundamental form in
$C^k$ by \cref{prop47} we cannot use them for the blow-up at infinity owing to their
dependence on the curvature of the initial surface. Thus, in case of singularities at
infinity, we need \cref{interiorestimates} as a second version.
Compared to \cite{KS01},
Theorem 3.5, our estimates are slightly worse, and more complicated to prove. First, we
were not able to prove that the curvature in (\ref{intestimate}) grows at most as fast as
the square root of the local energy concentration (\ref{enconcent}). This is mainly due
to the existence of the
volume-term in (\ref{434}) of the integral estimates that we cannot control. For the
same reason, we only have estimates depending on initial conditions. Namely, we have
to assume  (\ref{initialcond}). Second, the proof is more complicated since we have to
localize in time twice. This is because of the existence of higher homogeneity
curvature terms in (\ref{434})
that itself results from the interpolation of the curvature terms $Q^{m+2,1}$ in 
(\ref{hiergibtsQ}).
\notiz{noch durchsehen auf die bedingungen an lambda k usw}
\notiz{das c oben k unten l entfernen}
\begin{lemma}{\rm (Interior estimates II).}\label{interiorestimates}~Let
  $(M,g)$ be a Riemannian manifold of bounded geometry and $p\in M.$ Let $\r>0$ be
  chosen such that
  \begin{eqnarray*}
    \varrho<\min\Big\{\frac\pi{2\kappa},\inj_{(M,g)}\Big\}.
  \end{eqnarray*}
  Let $f:\Sigma\times[0,T]\to (M,g)$ be a Willmore flow of a closed surface satisfying 
  $T\le C(n)\r^4$ and
  \begin{eqnarray}
    \sup\limits_{0\le t\le T}\int\limits_{\Sigma_{\r}(p)}\bet{A}^2\mm
      \le\e_0(n)\label{enconcent}
  \end{eqnarray}
  for small universal $\e_0>0,$ where $\Sigma^g_\r(p):=f^{-1}\big(B^g_\r(p) \big).$ 
  Moreover, choose $\Lambda_K<\infty$ such that
  \begin{eqnarray}
    \varrho(\varrho^2+\,\mathfrak M_f)\sum\limits_{i=1}^{K}\norm{D^i\!R^{M\!}}
      {\frac{3}{i+2}}{L^\infty(M,g)}{}
      <\lam K.\label{initialcond}
  \end{eqnarray}
  for some $K\in\N.$ Then we have for all $t\in (0,T]$
  \begin{equation}
    \norm{\na^{k}\!A}{}{L^{2}(\Sigma_{\varrho/2})}{}
    \le c\big(n,k,C,\lam{k+1}\big)\,t^{-\frac{k}{4} }\label{intestimate}
  \end{equation}
  and
  \begin{equation*}
    \norm{\na^{k}\!A}{}{L^{\infty}(\Sigma_{\varrho/2})}{}
    \le c\big(n,k,C,\lam{k+3}\big)\,t^{-\frac{k+1}{4}}.
  \end{equation*}
\end{lemma}
\bf{Remark:}\rm~Note that $\Lambda_K<\infty$ can always be chosen since
$\mathfrak M_f\le\sqrt{2T}\W(f_0)+\mue{f_0}\:(\Sigma)<\infty$ by (\ref{TFiniteBound}).
\er\\
\bb
After scaling, we may assume that $\rho=1.$ Let $c_k:=c(\Lambda_k).$
Using \cref{452} we choose a cutoff function
$\w\g\in C^2(M)$ with $\absch{63/64}~\le\w\g\le\absch{1}$\!\!\!and
$\norm{D^i\w\g}{}{\infty}{}\le c(n)~(i=1,2).$
From \cref{prop44} we then get 
\begin{eqnarray}
  \int\limits^T_0\int\limits_{\Sigma_{63/64}}\bet{\na^2\!A}^2+\bet A^6\,d\mu\,dt
    \le\e_0+c_1\e_0 T+c_1 T\le c_1,\label{53ia}
\end{eqnarray}
where we used
\begin{eqnarray}
  \norm{DR}{2}{\infty}{}\sup\limits_{0\le t\le T}
    \bet{\Sigma}^{}\le
    c_1\norm{DR}{}{\infty}{}\mathfrak M_f\le c_1.\label{similar}
\end{eqnarray}
Using (\ref{423}) from \cref{lemma43} we get for $\w\g$ as above
with $\absch{31/32}~\le\w\g
\le\absch{63/64}$
\begin{eqnarray}
  \int\limits_0^T\norm{A}{4}{\infty,\Sigma_{31/32}}{}dt
    \le c\e_0\c\big(c_1+c\e_0 T\big)
    \le c_1.\label{54ia}
\end{eqnarray}
Now let $\g:=\w\g\cf$ for a cutoff function with 
$\absch{15/16}\le\w\g\le\absch{31/32}\quad.$
Furthermore, define a cutoff function in time by
\begin{eqnarray*}
  \psi(t):=\frac{t}{T}
\end{eqnarray*}
such that
\begin{eqnarray}
  \dot\psi(t)=\frac{1}{T}.\label{55ia}
\end{eqnarray}
Introducing the notation 
$\a_j(t):=c_{2j+1}(1+\norm{A}{4}{\infty,\Sigma_{31/32}}{})$ and 
$F_j(t):=\int_{\Sigma}\bet{\na^{2j}\!\!A}^2\g^{4j+4}\mm$ we obtain from \cref{prop45}
for $j\in\{0,1\}$
\begin{eqnarray}
  \dfrac{d}{dt}F_j(t)+\frac{1}{2}F_{j+1}(t)
    \le \a_j(t)F_j(t)+\a_j(t),
\end{eqnarray}
where we used an estimate similar to (\ref{similar}). Further we compute with 
$h_j(t):=\psi(t)F_j(t)$
\begin{eqnarray}
  \dfrac{d}{dt}h_j(t)+\frac{t}{2T}F_{j+1}(t)\le \a_j(t)h_j(t)+\frac{t}{T}\a_j(t)
  +\frac{1}{T}F_j(t).\label{56ia}
\end{eqnarray}
Now let $\absch{3/4}\le\w\g\le\absch{7/8}$ and $0\le j\le m.$ We define cutoff functions
in time similar to \cite{KS01} by
\begin{eqnarray}
 \chi_j(t):=\left\{\begin{array}{ll} 0& \mbox{for }t\le (j-1)\frac{T}{2m}+\frac{T}{2}\\
   \frac{2m}{T}\big(t-\frac{T}{2}-(j-1)\frac{T}{2m} \big) & \mbox{in between}\\
   1&\mbox{for }t\geq j\frac{T}{2m}+\frac{T}{2},\end{array}\right.\label{chidefi}
\end{eqnarray}
if $j\geq 1,$ and $\chi_{_0}(t):\equiv 1.$ 
Now for $j\geq 0$ define analogously $e_j(t):=\chi_j(t)E_j(t),$ where $E_j$ equals to 
$F_j$ except that we now have $\absch{3/4}\le\w\g\le\absch{7/8}$ and $0\le j\le m.$
Restricting to 
$j\geq 1,$ since $0\le\dot\chi_{j(t)}\le \frac{2m}{T}\chi_{j-1}$ we again get from
\cref{prop45}
\begin{eqnarray}
\dfrac{d}{dt}e_j(t)+\frac{1}{2}\chi_j(t)E_{j+1}(t)
  \le\hat\a_j(t)e_j(t)+\chi_j(t)\hat\b_j(t)+\frac{2m}{T}\chi_{j-1}(t)E_j(t)
\end{eqnarray}
if we now define 
$\hat\a_j(t):=c_{2j+1}\c\big(1+\norm{A}{4}{\infty,\Sigma_{7/8}}{}\big)$ and
$\hat\b_j(t):=c_{2j+1}\c\big(1+\norm{A}{4j}{\infty,\Sigma_{7/8}}{}\big).$
Applying Gronwall's Lemma yields for any $t\in[T/2,T]$
\begin{eqnarray}
  \lefteqn{\hspace{-0em}e_j(t)+\frac{1}{2}\int\limits_{T/2}^t\chi_j(s)E_{j+1}(s)ds}\no\\
  &\le&\exp\bigg(\int\limits_{T/2}^T\hat\a_j(s)ds \bigg)
    \Big[\underbrace{e_j(T/2) }_{=0}+\int\limits^t_{T/2}\chi_j(s)\hat\b_j(s)
    +\frac{2m}{T}\chi_{j-1}(s)E_j(s)ds\Big]\no \\
  &\le&\exp\bigg(c_{2j+1}
    +c_{2j+1}\int\limits^T_{T/2}\norm{A}{4}{\infty,\Sigma_{7/8}}{}ds\bigg)\no\\
  &&\c\Big[c_{2j+1}+c_{2j+1}\int\limits^t_{T/2}
    \norm{A}{4j}{\infty,\Sigma_{7/8}}{}
    +\frac{2m}{T}\chi_{j-1}(s)E_j(s)ds \Big].\label{31ia}
\end{eqnarray}
To be able to estimate $\int_{T/2}^T\norm{A}{4j}{\infty,\Sigma_{7/8}}{}ds,$
we apply Gronwall's Lemma to (\ref{56ia}) and get for $t\in[0,T]$
\begin{eqnarray*}
  \lefteqn{\hspace{-3em}\frac{t}{T}F_1(t)+\frac{1}{2T}\int\limits^t_0sF_2(s)ds}\\
  &=&\frac{t}{T}\is\bet{\na^2\!A}^2\g^8\mm
    +\frac{1}{2T}\int\limits_0^t s\is\bet{\na^{4}\!A}^2\g^{12}\mm\,ds \\
  &\le&\exp\bigg(c_3\c(T+\int\limits_0^T\norm{A}{4}{\infty,\Sigma_{31/32}}{}ds) \bigg)
    \Big[\int\limits_0^T\frac{s}{T}\a_1(s)ds
     +\frac{1}{T}\int\limits^T_0
     \int\limits_{\Sigma_{31/32}}\bet{\na^2\!A}^2\mm\,ds\Big]\\
  &\stackrel{(\ref{54ia})}{\le}&c_3\Big[c_3
   +c_3\int\limits_0^T\norm{A}{4}{\infty,\Sigma_{31/32}}{}ds
   +\frac{1}{T}c_1 \Big],
\end{eqnarray*}
where we used (\ref{53ia}) for the last term. With (\ref{54ia}) and 
$T\le C(n)$ we further obtain
$t F_1(t)\le c_3$
such that we finally get on $[T/2,T]$
\begin{eqnarray*}
  \int\limits_{\Sigma_{15/16}}\bet{\na^2\!A}^2\mm\le\frac {c_3}{T}.
\end{eqnarray*}
Employing \cref{lemma43} for $\absch{7/8}\le\w\g\le\absch{15/16}$\quad yields
\begin{eqnarray}
  \norm{A}{4}{\infty,\Sigma_{7/8}}{}
  \le c\e_0\c\Big(\frac{c_3}{T}+c\e_0 \Big)\le\frac{c_3}{T}\label{41ia}
\end{eqnarray}

for all $t\in[T/2,T].$ Using this, we can go on estimating (\ref{31ia})
\begin{eqnarray}
  \lefteqn{\hspace{-5em}e_j(t)+\frac{1}{2}\int\limits_{T/2}^t\chi_j(s)E_{j+1}(s)ds}\no \\
 &\le&\exp\big(c_{2j+1} \big)\Big[c_{2j+1}+c_{2j+1}\frac{T}{2}\frac{c_{3}}{T^{j}}
    +\frac{2m}{T}\int\limits_{T/2}^t\chi_{j-1}(s)E_j(s)ds \Big]\no \\
 &\le&\frac{c_{2j+1}}{T^{j-1}}
    +c_{2j+1}\frac{2m}{T}
    \int\limits^{t}_{T/2}\chi_{j-1}(s)E_j(s)ds.\label{annahmeia}
\end{eqnarray}
Now we show by induction that for $0\le j\le m$ and all $t\in[T/2,T]$
\begin{eqnarray*}
  e_j(t)+\frac{1}{2}\int\limits_{T/2}^t\chi_j(s)E_{j+1}(s)ds
  \le c_{2j+1}(m)\dfrac{1}{T^{j}}.
\end{eqnarray*}
For $j=0$ we obtain using (\ref{53ia}) (recalling $\absch{3/4}\!\le\w\g\le\absch{7/8}$)
\begin{eqnarray*}
  \is\bet{A}^2\g^4\mm+\frac{1}{2}\int\limits_{T/2}^t\is\bet{\na^2\!A}^2\g^8\mm\le\e
    +c_1\le c_1.
\end{eqnarray*}
For $j\geq 1$ we can estimate, using (\ref{annahmeia}),
\begin{eqnarray*}
  e_{j}(t)+\frac{1}{2}\int\limits_{T/2}^t\chi_{j}(s)E_{j+1}(s)ds
  \le&\dfrac{c_{2j+1}}{T^{j-1}}
    +c_{2j+1}\dfrac{4m}{T}\dfrac{c_{2j-1}(m)}{T^{j-1}}
  \le\dfrac{c_{2j+1}(m)}{T^{j}}.
\end{eqnarray*}
Thus we have at time $t=T$
\begin{eqnarray*}
  \is\bet{\na^{2m}\!A}^2\g^{4m+4}\mm
  \le\dfrac{c_{2m+1}(m)}{T^{m}}.
\end{eqnarray*}
To get an estimates for odd order derivatives one may simply apply
\cref{lemma61} with $r=1,~p=q=2,~\a=1,~\b=0,~s=4m+6$ and 
$t:=s^{-1}\in[-\frac 12,\frac 12]$ to obtain
\begin{eqnarray}
  \is\bet{\na^{2m+1}\!A}^2\g^{4m+6}\mm\le \dfrac{c_{2m+3}(m)}{T^{m+1/2}}.
\end{eqnarray}
In contrast to the time exponent, the argument does not close for the order of the
bounded geometry. More precisely, when using the interpolation inequality we have a loss
of one order with respect to the geometry bounds $c_k=c(\Lambda_k).$ Instead, it is
possible to estimate the odd order derivatives directly analogously to the above
estimates: We let $\absch{3/4}\le\w\g\le\absch{7/8},$
$j\geq1,$ $~k_j(t):=\chi_j(t)K_j(t),$
$\chi_j(t)$ as in (\ref{chidefi}), 
$$K_j:=\int\limits_\Sigma\bet{\na^{2j-1}\!A}^2\g^{4j+2}\mm,~~
\tilde\a_j(t):=c_{2j}\c(1+\norm{A}{4}{\infty,\Sigma_{7/8}}{} )~~
\mbox{and}~~\tilde\b_j(t):=
  c_{2j}\c(1+\norm{A}{\max\{4,4j-2 \}}{\infty,\Sigma_{7/8}}{}).$$
From (\ref{41ia}) we obtain for $j\geq2$
$$c_{2j}\int\limits^{t}_{T/2}\norm{A}{\max\{4,4j-2\}}{\infty,\Sigma_{7/8}}{}ds
\le\frac{c_{2j}}{T^{j-3/2}}$$
and thus
\begin{eqnarray*}
  k_j(t)+\frac{1}{2}\int\limits_{T/2}^t\chi_j(s)K_{j+1}(s)ds
    \le\frac{c_{2j}}{T^{j-3/2}}+\frac{c_{2j}(m)}{T}
    \int\limits^{t}_{T/2}\chi_{j-1}(s)K_j(s)ds.
\end{eqnarray*}
As above, we show by induction that for all $j\geq 1$ and all $t\in[T/2,T]$
\begin{eqnarray}
  k_j(t)+\frac{1}{2}\int\limits_{T/2}^t\chi_j(s)K_{j+1}(s)ds
    \le\frac{c_{2j}(m)}{T^{j-1/2}}.\label{indhypth}
\end{eqnarray}
For $j=1$ we again use \cref{lemma61} with 
$r=1,~p=q=2,~\a=1,~\b=0,~s=4m+6$ and 
$t:=s^{-1}\in[-\frac 12,\frac 12]$
to obtain
\begin{eqnarray*}
  \is\bet{\na\!A}^2\g^6\mm\le c\norm{A}{}{2}{0}\norm{\na^2\!A}{}{2}{0}
    +c\norm{A}{2}{2}{0}+\frac 12\is\bet{\na\!A}^2\g^6\mm.
\end{eqnarray*}
From this we get
\begin{eqnarray*}
  \int\limits^{t}_{T/2}\is\bet{\na\!A}^2\g^6\mm\,ds
    &\le& c\sqrt{\e_0T}\Bigg(
    \int\limits^{T}_{0}\int\limits_{\Sigma_{63/64}}\bet{\na^2\!A}^2\g^6\mm\,ds
    \Bigg)^{1/2}+c\e_0T
    \stackrel{(\ref{53ia})}{\le}c_1\sqrt T
\end{eqnarray*}
and therefore
\begin{eqnarray*}
  k_1(t)+\frac{1}{2}\int\limits_{T/2}^t\chi_1(s)K_{2}(s)ds
    \le c_2\c\Big(c_2+\frac{c_2(m)}{T}
    \int\limits^{t}_{T/2}\is\bet{\na\!A}^2\g^6\mm\,ds\Big)\le\frac{c_2}{T^{1/2}}.
\end{eqnarray*}
As above we use Gronwall's Lemma obtaining for all $j\geq2$
\begin{eqnarray*}
  k_j(t)+\frac{1}{2}\int\limits_{T/2}^t\chi_j(s)K_{j+1}(s)ds
    \le\frac{c_{2j}}{T^{j-3/2}}+\frac{c_{2j}(m)}{T}\frac{c_{2j-2}}{T^{j-3/2}}
    \le\frac{c_{2j}}{T^{j-1/2}}
\end{eqnarray*}
and therefore it follows from (\ref{indhypth})
\begin{eqnarray}
  \is\bet{\na^{2m+1}\!A}^2\g^{4m+6}\mm\le\frac{c_{2m+2}}{T^{m+1/2}}.
\end{eqnarray}
Therefore, we have now shown that for all $k\in\N_0$
\begin{eqnarray*}
  \norm{\na^k\!A}{2}{2,\Sigma_{3/4}}{}\strich{t=T}~~~
    \le\dfrac{c\big(n,k,\Lambda_{k+1}\big)}{T^{k/2}}.
\end{eqnarray*}
From (\ref{41ia}) and (\ref{421}) we finally obtain
\begin{equation}\label{test}
  \norm{\na^k\!A}{}{\infty,\Sigma_{1/2}}{}\strich{t=T}~~~
  \le\dfrac{c\big(n,k,\Lambda_{k+3}\big)}
    {T^{\frac{k+1}{4}}}
\end{equation}
and thus the claim follows after rescaling and renaming T into t.
\eb
\notiz{Ein paar Bemerkungen machen, in Bezug auf den Beweis von Kuwert-Sch\"atzle}
\section{Blow-up of singularities}
The goal of this section is to prove \cref{blowupthm}. For this
let $f:\Sigma\times[0,T)\to (M^n,g)$ be a smooth maximal Willmore flow defined on a
closed surface $\Sigma,$ $f_t:=f(\:\c\:,t)$ and $\mue {f}(t)=\mue {f_t}$ be the area
measure on $\Sigma$ induced by $\tildeg(t):=(f_t)^*g.$
For $t\in[0,T)$ we let 
$$\nu_{f_t}:=f_t\big(\mue{f_t}\:\llcorner\bet{A_{f_t}}^2 \big)$$
be the (finite) Radon measures
defined on the Borel $\sigma-$Algebra $\mathcal B(M)$ induced by the system of open sets
in $M$.
Let further
$$\chinullf{f}(r,t):=\sup_{p\in M}\nu_{f_t}\big(B_r^g(p) \big)\quad\mbox{and}\quad
\chif{f}(r,t):=\sup_{p\in M}\nu_{f_t}\big(\overline{B_r^g(p)} \big).$$
We want to make the arrangement that all geometric quantities without further
specification refer to $f_t$ (as a map between Riemannian manifolds), e.g.
$\mue{f_t}=\mu,~A_{f_t}=A,~\bet{\;\c\;}_{f_t}=\bet{\;\c\;},\ldots$
Furthermore, when considering a sequence of functions, we want to allow choosing
subsequences without further mention. At first, as remarked in \cite{lecKuw} in case of
$(\R^n,\d_{eucl})$ as the ambient
manifold, one can analogously show that
\begin{equation}\label{uppersemi}
  E(p,r,t):=\nu_{f_t}\big(\overline{B_r^g(p)} \big)
\end{equation}
is upper semi-continuous in all variables and hence the same also holds for $\chi.$ 
The author showed in his diploma thesis (\cite{dipl}, Hilfslemma 5.2), that
$\chinullf{}\!\!$ is lower semi-continuous in the time variable. Using a simple
covering argument, we get for $\Gamma=\Gamma(n)$ as in \cref{hilfslemma51}
$\chif{}\!\!(r,t)\le\Gamma\chi(r/2,t)\le\Gamma\chinullf{}\!\!(r,t)$
for all $r>0$ with $r^2\Vert\mbox{ricci}_{(M,g)}\Vert_{L^\infty(M,g)}\le1.$ Thus we
have in summary
\begin{eqnarray}
  \chinullf{}\!\!(r,t_0)\le\liminf_{t\to t_0}\chinullf{}\!\!(r,t)\le
  \limsup_{t\to t_0}\chinullf{}\!\!(r,t)
  \le\Gamma\chinullf{}\!\!(r,t_0)\label{halbstetcirc}
\end{eqnarray}
and
\begin{eqnarray}
  \Gamma^{-1}\chif{}\!\!(r,t_0)\le\liminf_{t\to t_0}\chif{}\!\!(r,t)\le
  \limsup_{t\to t_0}\chif{}\!\!(r,t)\le\chif{}\!\!(r,t_0).\label{halbstet}
\end{eqnarray}\notiz{bem, warum die bedingung an $\tau^-$ existiert}
The next Theorem is a precursor of \cref{blowupthm}.
\begin{theorem}{~\rm(Existence of a blow-up I).}\label{blowupthmvorl}~Let 
  $f:\Sigma\times[0,T)\to (M,g)$ be a maximal Willmore flow on a closed surface $\Sigma$
  into a Riemannian manifold of bounded geometry (of order \bg) with the property that
  the total area of $(\Sigma,\tildeg(t))$ is uniformly bounded on $(0,T),$ i.e.
  $$\mathfrak M_f:=\sup_{t\in(0,T)}\mue f(\Sigma)<\infty.$$
  Let further $\{t_j\},\{r_j \}$ and $\{p_j \}$ $(j\in\N)$ be given sequences satisfying
  $t_j\nearrow T,~r_j\searrow 0$ and $p_j\in M.$
  There exists constants $\e_{0/1}=\e_{0/1}(n)>0$ such that if 
  $$\liminf_{j\to\infty}\chi(r_j,t_j)\le\e_1^2$$
  and
  $$\tau^-:=\limsup_{j\to\infty}\frac{t_j^--t_j}{r_j^4}<0,\quad\mbox{\rm where}\quad
  t_j^-=\inf\big\{t\in[0,t_j]:\chi(r_j,\:\cdot\:)<\e_0^2 \mbox{\rm~on~}(t,t_j]\big\}$$
  then, after selection of a subsequence and reparametrization, the rescaled flows
  \begin{equation}
    f_j:(\Sigma,\tildeg_{\!j})\times\big[-r_j^{-4}t_j,r_j^{-4}(T-t_j)\big)\to(M,g_j),
    \qquad f_j(p,t):=f_{j,t}(p):=f(p,t_j+r_j^{4}t),\label{blowup}
  \end{equation}
  where $g_j:=r_j^{-2}g$ and $\tildeg_{\!\!j}(t)=(f_{j,t})^*(g_j),$
  converge in $C^{4,1}$ locally on $\widehat\Sigma\times(\tau^-,\tau^+)$  
  to a static  solution, given by a properly immersed Willmore surface 
  $\hat f_0:\widehat\Sigma\to\R^n.$\\\\
  More precisely there exists a sequence of radii $h_j\nearrow\infty,$ local coordinate
  charts  $\varphi_j:B^{g_j}_{h_j}(p_j)\to V_j\supset B_{2j^2}^n,$ a $2$-manifold
  $\widehat\Sigma$ without boundary (possibly not connected or even empty),
  open sets $U_j\subset\!\subset U_{j+1}$ with 
  $\bigcup_{j=1}^{\infty}U_j=\widehat\Sigma,$
  diffeomorphisms $\phi_j:U_j\to f_{j,0}^{-1}\big(\varphi_j^{-1}(B_j^n)\big),$
  a time interval $(\tau^-,\tau^+)\ni0,$
  open sets $\Sigma_j\subset\Sigma$ satisfying
  $\phi_j(U_j)\subset\Sigma_j$ and 
  $f_j(\Sigma_j,J)\subset B^{g_j}_{h_j}(p_j)$ for each
  $J\subset\!\subset(\tau^-,\tau^+)$ provided $j\geq j_0(J),$
  such that for $\bar f_j\circ\phi_j:=\bar f_j\big(\phi_j(\:\cdot\:),\:\cdot\:\big)$
  being well defined locally on $\widehat\Sigma\times(\tau^-,\tau^+),$
  where $\bar f_j:=\varphi_j\circ f_j$
  \begin{eqnarray}
    \bar f_j\circ\phi_j\to\hat f_0\quad\mbox{\rm~locally in $C^{4,1}$ on }
    \widehat\Sigma\times(\tau^-,\tau^+).\label{repconv}
  \end{eqnarray}
  Moreover, if $\theta_j\circ\phi_j\to\theta$ locally uniformly on $\widehat\Sigma,$
  then
  \begin{eqnarray}
    \bar f_j(\mue{f_{j}}\;\llcorner\theta_j)
    \to\hat f_0(\mue{\hat f_0}\;\llcorner\theta)
    \quad\mbox{\rm~for all }\tau\in(\tau^-,\tau^+)\label{weakconv10}
  \end{eqnarray}
  weakly as Radon measures.
\end{theorem}
\bf{Remark:}\rm~Again, we tried to stick to the analysis and the way of proceeding
as in \cite{lecKuw}. Deviating from this, we assume the given arbitrary sequence $r_j>0$
to be decreasing owing to the restriction to blow-ups. Related to this, we do not prove
an upper bound $\tau^-<0$ for the backward lifespan for the rescaled flows making an
ad-hoc assumption instead. \er\\
\bf{Proof of \cref{blowupthmvorl}}:\rm~Let
$\e_0>0$ be as in \cref{lifespan} and $\e_1:=C\e_0/6$ for $C$ as in
(\ref{lifespan2}).
It is easy to see that $f$ and $f_j$ are equivalent maximal Willmore flows.
We prove lower bounds for the forward lifespan of the rescaled flows. Choosing a
subsequence, we may assume that
\begin{eqnarray}
  r_j\big(r_j^2(1+\W(f_0))+\mathfrak M_f\big)\sum_{i=0}^{\bg}
  \norm{D^{i\!}R}{\frac3{i+2}}{L^\infty(M,g)}{}+r_j\inj_{(M,g)}^{-1}\le\e_1^2
  \label{kleinheitsvssN}
\end{eqnarray}
and $\chi(r_j,t_j)\le2\e_1^2.$
Note that with respect to the scaled metric, (\ref{kleinheitsvssN}) holds for $r_j=1.$ 
Letting $$\tau^+:=\liminf_{j\to\infty}\frac{t_j^+-t_j}{r_j^4},\quad\mbox{where}\quad
t_j^+:=\sup\big\{t\in[t_j,T):\chi(r_j,\,\c\,)<\e_0^2\mbox{ on }[t_j,t)\big\},$$
\cref{lifespan} implies for $0<\e_1\le C\e_0/6<\e_0$
\begin{eqnarray*}
  \liminf_{j\to\infty}\frac{T-t_j}{r_j^4}\geq\tau^+\geq C\log\frac{C\e_0^2}{3\e_1^2}>0.
\end{eqnarray*}
Also, we obviously have for the backward lifespan
$$\lim_{j\to\infty}\Big(-\frac{t_j}{r_j^4}\Big)=-\infty.$$
Employing that $\chif{f_j}\:(1,\,\c\,)\le\e_0^2$ locally on $(\tau^-,\tau^+),$
we obtain from \cref{interiorestimates} that for $k\le\bgminusdrei$
\begin{eqnarray}
  \norm{\na^k\!A_j}{}{L^\infty(\Sigma,\tildeg_{\!j})}{}\le c(n)
  \quad\mbox{locally on}\quad(\tau^-,\tau^+).\label{curvboundscaled}
\end{eqnarray}
In what follows, we want to employ a compactness result stated in \cref{PBkorollar},
which is a slight variant of a
compactness theorem, where the latter has originally been proven by J. Langer in
\cite{JL} and then has been generalized
by P. Breuning (cf. \cite{PB}, Corollary 2.28).\notiz{check ob man unten $\Sigma$
braucht??}
\begin{theorem}\label{PBkorollar}\rm(cf. \cite{PB}).\it~Let
  $\bar f_j:\Sigma_j\to\R^n$ be a sequence of isometric $C^k-$immersions, where
  $\Sigma_j$ is  an $m$-manifold without boundary, such that
  $\bar f_j^{-1}(B_{2j}^n)\subset\Sigma_j\subset(\Sigma^m,\bar g_j)$
  and
  $\Sigma$ is a closed surface. 
  With the image measure $\mu^j=\bar f_j(\mue{\bar f_j}\:)$ assume that for all $j\in\N$
  \begin{eqnarray}
    \mu^j(B_R^n)&\le& C(R) \mbox{~~~~for all } 0<R<2j,\\
    \norm{\barnabla^{\,l}\!\barA_j}{}{L^\infty(B_R^n,\bar g_j)}{}&\le& C_l(R)
    \mbox{~~~for all~}0<R<2j\mbox{~and~} 0\le l\le k-2,
  \end{eqnarray}
  where $\barA_j$ denotes the second fundamental form of $\bar f_j.$
  Then there exists a proper immersion $\hat f:\widehat\Sigma\to\R^n,$ where 
  $\widehat\Sigma$ is again an $m$-manifold without boundary, such that after passing
  to a subsequence there is a sequence
  of diffeomorphisms $$\phi_j:U_j\to(\bar f_j)^{-1}(B_j)\subset\Sigma_j,$$
  where $U_j\subset\widehat\Sigma$ are open sets with $U_j\subset\!\subset U_{j+1}$ and 
  $\widehat\Sigma=\cup_{j=1}^{\infty}U_j,$ such that 
  \begin{eqnarray}
    \norm{\bar f_j\circ\phi_j-\hat f}{}{C^0(U_j)}{}\to 0,\label{patconvergence}
  \end{eqnarray}
  and moreover 
  $\bar f_j\circ\phi_j\to \hat f$ locally in $C^{k-1}$ on $\widehat\Sigma$.\eb
\end{theorem}
\bf{Remark:}\rm~The surface $\widehat\Sigma$ is empty, if the $\bar f_j$ diverge
uniformly, and is possibly not con\-nected.\er\\
To be able to apply this corollary locally to the blow-up sequence $\{f_{j,0} \},$ we
want use for any $j\in\N$ $g_j$-harmonic coordinates 
$\varphi_j:B_{h_j}^{g_j}(p_j)\to V_j\subset\R^n$ of radius $h_j:=4j^2$
centred at $p_j$ with $\varphi_j(p_j)=0.$ To
check the assumptions in
\cref{harmcoord1}, we may further assume that $r_j\searrow0$ is such that
$r_j\Lambda_{(M,g)}(\bg)\le(4j^3)^{-1}c(n)$ for some given small $c(n)>0,$ so that 
$$h_j\Lambda_{(M,g_j)}=4j^2r_j\Lambda_{(M,g)}\le c(n)j^{-1}$$
and thus we obtain the bounds
\begin{equation}
  \begin{array}{lll}
    i)~~~~(1-1/j)\delta\le G_j\le(1+1/j)\delta
    &&\mbox{on} V_j\\\label{harmabsch}
    ii)~~\sup_{V_j}\bet{\p^\g G_j}\le j^{-2\g}&&\mbox{for
    all } 1\le\g\le \bgpluseins\\
    iii)~\sup_{V_j}\bet{\p^\g\Gamma^{(j)} }\le j^{-2\g-2}
    &&\mbox{for all } 0 \le\g\le \bg\\
  \end{array}
\end{equation}
and $G_j\strichklein{x=0}\quad=(\d_{\a\b}).$
Here, $G_j$ denotes the matrix of $(\varphi_j)_*(g_j)$ with respect to the standard
coordinate frame $\{e_\a \}_{1\le\a\le n}$ and $\{\Gamma^{(j)}\}_{j\in\N}$ are the
associated Christoffel's symbols.
Now since $\Norm{\p_tf_j}_{L^\infty(\Sigma,\w g_j)}\le c(n)$
we get that for any $x\in\Sigma$
\begin{eqnarray}
  d^{g_j}(f_j(x,t_1),f_j(x,t_2))\le\Big\lvert \int_{t_1}^{t_2}
  \bet{\p_tf_j}_{g_j}(x,\tau) d\tau\Big\rvert
  \le c(n)\bet{t_2-t_1}.\label{timeintegral}
\end{eqnarray} 
Letting $0\in I:=[t_1,t_2]\subset(\tau^-,\tau^+)$ be a compact interval, we have a well
defined local coordinate representations for all $f_j$ (for $j$ large enough).
Namely, with $\Sigma_j:=f_{j,0}^{-1}(B_{4j}^{g_j}(p_j))\subset\Sigma,$ we
define one-parameter families of isometric immersions
\begin{equation}
  \tilde f_{j}:(\Sigma_j,\w g_j)\times I\to (V_j\subset \R^n,G_j)\label{tildefj}
\end{equation}
and
\begin{equation}
  \bar f_{j}:(\Sigma_j,\bar g_j)\times I
       \to (V_j\subset\R^n,\d_{eucl})\label{barfj}
\end{equation}
using the harmonic coordinates from above.
Here, $\tilde f_j$ is the coordinate representation of $f_j,$ and 
$\bar f_j:=\tilde f_j$ as maps between manifolds. For the rest of this section, a tilde
and index $j$ indicates that a geometric quantity is induced by $\tilde f_j,$ as in
$\w\na$ or $\w A_j$ etc. The same arrangement holds for barred quantities such as
$\barnabla$ or $\barA_j$ etc. 
By naturality of the geometric quantities and using (\ref{curvboundscaled}), we have that
$\tilde f_j$ is also a Willmore flow with
\begin{equation}
  \norm{\w \na^k\!\w A_j}{}{L^\infty(\Sigma_j,\w g_j)}{}
    \le c(n)\label{secondffbound}
\end{equation}
for all $t\in I$ and $k\le \bgminusdrei.$
From \cref{relativecurv} it follows that also
$$\norm{\barnabla^{k}\!\barA_j}{}{L^\infty(\Sigma_j,\bar g_j)}{}
  \le c(n,\Gamma_{\bgminusdrei},A_{\bgminusdrei})\le c(n)$$
for all $k\le \bgminusdrei,~t\in I,$ and in particular for $t=0.$ 
As in (\ref{detbil}) we see that $\frac{1}{2}\bar g_j\le\tildeg\!_j\le 2\bar g_j$ 
implies the coordinate invariant relation $\frac{1}{4}\det(\bar g_j)
\le\det(\tildeg\!_j)\le 4\det(\bar g_j)$ and thus from \cref{massbound3}
\begin{equation*}
  \bar f_j(\mue{\bar f_j}\:)(B_R^n)\le2\tilde f_j(\mue{\tilde f_j}\:)(B^n_R)\le C R^2
\end{equation*}
for all $R<2j.$ From \cref{PBkorollar} 
we get after passing to a subsequence that
\begin{equation*}
  \bar f_j(\phi_j,0)\to\hat f_0
\end{equation*}
locally in $C^{13},$ where $\hat
f_0:=(\widehat\Sigma,\hat g)\to(\R^n,\d_{eucl})$ is a 
properly and isometrically immersed 2-manifold without boundary.
Fix a countable atlas of $\widehat\Sigma$ comprising of Riemannian normal coordinate
charts
$\{\psi_k:B_{2r_k}^{\hat g}(x_k)\to B_{2r_k}^2\subset\R^2 \}_{k\in\N}$ with 
\begin{eqnarray}
  2^{-1}\d\le(\hat g_{im})=:\widehat G\le 2\d.\label{metequiv3}
\end{eqnarray}
For $k\in\N$ fixed, $r:=r_k,$ the harmonic
coordinates from above induce a sequence of coordinate representations of
$f_j\circ \phi_j=f_j(\phi_j(\;\c\;),\;\c\;)$
\begin{equation*}
  \w F_{j}:(B_{r}^2,\w{G}_j(t))\times I\to (V_j\subset\R^n,G_j)
\end{equation*}
and
\begin{equation*}
  \bar F_{j}:(B_r^2,\barG_j(t))\times I\to (V_j\subset\R^n,\d_{eucl}) 
\end{equation*}
for $j$ large enough, such that $\w F_{j,t}$ and $\bar F_{j,t}$ are isometric immersions
and $\w F_j\equiv\bar F_j$ as maps between manifolds.
We now want to show that $\bar f_j\circ\phi_j\to\hat f$ in $C^{5,1}$ locally on
$\widehat\Sigma
\times (\tau^-,\tau^+).$ For this, we show bounds for lower order coordinate
derivatives of $\bar f_j\circ\phi_j$ similar to (\ref{coordbound2}). We check the
assumptions in
\cref{zeitableitungen}. At first, we note that $\w{G}_j(0)
\rightrightarrows\widehat G$ (i.e. uniformly) on $B_r^2\subset\!\subset B_{2r}^2$
and thus, using (\ref{metequiv3}), we have
\begin{itemize}
  \item [i)] $\frac12\d\le G_j\le 2\d$\quad and \quad$\frac13\d\le\w G_j(0)\le 3\d$ 
\end{itemize}
and from (\ref{harmabsch}),(\ref{secondffbound}) we have
\begin{itemize}
  \item [ii)] $\bet{\p^\g\Gamma^{(j)}}\le\Gamma_{\bg}:=1$\quad for all\quad $\g\le\bg$
  \item [iii)] $\sum_{l=0}^{12}
               \norm{\w\na^l\!\w A_j}{}{L^\infty(\w G_j(t),B_r^2)}{}\le A_{12}:=c(n)$
               \quad for all $t\in I.$
\end{itemize}
As in (\ref{metricevol4}) we get for the coordinate functions  
$\p_t(\w G_j)_{rs}=(\w\na_{\!\p_t}\w G_j)_{rs}$
\begin{itemize}
  \item [iv)] $\bet{\p_t\w G_j}\le c(n).$
\end{itemize}
If $A\in\R^{n\times n}$ is nonsingular, the inverse matrix may be represented by the
formula $A^{-1}=\frac{S(A)}{\det A}$ for some $S\in C^\infty(\R^{n^2},\R^{n^2}),$
Since $\w G_j(0)\rightrightarrows\widehat G$ (uniformly) on $B_r^2,$
and because $\frac13\d\le\w G_j(0)$ we infer that $\w\Gamma^{(j)}(0)\rightrightarrows
\widehat\Gamma$ and therefore $\bet{\w\Gamma^{(j)}(0)}\le\w\Gamma_0,$ where
$\wgm0\in\R^{>0}$ is a constant and $j$ is assumed to be chosen large enough. 
Differentiating Kronecker's symbol $\d^k_i=(\w G_j)_{il}(\w G_j)^{lk}$
with respect to the space variables we get using induction that
$\p^\g\w\Gamma^{(j)}(0)\rightrightarrows\p^\g\widehat\Gamma\mbox{ for }\g\le K$ and
hence 
\begin{eqnarray}
  \bet{\p^\g\w\Gamma^{(j)}(0)}\le\wgm{K}\strichklein{t=0}\quad<\infty\label{GamNullBound}
\end{eqnarray}
uniformly for $j$ sufficiently large. Now arguing
analogously to \cref{gammadotlemma}, we have the bounds
$$\bet{\p^\g(\p_t\w\Gamma^{(j)})},~
\bet{\p^\g\w\Gamma^{(j)}}\le c(n,\g,A_{\g+3},R_{\g+1},c_{\w G_j}
  ,\wgm{\g}\strich{t=0}\quad,\bet{I})$$
and therefore
\begin{itemize}
  \item [v)] $\bet{\p^\g\w\Gamma^{(j)}(t)}\le\w\Gamma_{9}
    := c(n,\wgm{9}\strichklein{t=0}\quad,\bet{I})$ for $\g\le 9$ and all $t\in
  I.$ 
\end{itemize}
From \cref{zeitableitungen} (with $k=0,~r=2,~s=0,~m=4$) it follows that on the
time-interval $I$ we have
$$\bet{\p^{p+1}(\bar f_j\circ\phi_j)}
\le c(n,\wgm{9},A_{9},\Gamma_{\!9},\bet{I})\le
c(n,\wgm{9}\strichklein{t=0}\quad,\bet{I})\quad(p\le 10)$$
and $$\bet{\p^p\p_t^l(\bar f_j\circ\phi_j)}\le c(n,\wgm{9},A_{12},     
\Gamma_{9},R_{10},\bet{I})
\le c(n,\wgm{9}\strichklein{t=0}\quad,\bet{I})\quad(p\le6\mbox{ and }l\le2).$$
By Arzel\`a-Ascoli a subsequence converges. Choosing a subsequence once more, we
have thus shown that
\begin{eqnarray*}
  \bar f_j\circ\phi\to\hat f\mbox{ in $C^{5,1}$ on  }U_{m}\times I\quad
    \mbox{for}\quad j\geq j_0(m,I)
\end{eqnarray*}
and choosing a diagonal sequence we have finally proven that
\begin{eqnarray}
  \bar f_j\circ\phi\to\hat f\mbox{ in $C^{5,1}$ locally on
    }\widehat\Sigma\times(\tau^-,\tau^+).\label{ck1conv}
\end{eqnarray}
Since $\w G_j(t)\rightrightarrows\widehat G(t)$ and 
$$\w G_j(t)\stackrel{(\ref{tildeab})}{\geq}c(n,\bet{t})\w G_j(0)\stackrel{i)}{\geq}
\frac13c(n,\bet{t})\d\qquad\mbox{for some $c>0$}$$
we see that $\hat f_t:(\widehat\Sigma,\hat g(t))\to\R^n$ is an
isometric immersion for all $t\in(\tau^-,\tau^+)$. Furthermore, \mbox{$\hat
f:(\widehat\Sigma,\hat
g)\times(\tau^-,\tau^+)\to\R^n$} is a Willmore flow.\\\\
To see that $\hat f$ is static, we let $\tau^-<\tau_1<\tau_2<\tau^+,$
$U\subset\!\subset\widehat\Sigma$ and compute for $j$ sufficiently large
\begin{eqnarray*}
  \int^{\tau_2}_{\tau_1}\int_U\bet{\bf{W}\it(\tilde f_j\circ\phi_j)\rm}^2
     \it d\mue{\tilde f_j\circ\phi_j}\qquad d\tau
  &=&\int^{\tau_2}_{\tau_1}\int_{\phi_j(U)}\bet{\bf{W}\it(\tilde f_j)\rm}^2
     \it d\mue{\tilde f_j}\;d\tau\\
  &\le&\int^{\tau_2}_{\tau_1}\int_{\Sigma}\bet{\bf{W}\it(f_j)\rm}^2
     \it d\mue{f_j}\:d\tau\\
  &\stackrel{(\ref{444})}{=}&\W(f_j)\strichklein{$\tau=\tau_1$}\;\:\quad
    -\W(f_j)\strichklein{$\tau=\tau_2$}\\
  &=&\W(f)\strichklein{$t=t_j+r_j^4\tau_1$}\qquad\qquad
    -\W(f)\strichklein{$t=t_j+r_j^4\tau_2$}\,\qquad\qquad.\\
\end{eqnarray*}
Since $\W(f)$ is monotonically decreasing with respect to time, and
$t_j+r_j^4\tau_1\rightarrow T,$ it follows that
$\bf{W(\rm\it\hat f)}\rm\equiv 0$ and hence $\hat f$ is static for $\tau\geq\tau_1.$\\\\
To prove (\ref{weakconv10}) we fix $\tau\in(\tau^-,\tau^+).$ Integrating the Willmore
flow equation with respect to
time as in (\ref{timeintegral}), using (\ref{kleinheitsvssN}) and
(\ref{curvboundscaled}), we get
\begin{eqnarray*}
  \norm{\bar f_j(\phi_j(\:\c\:),\tau)-\bar f_j(\phi_j(\:\c\:),0)}{}{C^0(U_j)}{}
    \le 2\sup_{x\in U_j}d^{\,G_j}\big(
    \tilde f_j(\phi_j(\:\c\:),\tau),\tilde f_j(\phi_j(\:\c\:),0)\big)\le C,
\end{eqnarray*}
and by passing to the limit 
$$\norm{\hat f(\:\c\:,\tau)-\hat f(\:\c\:,0)}{}{C^0(\widehat\Sigma)}{}\le C.$$
In particular $\hat f(\:\c\:,\tau):\widehat\Sigma\to\R^n$ is proper, and combining
with (\ref{patconvergence}) yields
$$\norm{\bar f_j(\phi_j(\:\c\:),\tau)-\hat f(\:\c\:,\tau)}{}{C^0(U_j)}{}\le C.$$
This shows that $\bar f_j\circ\phi_j$ is uniformly proper, i.e.
$\{x\in\widehat\Sigma:\bet{\bar f_{j,\tau}\circ\phi_j(x)}\le R \}\subset
K_{R,\tau}\subset\!\subset \widehat\Sigma.$
Using this fact together with (\ref{ck1conv}) we get, for compact
$K\subset U\subset\R^n$, $U$ open and bounded, at time
$\tau\in(\tau^-,\tau^+)$
$$(\bar f_j\circ\phi_j)^{-1}(K)\subset\hat f^{-1}(U)\quad
  \mbox{for $j$ sufficiently large.}$$
Using the equivalence (\ref{harmabsch}) we can conclude
\begin{eqnarray*}
  \bar f_j(\mue{f_j}\:\llcorner\theta_j)(K)
  &=&\int_{\bar f_j^{-1}(K)}\theta_j\,d\mue{f_j}\\
  &=&\int_{(\bar f_j\circ\phi_j)^{-1}(K)}\theta_j\circ\phi_j\,
    d\mue{f_j\circ\phi_j}\\
  &\le&\int_{\hat f^{-1}(U)}\theta_j\circ\phi_j\,
    d \mue{\tilde f_j\circ\phi_j}\qquad\quad\qquad\mbox{(U is bounded)}\\
  &\le&(1+1/j)\int_{\hat f^{-1}(U)}\theta_j\circ\phi_j\,
    d \mue{\bar f_j\circ\phi_j}\qquad. 
\end{eqnarray*}
Letting $j\to\infty$ and taking the infimum over all bounded $U\supset K,$ we see
$$\limsup_{j\to\infty}\bar f_j(\mue{f_j}\:\llcorner\theta_j)(K)
\le \hat f(\mue{\hat f}\,\llcorner\theta)(K).$$
Here, we used that $K$ can be approximated from above using only \it{bounded}\rm~open
sets. Furthermore, we also have $\hat f^{-1}(K)\subset(\bar
f_j\circ\phi_j)^{-1}(U)$ for $j$ large enough, which implies
\begin{eqnarray*}
  \frac j{j+1}\int_{\hat f^{-1}(K)}\theta_j\circ\phi_j\,d
    \mue{\bar f_j\circ\phi_j}\quad\;
  &\stackrel{(\ref{harmabsch})}{\le}&\int_{\hat f^{-1}(K)}
    \theta_j\circ\phi_j\,d\mue{\tilde f_j\circ\phi_j}\\
  &\le&\int_{(\bar f_j\circ\phi_j)^{-1}(U)}
    \theta_j\circ\phi_j\,d\mue{\tilde f_j\circ\phi_j}\\  
  &=&\int_{\bar f_j^{-1}(U)}\theta_j\,d\mue{f_j}\\
  &=&\bar f_j(\mue{f_j}\:\llcorner\theta_j)(U).
\end{eqnarray*}
For $j\to\infty$ we get, when taking the supremum with respect to K,
$$\hat f(\mue{\hat f}\,\llcorner\theta)(U)\le\liminf_{j\to\infty}
\bar f_j(\mue{f_j}\:\llcorner\theta_j)(U).$$
Note that it is sufficient that the last inequality is valid only for all
bounded $U\subset\R^n$ (cf. the proof of Theorem 1 in section 1.9 of \cite{EvGa}).
Alternatively, one may argue that the latter inequality also holds for all $U\subset\R^n$
open, not necessarily bounded by continuity from below of the measure $\hat f(\mue{\hat
f}\,\llcorner\theta)$. In any case, (\ref{weakconv10}) is now settled.\eb

Before we prove \cref{blowupthm} we want to formulate the proposed convergence more
precisely.\\\\ 
\bf{\cref{blowupthm}~}\rm(Part II): Under the assumptions and definitions of
\cref{blowupthm} the following convergence result holds:\;\label{partII}
  There exist sequences $R_j\nearrow\infty$ and $p_j\in M,$ local
  coordinate
  charts  $\varphi_j:B^{g_j}_{R_j}(p_j)\to V_j\supset B_{2j^2}^n,$ a $2$-manifold
  $\widehat\Sigma$ without boundary (possibly not connected),
  open sets $U_j\subset\!\subset U_{j+1}$ with 
  $\bigcup_{j=1}^{\infty}U_j=\widehat\Sigma,$
  diffeomorphisms $\phi_j:U_j\to f_{j,0}^{-1}\big(\varphi_j^{-1}(B_j^n)\big),$
  open sets $\Sigma_j\subset\Sigma$ satisfying
  $\phi_j(U_j)\subset\Sigma_j$ and 
  $f_j(\Sigma_j,J)\subset B^{g_j}_{R_j}(p_j)$ for each
  $J\subset\!\subset\R$ provided $j\geq j_0(J),$
  such that for $\bar f_j\circ\phi_j:=\bar f_j\big(\phi_j(\:\cdot\:),\:\cdot\:\big)$
  being well defined locally on $\widehat\Sigma\times\R,$
  where $\bar f_j:=\varphi_j\circ f_j$
  \begin{eqnarray}
    \bar f_j\circ\phi_j\to\hat f_0\quad\mbox{locally in $C^{4,1}$ on }
    \widehat\Sigma\times\R.\label{repconvIntro}
  \end{eqnarray}
\eb

\notiz{check ob $\W(f_0)$ \"uberhaupt gebraucht wird?}
\bf{Proof of \cref{blowupthm}:~}\rm 
Let $\e_0,\e_1$ be as in \cref{blowupthmvorl}, $0<\e<\min\{\e_T,\e_1^2/2\},$ and
$\{r_j\}_{j\in\N}$ be an arbitrary given sequence with
$r_j\searrow 0.$
Consider the times
$$t_j:=\sup\{t\geq 0:~\chi(r_j,\,\c\,)<\e^2\mbox{ on }[0,t) \}.$$
Clearly $t_j\le t_{j+1}.$ By definition of $\e_T$ in (\ref{epst}) and since
$0<\e<\e_T$ we also
get that $t_j<T.$ By smoothness of the flow we furthermore have $\lim_{j\to\infty}t_j=T.$
This could also be seen by using the upper semicontinuity of $\chi$ in the time
variable.
By maximality of $t_j$ and again by upper semicontinuity of $\chi(r_j,\:\c\:)$, we get 
$\chi(r_j,t_j)\geq \e^2.$
Next, we show the existence of so-called concentration points. For this, let $j\in\N$ be
arbitrary but fixed, and $K:=\{p\in M:\dist_g(p,f(\Sigma,t_j))\le 2\}.$
Since an upper semicontinuous function attains its supremum over compact sets, there
exist $p_j\in M$ such that
\begin{eqnarray}
  \e^2\le\chi(r_j,t_j)=\sup_{p\in M}E(p,r_j,t_j)=\sup_{p\in
K}E(p,r_j,t_j)=E(p_j,r_j,t_j).\label{curvconcent}
\end{eqnarray}
On the contrary, from (\ref{halbstet}) we get that
$$\chi(r_j,t_j)\le\Gamma\liminf_{t\nearrow t_j}\chi(r_j,t)
\le\Gamma\e^2\le\e_1^2$$
when restricting to $\e^2\le\e_1^2\Gamma^{-1},$ where $\e_1^2<\e_0^2$ is as in
\cref{blowupthmvorl}.
Moreover, $t_j^-=0$ and thus
$\tau^-=\limsup_{j\to\infty}\frac{t_j^--t_j}{r_j^4}=-\infty.$
From this we get the desired (sub-)convergence locally on the interval
$(-\infty,\tau^+)$ from \cref{blowupthmvorl}.
To see that $\tau^+=\infty,$ we claim that for any
$\tau\in(-\infty,\tau^+),~\varrho\in(0,1),$
\begin{eqnarray}
  \limsup_{j\to\infty}\chif{f_j}\:(\varrho,\tau)\le\e^2.\label{tauplusunendlich}
\end{eqnarray}
Otherwise, since $\chif{f_j}\:(\varrho,\tau)=\chi(r_j\varrho,t_j+r_j^4\tau),$ we find a
subsequence with
\begin{eqnarray}
  \int\limits_{f^{-1}(\overline{B^g_{r_j\varrho}(P_j)})}
  \bet{A}^2\mm\strich{$t=t_j+r_j^4\tau$}\qquad\qquad\geq\tilde\e^2,
\end{eqnarray}
for some $P_j\in M$ and $\tilde\e>\e.$ 
Rescaling parabolically, i.e. letting $f_j$ be defined as in (\ref{blowup}) with
$g_j:=(\varrho r_j)^{-1}g$, and
using
harmonic coordinates as above, with the exception to be centred at $P_j$ instead of
$p_j,$ we can again define coordinate representations $\w{\mathcal F}_j$
as in (\ref{tildefj}) and Riemannian metrics $\mathcal G_j$ as in (\ref{harmabsch}).
We have, using lower semi-continuity (\ref{halbstetcirc}) of 
$\chinullf{f_j}\:(1,\:\c\:),$
$$\int\limits_{\w{\mathcal F}_{j,0}^{-1}\big(B^n_{\sqrt{1-1/m}}\big)}
\!\!\!\!\!\!\bet{A_{\w{\mathcal F}_j}(0)}^2\mm_{\w{\mathcal F}_j}(0)\le\e^2
\quad\mbox{but}\quad
\int\limits_{\w{\mathcal F}_{j,\tau}^{-1}(\overline{B^n_{\varrho\sqrt{1+1/m}}})}
\!\!\!\!\!\!
\bet{A_{\w{\mathcal F}_j}(\tau)}^2\mm_{\w{\mathcal F}_{j}}(\tau)\geq\tilde\e^2$$
for $j\geq m,$ since by equivalence of the metrics 
$(1-1/j)\d\le\mathcal G_j\le(1+1/j)\d$ it holds 
$$B^n_{\sqrt {1-1/m}}\subset B_1^{\mathcal G_j}(0)\quad\mbox{and}\quad
\overline{B_\varrho^{\mathcal G_j}(0)}\subset\overline{B^n_{\varrho\sqrt{1+1/m}}}.$$
Analogously to what is shown above, the $\w{\mathcal F}_j$ (sub-)converge after
reparametrization locally in $C^{5,1}$ on $\widehat\Sigma^*\times(-\infty,\tau^+).$ In
particular, convergence (weakly as Radon measures) (\ref{weakconv10}) now
implies
$$\int\limits_{\widehat{\mathcal F}_0^{-1}(B^n_{\sqrt{1-1/m}})}
\!\!\!\!\!\!\bet{A_{\widehat{\mathcal F}_0}}^2\mm_{\widehat{\mathcal F}_0}\le\e^2
<\tilde\e^2\le
\int\limits_{\widehat{\mathcal F}_0^{-1}(\overline{B^n_{\varrho\sqrt{1+1/m}}})}
\!\!\!\!\!\!\bet{A_{\widehat{\mathcal F}_0}}^2\mm_{\widehat{\mathcal F}_0}<\infty$$
for all $m\in\N,$ recalling that $\widehat{\mathcal F}_0$ is proper by
\cref{PBkorollar}, and the reparametrizations of $\w{\mathcal F}_{j,\tau}$ are uniformly
proper. Letting $m\to\infty$ we get, by continuity from above and from below
of the Radon measure $\nu:=\widehat{\mathcal F}_0\big(\mu_{\widehat{\mathcal
F}_0}\llcorner
\bet{A_{\widehat{\mathcal F}_0}}^2\big),$ the contradiction
$$\nu\big(B_1^n \big)\le\nu\big(\overline{B_\varrho^n} \big).$$ Therefore we have
now shown (\ref{tauplusunendlich}).
Now (\ref{lifespan2}) implies for any $\tau\in(-\infty,\tau^+)$
$$\tau^+=\liminf_{j\to\infty}\frac{t_j^+-t_j}{r_j^4}\geq\tau+C\log\frac{C\e^2_0}
{2\e_1^2}\geq\tau+c\log 2,$$
by definition of $\e_1$ in 
Letting $\tau\nearrow\tau^+$ we
conclude that $\tau^+=\infty.$\\\\
Using (\ref{curvconcent}) we can further estimate
\begin{eqnarray*}
  \e^2
  &\le&\!\!\!\!\!\int\limits_{f^{-1}
    (B^g_{r_j}(p_j))}\!\!\!\!\!\!\bet{A_f}^2d\mue f\strich{$t=t_j$}\\
  &=&\int\limits_{f_j^{-1}(B^{g_j}_1(p_j))}\!\!\!\!\!\!
    \bet{A_{f_j}}^2d\mue{f_j}\;\strich{$t=0$}\\
  &\le&\!\!\!\int\limits_{
    \bar f_{j}^{-1}(\overline{B^{n}_{\sqrt{1+1/m}}})}\!\!
    \bet{A_{\tilde f_{j}}}^2d \mue{f_{j}}\;\strich{t=0}
\end{eqnarray*}
for $j\geq j_0(m)$ and thus, again using the weak convergence as Radon measures
(\ref{weakconv10}) and continuity from above of the limit measure we get
\begin{eqnarray}
  0<\e^2\le\int\limits_{\hat f_0^{-1}(\overline{B_1^n})}
  \bet{A_{\hat f_0}}^2d\mue{\hat f_0}\;\,.\label{nontrivlimit}
\end{eqnarray}
For any compact set $C\subset\widehat\Sigma$ we have
$\hat f_0(C)\subset B^n_{R/2}$ for some $R<\infty,$ and since 
$\tilde f_{j,0}\circ\phi_j\to\hat f_0$ in $C^0$ locally on $\widehat\Sigma$ we see that
also $$\tilde f_{j,0}\circ\phi_j(C)\subset B^n_{R/2}\subset B_R^{G_j}(0)
$$
for large $j$. Using this, we further get
\begin{eqnarray*}
  \chi(\varrho,t_j)&\geq&\chi(r_jR,t_j)\\
  &\geq&\int\limits_{f^{-1}
    (B^g_{r_jR}(p_j))}\!\!\!\!\!\!\bet{A_f}^2d\mue f\strich{$t=t_j$}\\
  &=&\int\limits_{f_j^{-1}(B^{g_j}_R(p_j))}\!\!\!\!\!\!
    \bet{A_{f_j}}^2d\mue{f_j}\;\strich{$t=0$}\\
  &=&\!\!\!\!\!\!\int\limits_{(\tilde f_{j,0}\circ\phi_j)^{-1}(B^{G_j}_R(0))}\!\!\!\!\!\!
    \!\!\!\!\!\!\bet{A_{\tilde f_{j,0}\circ\phi_j}}^2d\mue{
    \tilde f_{j,0}\circ\phi_j}\\    
  &\geq&\int\limits_{C}
    \bet{A_{\tilde f_{j,0}\circ\phi_j}}^2d\mue{\tilde f_{j,0}\circ\phi_j}\;\qquad. 
\end{eqnarray*}
Again, since $\tilde f_{j,0}\circ\phi_j\to\hat f_0$ in $C^2$ locally on $\widehat\Sigma$
and using (\ref{harmabsch}) we get, letting $j\to\infty,$
$$\liminf_{j\to\infty}\chi(\varrho,t_j)\geq\int\limits_C
\bet{A_{\hat f_0}}^2d\mue{\hat f_0}$$ 
and further, letting $C\nearrow\widehat\Sigma,$ we
obtain $$\e_T^2\geq\int\limits_{\widehat\Sigma}\bet{A_{\hat f_0}}^2d\mue{\hat f_0}\;\,.$$
The second estimate in (\ref{epstIntro}) follows from the nontriviality
(\ref{nontrivlimit}) of the limit surface $\widehat\Sigma.$ The last estimate in
(\ref{epstIntro}) follows from the fact that if $f:\Sigma\to\R^n$ is a properly
immersed Willmore surface with $\lVert A\rVert_{L^2}<\eps_0,$ then f is a standard plane
(\cite{lecKuw}, Corollary 2.1).
Finally, (\ref{nontrivialIntro}) follows when choosing 
$\e:=\min\{\mathcal E_n/2,\e_1/4 \}=c(n)>0.$
Using the mass-density estimate (\ref{massdensity1}) we obtain
\begin{eqnarray*}
  0<c(n)&\le&\int\limits_{f^{-1}(\overline{B_{r_j}^g(p_j)})}\!\!\!\!
    \bet{A_f}^2d\mue f\strich{$t=t_j$}\\
  &\le&\norm{A_f(t_j)}{2}{L^\infty(\Sigma,\tildeg(t_j))}{}
  \mue{f_{t_j}}\;\,\big(f_{t_j}^{-1}(\overline{B_{r_j}^g(p_j)}) \big)\\
  &\le& C^*r_j^2 \norm{A_f(t_j)}{2}{L^\infty(\Sigma,\tildeg(t_j))}{},
\end{eqnarray*}
where $C^*<\infty$ abbreviates the right-hand side of (\ref{massdensity1}). Therefore the
equation $\tau^+=\infty$ implies
$$(T-t_j)^{1/4}\norm{A_f(t_j)}{}{L^\infty(\Sigma,\tildeg(t_j))}{}\geq
c\:\c\Big(\frac{t_j^+-t_j}{r_j^4} \Big)^{1/4}\to\infty.$$
As in \cite{lecKuw} for the euclidean case, we may say, adopting the terminology from the
mean curvature flow, that a finite time singularity of the Willmore flow in a manifold of
bounded geometry is always of type II.\\\\
The statement in the Remark below \cref{blowupthm} is part of \cref{massbound3}.\eb\\
With the same topological argument as in \cite{KS01}, Lemma 4.3, we finally have the
following Lemma.
\begin{lemma}\rm(cf. \cite{KS01}).\it~Let
  $\hat f:\widehat\Sigma\to\R^n$ be the blow-up constructed above. If $\widehat\Sigma$
  contains a compact component C, then in fact $\widehat\Sigma=C$ and $\Sigma$ is
  diffeomorphic to C.
\end{lemma}
\bb For $j$ sufficiently large, $\phi_j(C)$ is both open and closed in $\Sigma.$ Hence
by connectedness of $\Sigma$ we have $\Sigma=\phi(C)$ and thus $\widehat\Sigma=C.$\eb\\
Thus, the blow-up surface is either an embedded round sphere, or contains at least one
component that is a nonumbilic (compact or noncompact) Willmore surface.  
As a consequence, provided that the Willmore energy of the initial surface is strictly
below $8\pi,$ the case of an embedded sphere can be excluded by an inequality of
Li and Yau \cite{LiYau}. Alternatively, embedded spheres can also be excluded if
it is possible to establish a lower area bound (cf. \cite{KS01}, Theorem 4.4).
\chapter*{Appendix}
\addcontentsline{toc}{chapter}{Appendix}
\markboth{Appendix}{\rightmark}
\begin{appendix}
\setcounter{chapter}{1}
\setcounter{section}{0}
\setcounter{equation}{0}
\setcounter{lemma}{0}
\section{Coordinate estimates}\notiz{achtung notation: remark: $\w A\neq Df\c A$ usw}
In this chapter we want to derive coordinate estimates depending on the bounds of the 
geometry, which are used many times in this work. For this, assume that for a given
Riemannian metric $g,$ open sets $U\subset\R^d$ and $V\subset\R^n,$ and a time interval
$0\in I=[t_1,t_2]\subset\R$ we have smooth
coordinate representations
\notiz{definiere die koordinatenableitungen : quer p usw. }
\begin{eqnarray*}
  \begin{array}{rcl}
     \tilde f&:&(U,\widetilde G(t))\times I\to
     (V,g)
  \end{array}
\end{eqnarray*}
and
\begin{eqnarray*}
\begin{array}{rll}
     \bar f&:&(U,\barG(t))\times I\to
     (V,\d_{eucl})
  \end{array}
\end{eqnarray*}
such that for any $t\in I~~\tilde{\bar f}_t:=\tilde{\bar f}(\;\c\;,t)$ is 
an isometric immersion. Here, $\tilde f\equiv\bar f\equiv:f$ as maps between
manifolds. In addition, we want to assume the following:
\begin{enumerate}
\item [i)] $c^{-1}\delta\le g\le c\delta$~~~~~~~and~~~~~~~
   $c^{-1}\d\le\w G_{t=0}\le c\d$\hfill
   (\stepcounter{equation}\setcounter{app1}{\arabic{equation}}\theequation)
\item [ii)] $\bet{\prn\gamma \Gamma_{(g)}}\le\gm B$ for $\bet\g\le B$\hfill
   (\stepcounter{equation}\setcounter{app2}{\arabic{equation}}\theequation)
\item [iii)] $\sum_{l=0}^{L}
  \norm{\widetilde\na^l\!\widetilde A}{}
       {L^\infty(\widetilde G(t),U)}{}\le A_L$
     \hfill
     (\stepcounter{equation}\setcounter{app3}{\arabic{equation}}\theequation)
\item [iv)] $\bet{\partial_t\widetilde G}_{\w G(t)}\le P~\forall t\in(t_1,t_2)$\hfill
     (\stepcounter{equation}\setcounter{app4}{\arabic{equation}}\theequation)
\item [v)] $\bet{\partial^\gamma\widetilde\Gamma}\le\wgm N$ for 
   $\bet\g\le N$ and $\forall t\in I.$\hfill
     (\stepcounter{equation}\setcounter{app5}{\arabic{equation}}\theequation)
\end{enumerate}
Note that for negative indices in $\Gamma_B,A_L$ and $\wgm{N}$ we have the empty
condition.\\\\
From (\Alph{chapter}.\arabic{app1}) it follows that $\bet{g^{-1}}\le c$
using polarization, and further from this, using (\Alph{chapter}.\arabic{app2}),
\begin{eqnarray}
 \bet{\p^\g g},\bet{\p^\g(g^{-1})}\le c(m,n,\gm m)\mbox{ for }0
     \le\bet\g\le m+1,\label{gbound}
\end{eqnarray}
where the second estimate follows from the first by differentiating Kronecker's
symbol\\
$\d_{\b}^\a=g^{\a\g}g_{\g\b}.$
Further, since $\bet{\pt\w G}_{\w G(t)}\le P,$
we have that for $c=c(P,\bet{I})$
\begin{equation}\label{tildeab}
  c^{-1}\d\le c^{-1}\w G_{t=0}\le\w G\le c\,\w G_{t=0}\le c\d
\end{equation}
from (\Alph{chapter}.\arabic{app4}) and  Lemma 14.2 of \cite{HAM82}.
Thus by definition of $\barG$ we also have
\begin{eqnarray}
  c^{-1}\d\le\barG\le c\d.\label{barab}
\end{eqnarray}
As in (\ref{gbound}) we note that
\begin{eqnarray}
 \bet{\p^\g\w G},\bet{\p^\g(\w G^{-1})}\le c(m,d,\wgm m)\mbox{ for }
     1\le\bet\g\le m+1.\label{barGbound}
\end{eqnarray}
\begin{lemma}\label{ortsabl}
  Under the above assumptions, the following quantities are uniformly bounded on
  $U\times I$ for any $s\in\N_0$ and any $k+r=m\in\N_0:$
\begin{itemize}
 \item [a)] $\bet{\p^k\widetilde\na^r\!\widetilde A}
   \le c(\wgm{k-1},A_{k+r},\gm{k-1})$
 \item [b)] $\bet{\p^{m+2}f}\le c(\wgm m,A_m,\gm m)$~~~~~$(m\geq-1)$
 \item [c)] $\bet{\p^{m+1}\barG},~\bet{\p^{m+1}(\barG~\!^{\!-1})}
            \le c(\wgm m,A_m,\gm m)$~~~~~$(m\geq-1)$
 \item [d)] $\bet{\p^m\barGamma}\le c(\wgm m,A_m,\gm m)$
 \item [e)] 
            $\bet{\barD^m\!\barA}_{\bar f}+\bet{\barnabla^m\!\barA}_{\bar f}
            \le c(\wgm m,A_m,\gm m)$
 \item [f)] $\bet{\p^{m+2}(\p^s g)\cf}
            \le c(\wgm m,A_m,\gm{m+s+1})$~~~~~$(m\geq-2)$
 \item [g)] $\bet{\p^{m+2}(\p^s \Gamma)\cf}
            \le c(\wgm m,A_m,\gm {m+s+2})$~~~~~$(m\geq-2)$
 \item [h)] $\bet{\p^{m+2}(D^{s\!} R)\cf}
            \le c(\wgm m,A_m,\gm{m+1},R_{m+s+2})$~~~~~$(m\geq-2),$
\end{itemize}
\notiz{f\"ur r einen anderen buchstaben}
 where $\sum_{i\le m}\norm{D^iR }{}{g}{}\le R_m$ and $c$ may also depend on 
 $d,n,m,s,P,\bet{I}$ and the constants in (\Alph{chapter}.\arabic{app1}).\\\\ Here,
$\barGamma$ are the Christoffel's symbols induced by $\barG,$ $\barD$ 
 is the respective connection induced by $\bar f,$ and $\barA$ is the second fundamental
form of $\bar f.$
\notiz{1) checke, ob man die absch\"atzung f\"ur $\Gamma$ und R braucht!!!!}
\notiz{remark: wozu und warum braucht man das alles ... auch bei durchdiffernezieren}
\end{lemma} 
\bb To show a) - g) we use induction over $L\in\N_0,$ where $m\le L.$ We let the
constants $c$ also depend on $d,n,L,s,P,\bet{I}$ and the constants in
(\Alph{chapter}.\arabic{app1}).\\ 
Ad a) \& b):
\begin{eqnarray*}
  \sum\limits_{i\a}\bet{\p_if^\a}^2=\norm{\p f}{2}{\d}{}
   \le c\norm{Df}{2}{\w f}{}\le c
\end{eqnarray*}
and analogously
\begin{eqnarray*}
 \bet{\w A}^2\stackrel{}{\le}c\norm{\w A}{2}{\w f}{}\le c(A_0)
\end{eqnarray*}
using the equivalences (\Alph{chapter}.\arabic{app1}) and (\ref{tildeab}).
With Gau\ss-Weingarten, i.e. $\w A=\p^2\!f+\p f\c\p f\c\Gamma\cf+$ $\p f\c\w\Gamma,$
we estimate
$$\bet{\p^2\!f}\le \bet{\w A}+c\big(\norm{\p f}{2}{\d}{}\norm{\Gamma}{}{\d}{}
  +\norm{\p f}{}{\d}{}\norm{\w \Gamma}{}{\d}{}\big)
    \le c(\wgm 0,A_0,\gm 0).$$
Ad c): By polarization, we get from (\ref{barab})
$\bet{\barG},\bet{\barG^{-1}}\le c.$ Since we
may write $\barG=\p f\c\p f$ we can now estimate
\begin{eqnarray*}
 \bet{\p\barG}=\bet{\p f\c\p^2\!f}\le c\norm{\p f}{}{\d}{}\norm{\p^2\!f}{}{\d}{}
    \le c(\wgm 0,A_0,\gm 0)
\end{eqnarray*}
The estimate for the inverse metric again follows from differentiating Kronecker's
symbol.\\
Ad d): $\bet{\barGamma}=\bet{\barG^{-1}\c\p\barG}\le c(\wgm0,A_0,\gm0)$
from the above.\\
Ad e): Again with Gauss-Weingarten we estimate
\begin{eqnarray*}
  \bet{\bar A}\le\norm{\p^2\!f}{}{\d}{}+c\norm{\p f}{}{\d}{}
      \norm{\barGamma}{}{\d}{}\le c(\wgm0,A_0,\gm0),
\end{eqnarray*}
and thus we have the same bound for
\begin{eqnarray*}
  \norm{\bar A}{}{\bar f}{}\le c\norm{\bar A}{}{\d}{}\le c(\wgm0,A_0,\gm0).
\end{eqnarray*}
Ad f): $\bet{(\p^sg)\cf}\le c(\gm{s-1}),~\bet{\p(\p^s g)\cf}
  \le c(\gm s)\mbox{~and}
  ~\bet{\p^2(\p^s g)\cf}\le c(\wgm0,A_0,\gm{s+1})$ follow from the chain rule.\\
Ad g): As above.\\\\
Now assume that a) - g) hold true for any $s\in\N_0$ and any $m=k+r\le L\in\N_0$
and let $m=k+r=L+1.$\\
Ad a): To estimate $\bet{\p^k\w\na^{r}\!\w A}$ we use a telescope sum to get
\begin{eqnarray*}
  \p^k\w\na^{r}\!\w A -\w \na^{L+1}\!\w A
  &=&\p^k\w \na^{r}\!\w A -\w \na^{k}\w \na^{r}\!\w A\\
  &=&\sum_{j=1}^{k}\big(\p^j\w\na^{k+r-j}\!\w A -\p^{j-1}\w\na^{k+r+1-j}\big)\w A\\
  &=&\sum _{j=1}^{k}\p^{j-1}(\p-\w\na)\w\na^{k+r-j}\!\w A
\end{eqnarray*}
so that we have to compute $(\p-\w\na)S$ for $S:=\w \na^{k+r-j}\!\w A.$ But for any
$N_{\!\tilde f\,}$-valued tensor field it holds 
\begin{eqnarray}
  \w\na S&=&\w DS-P^\top\w DS\no\\
    &=&\p S+\w\Gamma\c S+\Gamma\cf\c\p f\c S
      +\w G^{-1}\c g\cf\c S\c\w A\c\p f.\label{norableitung}
\end{eqnarray}
Substituting, we therefore get for $1\le j\le k$ using the induction hypothesis
\begin{eqnarray*}
  \bet{\p^{j-1}(\p-\w\na)\w\na^{k+r-j}\!\w A}
   &=&\bet{\p^{j-1}(\w\Gamma\c S+\Gamma\cf\c\p f\c S+\w G^{-1}\c g\cf\c S\c\w A\c\p f)}\\
   &\le&c(\wgm {k-1},A_{L},\gm {k-2},\gm {k-1})\\
   &\le&c(\wgm {k-1},A_{L},\gm {k-1}).
\end{eqnarray*}
and thus in summary
\begin{eqnarray*}
  \bet{\p^k\w\na^{r}\!\w A}&\le&\bet{\w\na^{L+1}\!\w A}+c(\wgm {k-1},A_{L},\gm {k-1})\\
  &\le&c(\wgm {k-1},A_{L+1},\gm {k-1}).
\end{eqnarray*}
Ad b): Again with Gauss-Weingarten we can estimate
\begin{eqnarray*}
  \bet{\p^{L+3}f}&=&\bet{\p^{L+1}(\w A+\p f\c\p f\c\Gamma\cf+\p f\c\w\Gamma)}\\
   &\le&c(\wgm L,A_{L+1},\gm L,\gm {L+1},\wgm {L+1})\\
   &=&c(\wgm {L+1},A_{L+1},\gm {L+1}).
\end{eqnarray*}
Ad c): Differentiating we get the bounds
\begin{align*}
  \bet{\p^{L+2}\barG}&=\bet{\p^{L+2}(\p f\c\p f)}\le c(n,L,P,\wgm{L+1},A_{L+1},\gm{L+1})
\end{align*}
Analogously to the estimates for $m\le 0,$ these bounds also hold for the inverse 
metrics.\\
Ad d): Similarly
\begin{eqnarray*}
  \bet{\p^{L+1}\barGamma}&=&\bet{\p^{L+1}(\barG^{-1}\c\p\barG)}\\
  &\le&c(\wgm L,A_L,\gm L,\wgm{L+1},A_{L+1},\gm{L+1})\\
  &\le&c(\wgm{L+1},A_{L+1},\gm{L+1}).
\end{eqnarray*}
Ad e): By induction it follows that
\begin{eqnarray*}
  \barD^{L+1}\!\barA=\p^{L+1}\!\barA
  +\sum_{u=1}^{L+1}\quad\sum_{k_1+\ldots+k_u+i\le L}\hspace{-1.7em}\p^{k_1}
   \barGamma\c\ldots\c\p^{k_u}\barGamma\c\p^i\!\barA\label{barAabsch}
\end{eqnarray*}
(where we have used $\Gamma^{\R^n}\equiv0$) we estimate using Gauss-Weingarten, i.e.
$\barA=\p^2 f+\p f\c\barGamma,$
\begin{eqnarray*}
  \bet{\barD^{L+1}\!\barA}&\le&\bet{\p^{L+1}\!\barA}
     +c(\wgm L,A_L,\gm L)\max_{i}\bet{\p^i\!\barA}\\
  &\le&c(\wgm{L+1},A_{L+1},\gm{L+1}).
\end{eqnarray*}
Since $\barD^{m}\!\barA$ is tensorial, the proposed estimate for the first
summand follows by equivalence of the metrics (\ref{barab}). For the second, we show by
induction over $j\in\N_0$ that for all $i,j\in\N_0$ we have
$$\bet{\p^i\barnabla^j\!\barA}\le c(\wgm{i+j},A_{i+j},\Gamma_{i+j}).$$
Case $j=0$ follows from Gauss-Weingarten and from what has been shown above. For the
induction step we note that for a $N_{\!\bar f}$-valued tensor field we have
$$\barnabla S=\barD S-P_{\bar f}^\top \barD S=\p S+\barGamma\!\c S+\barG^{-1}
\c\barA\c S\c\p f$$
and thus
\begin{eqnarray*}
  \bet{\p^i\barnabla^{j+1}\!\barA}
  &\le&\bet{\p^i(\p\barnabla^j\!\barA+\barGamma\!\c \barnabla^j\!\barA+\barG^{-1}
    \c\barA\c \barnabla^j\!\barA\c\p f )}\\
  &\le&c(\wgm{i+j+1},A_{i+j+1},\Gamma_{i+j+1}).
\end{eqnarray*}
For $i:=0$ it follows e), again using the equivalence (\ref{barab}) and that
$\barnabla^r\!\barA$ is tensorial.\\
Ad f): Case $m=-2$ follows from (\ref{gbound}). Using induction one can show that
for $l\ge 1$
\begin{eqnarray}
  \p^{l}\big((\p^s g)\cf \big)
  =\sum_{u=1}^{l}\quad\sum_{k_1+\ldots+k_u=l}\hspace{-1.1em}(\p^{u+s}g)\cf\c\p^{k_1}f
     \c\ldots\c\p^{k_u}f\label{chainruleformula}
\end{eqnarray}
and thus
$$\bet{\p^{m+2}\big((\p^s g)\cf \big)}\le c(\gm{m+s+1},\wgm{m},A_{m}).$$
Ad g): Analogously to f), so that we have now shown a) - g).\\
Ad h): As in e), one can show by induction (over l) that for any
$l,s\in\N_0$
\begin{eqnarray}
  D^{l+s}R^{M}=\p^{l}D^sR +\sum_{u=1}^{l}\quad
   \sum_{k_1+\ldots+k_u+k\le l-1}
   \hspace{-2.2em}\p^{k_1}\Gamma\c\ldots\c\p^{k_u}\Gamma\c\p^{k}D^sR\label{riemformel}
\end{eqnarray}
and thus
\begin{eqnarray*}
  \bet{\p^{u}D^{s\!}R }&\le&\norm{D^{u+s\!}R }{}{\d}{}+c(n,u,s,\gm{u-1})
  \sum_{i=0}^{u-1}\norm{\p^iD^{s\!} R }{}{\d}{}\\
    &\le&c(n,u,s,\gm{u-1})\sum_{j=0}^{u+s}\norm{D^{j\!}R }{}{\d}{}\\
    &\stackrel{(\Alph{chapter}.\arabic{app1})}{\le}&
       c(n,u,s,\gm{u-1})\sum_{j=0}^{u+s}\norm{D^{j\!}R }{}{g}{}\\
    &\le& c(n,u,s,\gm {u-1},R_{u+s}).
\end{eqnarray*}
From this estimate we get for $m\ge -1$ 
\begin{eqnarray*}
  \bet{\p^{m+2}\big((D^{s\!} R )\cf \big)}
     \le c(\gm{m+1},R_{m+s+2},\wgm{m},A_{m}),
\end{eqnarray*}
where we used a formula analogous to (\ref{chainruleformula}) with 
$D^sR $ instead of $\p^sg,$ i.e. for $l\geq 1$
\begin{eqnarray}
  \p^{l}\big((D^{s\!}R)\cf \big)
  =\sum_{u=1}^{l}\quad\sum_{k_1+\ldots+k_u=l}
    \hspace{-1.1em}(\p^{u}D^{s\!}R)\cf\c\p^{k_1}\!f
     \c\ldots\c\p^{k_u}\!f.
\end{eqnarray}
\eb
\begin{lemma}\label{zeitableitungen}
Assume that (\Alph{chapter}.\arabic{app1}) to (\Alph{chapter}.\arabic{app5}) holds and
assume further that we have a local flow 
$$f:(U,\widetilde G(t))\times I\to(V,g)
  \mbox{~~~~~}(U\subset\R^d,V\subset\R^n\mbox{~open})$$
of isometric immersions $f_t$ such that 
\begin{eqnarray}
  \p_t f=-\bf W,\label{pde}
\end{eqnarray}
 where $\bf W$ is a universal
linear combination of elements in
$$\{\p^k\na^{r}\!A,\w G,\w G^{-1},\p f, g\cf,(D^{s\!}R)\cf \}$$
for $k+r+2,s\le m\in\N.$ Then for any $p\in\N_0,~l\in\N$
\begin{eqnarray}
  \bet{\p^{p+1}f}\le c(\wgm{p-1},A_{p-1},\Gamma_{p-1})\label{vonoben}
\end{eqnarray}
and 
\begin{equation}
  \bet{\p^p\p_t^lf}
    \le c(\wgm{p+(l-1)m+k-1},A_{p+lm-2},
      \Gamma_{p+(l-1)m+k-1},R_{p+(l-1)m+s}).\label{rhside}
\end{equation}
Here, the constant c may also depend on $n,d,m,p,l,P,\bet{I}$ and the constants in
(\Alph{chapter}.\arabic{app1}).
\end{lemma}\notiz{s kleiner m braucht man wohl nicht... -entfernen!}
\bb We let the constant c may also depend on $n,d,m,p,l,P,\bet{I}$ and the constants in
(\Alph{chapter}.\arabic{app1}).
Assume that $\bet{\p^i\p_t^jf}\le f^{(l)}(p)<\infty$ $(0\le i\le p,~1\le j\le l)$
for an auxiliary function $f^{(l)}(p)$ to be determined. 
By induction over $l$ one can show that 
\begin{eqnarray*}
  \bet{\p^p\p_t^l\big((\p^sg)\cf \big)}&\le& c(\wgm{p-2},A_{p-2},\Gamma_{p+l+s-1}
    ,f^{(l)}(p))\\
  \bet{\p^p\p_t^l\big((\p^s\Gamma)\cf \big)}&\le& c(\wgm{p-2},A_{p-2},\Gamma_{p+l+s}
    ,f^{(l)}(p))\\
  \bet{\p^p\p_t^l\big((D^{s\!}R)\cf \big)}&\le& c(\wgm{p-2},A_{p-2},\Gamma_{p+l-1}
    ,f^{(l)}(p),R_{p+s+l})\\
  \bet{\p^p\p_t^l\w G},~\bet{\p^p\p_t^l(\w G^{-1})}
    &\le& c(\wgm{p-1},A_{p-1},\Gamma_{p+l-1},f^{(l)}(p+1))\\
  \bet{\p^p\p_t^l\w\Gamma}
    &\le& c(\wgm{p},A_{p},\Gamma_{p+l},f^{(l)}(p+2)).
\end{eqnarray*}
For $l=1,$ the above inequalities follow after differentiation in time using
\cref{ortsabl}, and then using the induction hypothesis. Taking the above inequalities
into account and now using induction over $r\in\N_0,$ one can show that
\begin{eqnarray*}
    \bet{\p^{p+k}\p_t^l\na^r\!A}
    \le c(\wgm{p+k+r},A_{p+k+r},\Gamma_{p+k+r+l},f^{(l)}(p+k+r+2)),
\end{eqnarray*}
where for the induction step we used (\ref{norableitung}) with $S:=\na^{r-1}\!A.$
We now want to estimate $\bet{\p^p\p_t^lf}.$
Due to the special structure of $\bf W\rm$ and estimating $k+r\le m-2,$ the above
inequalities reduce to
\begin{eqnarray}
  \bet{\p^p\p_t^l(g\cf)}&\le& c(\wgm{p-2},A_{p-2},\Gamma_{p+l-1}
    ,f^{(l)}(p))\label{une}\\
  \bet{\p^p\p_t^l\big((D^{s\!}R)\cf \big)}&\le& c(\wgm{p-2},A_{p-2},\Gamma_{p+l-1}
    ,f^{(l)}(p),R_{p+s+l})\label{deux}\\
  \bet{\p^p\p_t^l\w G},~\bet{\p^p\p_t^l(\w G^{-1})}
    &\le& c(\wgm{p-1},A_{p-1},\Gamma_{p+l-1},f^{(l)}(p+1))\label{trois}\\
    \bet{\p^{p+k}\p_t^l\na^r\!A}
    &\le&c(\wgm{p+m-2},A_{p+m-2},\Gamma_{p+m-2+l},f^{(l)}(p+m))\label{quattre}
\end{eqnarray}
and from \cref{ortsabl}
\begin{eqnarray}
  \bet{\p^{p+k}\na^r\!A}
    &\le&c(\wgm{p+k-1},A_{p+m-2},\Gamma_{p+k-1}),
\end{eqnarray}
where (\ref{une}) to (\ref{trois}) also hold for $l=0$ by \cref{ortsabl}, if we set
$f^{(0)}(p):=0.$ We let 
$$f^{(l)}(p):=c(\wgm{p+(l-1)m+k-1},A_{p+lm-2},\Gamma_{p+(l-1)m+k-1},R_{p+(l-1)m+s}).$$
Clearly, 
$f^{(l)}(p)$ is increasing in $l$ and $p,$ and 
$f^{(l)}(p+m)\le f^{(l+1)}(p)$ in the obvious sense.\\ 
(\ref{rhside}) is proved if we can show that
$$\bet{\p^p\p_t^lf}\le f^{(l)}(p)$$
for $1\le l\le L\in\N.$
Using \cref{ortsabl} and the remarks above it, it is easy to see that for $L=1$
\begin{eqnarray}
  \bet{\p ^p\p_tf}\le c(\wgm{p+k-1},A_{p+m-2},\Gamma_{p+k-1},R_{p+s})
    =f^{(1)}(p).\label{leins}
\end{eqnarray}
To estimate $\bet{\p^p\p_t^{L+1}f},$ we use the equation $\p_tf=\bf W\rm.$ Let $\mu\le
p$ and $\lambda\le L.$ Using
(\ref{une}) to (\ref{quattre}) for $l:=\lambda\le L,$ we get for the induction step
\begin{eqnarray*}
  \bet{\p^\mu\p_t^\lambda\w G},~\bet{\p^\mu\p_t^\lambda(\w G^{-1})}
    &\le& c(\wgm{ \mu-1},A_{ \mu-1},\Gamma_{ \mu+\lambda-1},f^{(\lambda)}( \mu))\\
  &\le& f^{(\lambda+1)}( \mu)\le f^{(L+1)}(p),
\end{eqnarray*}
\begin{eqnarray*}
  \bet{\p^{\mu+1}\p_t^{\lambda}f}&\le& f^{(\lambda)}(\mu+1)\le f^{(\lambda)}(p+m)
    \le f^{(L+1)}(p)\mbox{  for }(\lambda\geq 1),\mbox{ and}\\
  \bet{\p^{\mu+1\!}f}&\le& c(\wgm{\mu-1},A_{\mu-1},\Gamma_{\mu-1})\le f^{(L+1)}(p),
\end{eqnarray*}
\begin{eqnarray*}
  \bet{\p^\mu\p_t^{\lambda}(g\cf)}
    &\le& c(\wgm{\mu-2},A_{\mu-2},\Gamma_{\mu+\lambda-1},f^{(\lambda)}(\mu))\\
  &\le& f^{(\lambda+1)}(\mu)\le f^{(1)}(\mu) \le f^{(L+1)}(p),
\end{eqnarray*}
\begin{eqnarray*}
  \bet{\p^\mu\p_t^{\lambda}\big((D^{s\!}R)\cf \big)}
  &\le& c(\wgm{\mu-2},A_{\mu-2},\Gamma_{\mu+\lambda-1},f^{(\lambda)}(\mu),
    R_{\mu+s+\lambda})\\
  &\le& f^{(\lambda+1)}(\mu)\le f^{(L+1)}(p),
\end{eqnarray*}
\begin{eqnarray*}
  \bet{\p^{\mu+k}\p_t^{\lambda}\na^{r}\!A}
    &\le& c(\wgm{\mu+m-2},A_{\mu+m-2},\Gamma_{\mu+m-2+\lambda},f^{(\lambda)}(\mu+m))\\
  &\le& f^{(\lambda+1)}(\mu)\le f^{(L+1)}(p)\mbox{  for }(\lambda\geq 1),
  \mbox{ and finally}
\end{eqnarray*}
\begin{eqnarray*}
  \bet{\p^{\mu+k}\na^{r}\!A}
    &\le& c(\wgm{\mu+k-1},A_{\mu+m-2},\Gamma_{\mu+k-1})\\
  &\le& f^{(1)}(\mu)\le f^{(L+1)}(p).
\end{eqnarray*}
Therefore we have shown $\bet{\p^p\p_t^{L+1}f}\le f^{(L+1)}(p).$ 
Actually, if $f^{(1)}(p)$ is given, $f^{(l)}(p)$ is the
optimal function having the property $f^{(l)}(p+m)\le f^{(l+1)}(p)$, because equality
holds.
In (\ref{vonoben}) we merely stated again \cref{ortsabl}.\eb
\notiz{warum a.3?: weil wir keine lokale schranke an $\tilde\Gamma$ brauchen wie in a.1}
\begin{lemma}\label{relativecurv}
  Let $U$ be a manifold without boundary, and let for $V\subset\R^n$ open 
  $$\w f:(U,\tildeg)\to(V,g)$$
  and $$\bar f:(U,\bar g)\to(V,\d_{eucl})$$ be isometric immersions such that
  $\w f=\bar f$ as maps between manifolds. Let further $L\in\N_0$ and assume that (with
respect to the standard coordinates of $\R^n$) we have the bounds
  $$c_g^{-1}\d\le g\le c_g\d,\quad\sum_{i=0}^{L}\sup_{V}\bet{\p^i\Gamma}\le\Gamma_L
  ,\quad\sum_{i=0}^{L}\norm{\w\na^i\!\w A}{}{L^\infty(U,\tildeg)}{}\le  A_L.$$
  Then 
  \begin{eqnarray}
    \sum_{i=0}^{L}\big(\norm{\barD^{i}\!\barA}{}{L^\infty(U,\bar g)}{}
      +\norm{\barnabla^{i}\!\barA}{}{L^\infty(U,\bar g)}{}\big)\le c,\label{barnaabsch2}
  \end{eqnarray}
  where c only depends on $c_g,\Gamma_L, A_L$ and the dimension n.
\end{lemma}
\bb We proceed as follows: Since the extrinsic and relative curvature of $\w f$
is bounded we also get a bound for the intrinsic curvature $\w R$ of $(U,\tildeg)$. Now
(\ref{barnaabsch2}) can be estimated by a uniform pointwise estimate. Namely,
we can control $\p^i\tildeg$ at the origin of Riemannian normal coordinates in terms
of the intrinsic curvature and its derivatives. From this we get uniform bounds for
$\p^i\!f^\a$ and thus also for the relative curvature of $\bar f.$\\\\
To begin, since $R=\p\Gamma+\Gamma\c\Gamma$ we get, using (\ref{riemformel}) for $s=0$,
that $\bet{D^{l\!}R}\le c(\Gamma_{l+1})$ and thus
$$\norm{D^{l\!}R}{}{L^\infty(V,g)}{}\le c(c_g,\Gamma_{\!L})\qquad(l\le L-1).$$
From the equations of Gauss (\ref{135}) we know 
$$\w R=R\cf\star Df\star Df\star Df\star Df+\w A\s\w A=\Q 00{}+P^0_2(\w A).$$
Differentiating covariantly we get $\w\na^{m\!}\w R=\Q {m-1}1{}+P^{m}_2(\w A)$ for
$m\geq1,$ and therefore 
$$\sum_{i=0}^{L-1}\norm{\w\na^{i\!}\w R}{}{L^\infty(U,\tildeg)}{}
  =:\w{\mathcal R}_{L-1}
\le c(n,c_g,\Gamma_L, A_{L-1}).$$
For arbitrary but fixed $x_0\in U$ we choose Riemannian
normal coordinates with respect to $\tildeg$ centred at $x_0$ with 
$2^{-1}\d\le(\tildeg\!_{ij})\le 2\d$ locally as bilinear forms. 
It is not hard to show that $\w\na^{m\!}\w R(x_0)$ in the above coordinates is a
universal linear combination of elements in the set $\{\p^2\tildeg(x_0),\ldots,
\p^{m+2}\tildeg(x_0)\}.$ The other direction is also true, as it has been proven in 
\cite{GJW}: Using a
system of symmetries for the
$m$-th partial derivative of $\tildeg$ obtained from differentiating the Lemma of
Gauss in Riemannian normal coordinates $\{x^i \}_{i=1,2}$, i.e.\\
$\d_{ik}x^k=\tildeg\!_{ik}x^k,~m+2$ times, it is
possible to isolate $\p^{m+2}\tildeg$ so that one arrives at
\begin{eqnarray*}
  \sum_{\sigma\in S_m}(\w\na^{m\!}\w R)_{k_{\sigma(m+2)}\cdots k_{\sigma(3)}r
  k_{\sigma(2)}sk_{\sigma(1)}}\strichklein{$x_0$}
  =C(n)\Big(\p^{m+2}_{k_{m+2}\cdots k_1}\tildeg\!_{rs}
  +\substack{\mbox{a polynomial in lower}\\\mbox{order partials of }}~\tildeg\Big).
\end{eqnarray*}
Using induction it finally follows that
\begin{eqnarray*}
  \p^{m+2}_{k_{m+2}\cdots k_1}\tildeg\!_{rs}\strichklein{$x_0$}
  =C(n)\Big(\sum_{\sigma\in S_m}(\w\na^{m\!}\w R)_{k_{\sigma(m+2)}\cdots k_{\sigma(3)}r
  k_{\sigma(2)}sk_{\sigma(1)}}
  +\substack{\mbox{a polynomial in lower order}\\\mbox{covariant derivatives of}
  }~\w R\Big).
\end{eqnarray*}
The last equation can also be obtained when differentiating (5.1) in \cite{GJW}.
Therefore, in $x_0$ we have $$\w\Gamma=0\quad\mbox{and}\quad
\sum_{i=1}^{L}\bet{\p^i\w\Gamma}\le c(n,\w{\mathcal R}_{L-1})\le c(n,c_g,\Gamma_L,
A_{L-1}).$$
Now it is easy to see that for a) - h) in \cref{ortsabl} we more generally have
pointwise estimates, i.e.
a) $\bet{\p^k\widetilde\na^r\!\widetilde A}(x)\le c(\wgm{k-1}(x),A_{k+r},\gm{k-1})
,\ldots,\mbox{g) }\bet{\p^{m+2}(D^{s\!} R)\cf}(x)\le c(\wgm
m(x),A_m,\gm{m+1},R_{m+s+2}),$ provided
$\bet{\partial^i\widetilde\Gamma}(x)\le\wgm N(x)$ for $i\le N$ and $\forall t\in I.$
Thus from e) in this Lemma with $t_1=t_2=0$ 
and by naturality of the second
fundamental form $\barA$, $\barnabla$ and $\barD,$ we infer that for $i\le L$
\begin{eqnarray*}
  \bet{\barD^{i}\!\barA}_{\bar f}(x_0)+\bet{\barnabla^{i}\!\barA}_{\bar f}(x_0)
  &\le& c(n,c_g,\wgm L(x), A_L,\Gamma_L)\\
  &\le& c(n,c_g,\Gamma_L, A_L),
\end{eqnarray*}
which proves the Lemma.\eb
   
\section[Bounds on the metric in special coordinates]
{Bounds on the metric in harmonic coordinates and\\ Riemannian normal
coordinates}\notiz{unteres lemma macht [HE99] scale-invariant. $->$ check!}
\begin{lemma}\rm (Harmonic Coordinates; cf. \cite{HE99},
  Theorem 1.2)\it\label{harmcoord1}. Let
  $(M,g)$ be a smooth, complete Riemannian n-manifold. Assume that for some 
  $k\in\N_0$ and a given universal constant $c(n,k),~R>0$ can be chosen such that
  $$R\Lambda(k)\le c$$ where $\Lambda_{(M,g)}(k):=\sum_{i=0}^{k}
    \norm{D^i\ric_{(M,g)}}{\frac{1}{i+2}}{L^\infty(M,g)}{}+\inj(M,g)^{-1}.$
  Then for any $p\in M$ there exists a g-harmonic chart $$\varphi:B_R^g(p)\to\R^n$$ and 
  a universal $C=C(n,k)$ such that in these coordinates $g_{\a\b}(p)=\d_{\a\b}$ and
  \begin{equation*}
    \begin{array}{rlll}
      i)&1-C\Lambda R)\delta\le (g_{\a\b})
        \le(1+C\Lambda R)\delta_{\a\b}&&\mbox{on $B_R ^g(p)$},\\
      ii)&\sup_{B_R^g(p)}\bet{\p^\g g_{\a\b}}\le C\Lambda^{\bet{\g}}&&\mbox{for all } 
        1\le \bet{\g}\le k+1,\\
      iii)&\sup_{B_R^g(p)}\bet{\p^\g\Gamma_{\b\d}^\a }\le C\Lambda^{\bet{\g}+1}&&
       \mbox{for all } 0\le \bet{\g}\le k.
    \end{array}
  \end{equation*}\eb
\end{lemma}
\begin{lemma}\rm (Riemannian normal coordinates; cf. \cite{HE99},
    Theorem 1.3)\label{riemcoord1}.\it~Let
  $(M,g)$ be a smooth, complete Riemannian n-manifold. Assume that for a given 
  universal constant
  $c(n),~R>0$ can be chosen such that $$R\Lambda\le c$$ where 
  $\Lambda_{(M,g)}:=\norm{R}{1/2}{L^\infty(M,g)}{}+\norm{DR}{1/3}{L^\infty(M,g)}{}
    +\inj(M,g)^{-1}.$
  Then for any $p\in M$ there exist Riemannian normal coordinates with respect to g 
  (of radius R around p) and a universal constant $C=C(n)$ such that in these
  coordinates $g_{\a\b}(p)=\d_{\a\b}$ and
  \begin{equation*}
    \begin{array}{rlll}
      i)&(1-C r^2\Lambda^2)\delta
        \le (g_{\a\b})\le(1+C r^2\Lambda^2)\delta&&\mbox{on $B_r^g(p)$ for all } 
        0\le r<R,\\
      ii)&\sup_{B_r^g(p)}\bet{\p_\g g_{\a\b}}\le C r\Lambda^2&&
        \mbox{for all } 0\le r<R,\\
      iii)&\sup_{B_r^g(p)}\bet{\Gamma_{\b\d}^\a }\le C r\Lambda^2&&
        \mbox{for all } 0\le r<R.
    \end{array}
  \end{equation*}
\end{lemma}
\eb
\section{Construction of cutoff functions and a partition of unity}
In the next Proposition we refer to a basic comparison estimate for the Hessian of the
distance function taken from (\cite{Pet}, Corollary 2.4). In this section we denote by
$\{\partial_\a\}_{2\le\a\le n}$ the spherical coordinate 
frame of Riemannian polar coordinates. Exceptionally, we sum from
$\a=2,\ldots,n.$
\notiz{checke letzten satz!}

\begin{prop}{\rm(\cite{Pet}, Corollary 2.4).}\,\label{Pettersen}~Let
  $(M,g)$ a Riemannian manifold with bounded sectional curvature 
  $k\le$ sec$_M\le K.$ Denoting by $(S^\b_\a)$ the Hessian of $d_p:=d_g(p,\cdot)$ in
  Riemannian polar coordinates, we have (as long as $\mathfrak{sn}_K(r)>0$)
  $$\sqrt{K}\mathfrak{ct}_{K}(r)
    \le\left(S^{\beta}_{\alpha}(r,\theta)\right)_{2\le\alpha,\beta\le n}
    \le\sqrt{k}\mathfrak{ct}_{k}(r),$$
  where
  $\sqrt{k}\mathfrak{ct}_{k}(r)
    =\frac{\mathfrak{sn}^{\prime}_{k}(r)}{\mathfrak{sn}_{k}(r)}
    =\sqrt{k}\frac{\mathfrak{cs}_{k}(r)}{\mathfrak{sn}_{k}(r)}$ and
   $\mathfrak{sn}_{k}(r):=\frac{1}{\sqrt{k}}\sin\left(\sqrt{k}r\right)$ for 
   $k\in\mathbb R\backslash\{0\}$.\eb
\end{prop}
\begin{lemma}{\rm (Cutoff functions).}\label{452}~Let
  $(M^n,g)$ be a Riemannian manifold and assume that \mbox{$\varrho>0$} can be
  chosen such that $\varrho<\inj(M,g).$ Then for any $\d\in\R^+$ satisfying 
  $$0<\d<\varrho<\inj(M,g)$$ there exists for arbitrary $p\in M$
  a cutoff function $\wg:=\wg_{p,\d,\r}\in C^{\infty}(M)$ with
  \begin{equation}\label{abschneidefunktion}
    \chi_{B_\d(p)}\le\w\g\le\chi_{B_\r(p)}\mbox{~~~~~and~~~~~}
    \norm{D^j\wg}{}{\infty}{}\le \dfrac{c}{(\r-\d)^j}
  \end{equation}
  for $j=1$ and $c=c(n).$ If additionally 
  $\bet{\sec} \le\kappa^2<\infty,$ then for any $\d,\r\in\R^+$ satisfying 
  $$0<\d<\r<\min\{\inj(M,g),\frac \pi{2\kappa} \}$$
  (\ref{abschneidefunktion}) also holds for $j=2$.
\end{lemma}
\bb
Let $p\in M$ be arbitrary but fixed, $r:=\d/\r$ and assume that $\r=1$ (scale otherwise).
Also, it is enough to restrict to $r\geq 1/2,$ i.e. we have $1/2\le
r<1<\inj(M,g)=:i_0.$ Since the exponential function 
$\exp_p:B_{i_0}(0_p):=\{v\in T_pM:\bet v<i_0\}\rightarrow B_{i_0}(p)\subset M$ is a
diffeomorphism and $\dist_g(p,q)=\bet{\exp_p^{-1}(q)}$ for $q\in B_{i_0}(p),$
we know that
the function $d_p:=\dist_g(p,\cdot)$ is smooth on $B_{i_0}(p)\backslash\{p\}.$ Now let 
$h\in C^\infty(\R)$ be a cutoff function with $h\equiv 1$ on $B^{can}_{r}(0),$
spt $h\subset B_{(1+r)/2}^{can}(0),$ $\bet{h^{'}}\le c(1-r)^{-1},$
$\bet{h^{''}}\le c(1-r)^{-2}$ and define 
$\wg:=h\circ d_p.$ After extending $\wg$ by zero to the whole of $M$, we get 
$\wg\in C^\infty(M),$ and furthermore $\wg\equiv 1$ on $B_{r}(p)$ and spt $\wg\subset
B_{(1+r)/2}(p)$ by construction.
The first derivative can easily be estimated on
$B_{(1+r)/2}(p)\backslash\overline{B}_{1/4}(p)$ as follows:\\
\begin{eqnarray*}
  \bet{D\wg}&\le& c\bet{h^\prime\circ d_p}\bet{Dd_p}\\
            &\le& c\max\limits_{\R}\bet{h^\prime}\bet{grad_{ g}d_p}\\
            &\le& c(1-r)^{-1},
\end{eqnarray*}
since $d_p$ is a distance function, i.e. $\bet{grad_gd_p}\equiv 1$.\\ 
Now we assume that $0<1/2\le r<1<\min\{\inj(M,g),\frac\pi{2\kappa} \}.$ Then, since
\begin{eqnarray*}
  D^2\wg=h^{\prime\prime}\circ d_pDd_p\otimes Dd_p+h^{\prime}\circ d_pD^2d_p,
\end{eqnarray*}
we can estimate the second derivative on $B_{(1+r)/2}(p)\backslash\overline{B}_{1/4}(p)$
\begin{eqnarray}
  \bet{D^2\wg}&\le&\max\limits_{\R}\bet{h^{\prime\prime}}\bet{Dd_p}^2
    +\max\limits_{\R}\bet{h^{\prime}}\bet{D^2d_p}\\
  &\le&c(1-r)^{-2}+c(1-r)^{-1}\bet{D^2d_p}\no\\
  &\le&c(1-r)^{-2}(1+\bet{D^2d_p}).\label{454}
\end{eqnarray}
If $S:=Dgrad_{ g}d_p$ denotes the Hessian of $d_p,$ i.e. the Weingarten operator of the
distance spheres, we pointwise have for an adapted
local $ g$-orthonormal frame $\{e_\a\}$ (summing over $\a,\b$)
\begin{eqnarray*}
  \bet{D^2d_p}^2\2&=&\2(D^2_{e_\a,e_\b}d_p)^2
                =\left(D_{e_\a}(Dd_p\cdot e_\b)\right)^2
                =(D_{e_\a}( g(grad_{g}d_p,e_\b)))^2
                = g(D_{e_\a}grad_{g}d_p,e_\b)^2\\
                &=&\2\bet{Dgrad_{ g}d_p}^2
                =\bet{S}^2.
\end{eqnarray*}
To estimate $\bet{S}^2,$ we get from \cref{Pettersen}
in Riemannian polar coordinates $\{\p_\a \}_{1\le\a\le n}$
$$\kappa\cot\left(\kappa r\right)
   \le\left(S^{\beta}_{\alpha}(r,\theta)\right)_{2\le\alpha,\beta\le n}
   \le \kappa\coth\left(\kappa r\right),$$
since $\sin(iz)=i\sinh(z)$ for $z\in\mathbb C,$ and $\bet{\sec}\le\kappa^2.$ Note, that
$S_\a^\b=0$ if $\a\c\b=0.$ This follows when differentiating the identity 
$g\cf(grad_g d_p,grad_g d_p)\equiv 1$ and taking into account that $S$ is
self-adjoint
with respect to $g.$ Because 
$0\le \kappa\cot(\kappa r)$ for $0\le r \le\frac{\pi}{2\kappa},$ and
$\kappa\coth(\kappa r)\le 5/r\le 10$ for $\kappa/2\le \kappa r\le 2$ we can
estimate
\begin{eqnarray*}
  \bet{S}^2(r,\theta)&=& g\left(S(\partial_\a),S(\partial_\b)\right)g^{\a\b} \\
                     &=& g\left(S^2(\partial_\a),\p_\b\right)g^{\a\b} \\
                     &=&trace(S^2) \\
                     &\le&(n-1)\max\{\lambda^2:~\lambda\mbox{ is an
                     eigenvalue~of~S}\}\\
                     &\le& c,
\end{eqnarray*}
i.e. $\bet{D^2\wg}\le c(1-r)^{-2}$ from (\ref{454}). Finally, rescaling yields the claim 
of the Lemma.\eb\\
To construct an appropriate partition of unity (see \cref{456}), we use Lemma
1.1 from \cite{HE99}:\notiz{i) besser formuliernen}
\begin{lemma}{\rm (Covering Lemma I)}.\label{hebeylemma}~Let
  $(M,g)$ be a smooth, complete Riemannian n-manifold with 
  ricci$_{(M,g)}\geq (n-1)kg,$ and let $\sigma>0$ be given. There exists a sequence
  $(p_i)_{i\in\mathbb N}$ of points in $M$ such that for any $r\geq\sigma:$
  \begin{itemize}
    \item[i)] the family $\left(B_r(p_i)\right)_{i\in\mathbb N}$ is a uniformly locally
               finite covering of $M$, and there is an
               upper bound $N=N(n,\sigma,r,k)$ for the number of Balls intersecting a 
               previously given one
    \item[ii)] for any $i\neq j,~B_{\sigma/2}(p_i)\cap B_{\sigma/2}(p_j)=\emptyset$
  \end{itemize}
  where, for $p\in M$ and $r>0,~B_r(p)$ denotes the geodesic ball of centre p and
  radius r.\eb
\end{lemma}
\begin{lemma}{\rm(Existence of an appropriate partition of unity).}\label{456}~Let
  $(M,g)$ be a Riemannian manifold  and assume that $\r>0$ can be chosen with 
  $$\r\max\{c_1\norm{\mbox{\rm ricci}{}_{(M,g)}}{}{\infty}{}~\hspace{-1.2em}^{1/2}
  ,\inj(M,g)^{-1}\}< 1,$$ where $c_1>0$ is a constant. 
  Then there exists a uniformly locally finite covering $\{B_\r(p_i) \}_{i\in\N}$ of
  $M$, i.e. an
  arbitrary Ball $B_\r(p_i)$ intersects with at most $N=N(n,c_1)$ other Balls of this
  covering. Furthermore, there is a smooth partition of unity 
  $\{\widetilde\eta_i \}_{i\in\N}$
  subordinate to this covering with $0\le\widetilde\eta_j\le 1$ and 
  $\bet{D\widetilde\eta_i}\le\frac{c(n,N)}\r.$\eb
\end{lemma}
\bf Remark:\rm~ A bound $\rm ricci_{(M,g)}\it\geq (n-1)k g$ for $k\in\R$ instead of 
demanding a bound for $\norm{\rm ricci_{(M,g)}\it}{}{\infty}{}$ would
suffice.\er\\
\bb We may assume that $\r=1,$ i.e. $\inj(M,g)<1$ and 
$\norm{\rm ricci}{}{\infty}{}<c_1^{-2}.$ The bound on the Ricci curvature implies that
$\rm ricci\it\geq -c(n)c_1^{-2}g.$ Choosing $r:=\sigma:=1/2$ in \cref{hebeylemma} we
obtain a covering $\{B_{1/2}(p_i) \}_{i\in\N}$ of $M$ such that
each ball $B_1(p_{i_0})$ ($i_0\in\N$) of twice the radius intersects with at most
$N(n,1/2,1,c_1)=N(n,c_1)$ other balls $B_1(p_j)\in\{B_1(p_i)\}_{i\in\N}.$
Using \cref{452} we have a sequence of bump functions $\{\g_i \}_{i\in\N}$ with
$\chi_{B_{1/2}(p_i)}\le\g_i\le\chi_{B_1(p_i)}$ and
$\norm{D\g_i}{}{\infty}{}\le c.$ Since $\{B_{1/2}(p_i) \}_{i\in\N}$ is a covering of $M$,
we have $$1\le\sum\limits_{j=1}^\infty\g_i\le N.$$
Thus we can define
$$\widetilde\eta_i:=\dfrac{\g_i}{\sum^\infty_{j=1}\g_j}.$$
For $i\in\mathbb N$ and $q\in B_1(p_i)$ both arbitrary but fixed, let wlog (renumber 
otherwise) $\{B_1(p_1),\ldots,B_1(p_k)\}$ with $k\le N$ the only pairwise distinct
Balls containing q. Then we have in q
\begin{eqnarray*}
  \bet{D\widetilde\eta_i}&\le&\dfrac{\bet{D\g_i}}{\sum_{j=1}^k\g_j}
    +\dfrac{\g_i}{\left(\sum_{j=1}^k\g_j \right)^2}
    \Big\lvert \sum\nolimits^k_{j=1}D\g_j\Big\rvert\\
  &\le&c+\sum\nolimits^k_{j=1}\bet{D\g_j} \\
  &\le&c+cN.
\end{eqnarray*}
Since $\sum\widetilde\eta_i\equiv 1$ the lemma follows after rescaling.\eb\\
The next Lemma is in spirit of
\cref{hebeylemma} and uses a kind of Vitali's argument. It is needed in the context
of interior estimates.
\begin{lemma}{\rm (Covering Lemma II)}.\label{hilfslemma51}~Let
  $(M,g)$ be a smooth complete Riemannian manifold with $\ricci_{(M,g)}\geq kg$
  as bilinear forms for some $k\in\R.$ Then any closed geodesic ball 
  $\overline{B}_\varrho\subset M$
  with $\varrho\bet{k}^{1/2}\le1$ can be covered with at most
  $\Gamma(n)$ other balls $B_{\varrho/2}$ of radius $\varrho/2.$ If even
  $\norm{\ricci_{(M,g)}}{}{\infty}{}<\infty,$ then such a cover exists for any 
  $\varrho>0$ with
  $\varrho^2\Vert \ricci_{(M,g)}\Vert_{\infty}\le 1.$
\end{lemma}\notiz{wo werden die koordinaten gew\"ahlt, sodass man sich sp\"ater (siehe
ende des kapitels)}
\notiz{auf diese eigenschaft beruft? -$>$ check}
\bb For bounded $U\subset M$ open, we let $\bet{U}:=\mu_g(U).$ From the remark below
Theorem 1.1 of \cite{HE99}, which is based on a comparison theorem of Bishop and Gromov,
it follows that for any $\varrho>0$ and $p\in M$ we have 
\begin{eqnarray*}
  \bet{B_{2\varrho}(p)}&\le&\mbox{exp}
    \big(\sqrt{4(n-1)\bet{k}}\varrho \big)8^n\bet{B_{\varrho/4}(p)}\\
  &\le&\Gamma(n)\bet{B_{\varrho/4}(p)},
\end{eqnarray*}
when we define $\Gamma(n):=\mbox{exp}\big(\sqrt{4(n-1)}\big)8^n.$
Now for $p\in M$ and $\varrho>0$ arbitrary but fixed, consider the system of sets
$$\mathcal M:=\Big\{S:S=\bigcup\limits_{i\in\Lambda\subset\N}B_{\varrho/4}(q_i)
\mbox{~for a disjoint union of balls }B_{\varrho/4}(q_i)\subset B_\varrho(p)
\subset M \Big\}.$$ It is easy to see that $\mathcal M$ (is non-empty and) is partially
ordered by inclusion and every chain in $\mathcal M$ has an upper bound.
Hence, by Zorn's lemma, $\mathcal M$ contains a maximal element $\widehat{\mathcal S}$
in $\mathcal M.$ Now consider an arbitrary finite number N of Balls in 
$\widehat{\mathcal S}.$ After renumbering, we may assume that
$\{q_1,\ldots,q_N\}$ are the centres of such balls. Because 
\begin{eqnarray*}
  N\bet{B_{\varrho}(p)}&=&\sum\limits_{j=1}^N\bet{B_\varrho(p)}
  \le\sum\limits_{j=1}^N\bet{B_{2\varrho}(q_j)}
  \le\sum\limits_{j=1}^N\Gamma\bet{B_{\varrho/4}(q_j)}
  =\Gamma\bet{\cup_{j=1}^NB_{2\varrho}(q_j)}
  \le\Gamma\bet{B_\varrho(p)}
\end{eqnarray*}
we get, since $\bet{B_\varrho(p)}>0,$ that $N\le\Gamma,$ and thus $B_\varrho(p)$
contains at most $\Gamma$ disjoint balls of radius $\varrho/4.$ Now it is easy to see
that $$\overline{B}_\varrho(p)\subset\bigcup_{i\in\Lambda_{\widehat{\mathcal S}}}
B_{\varrho/2}(q_i)$$ and also that for the number of elements we have
$\bet{\Lambda_{\widehat{\mathcal S}}}\le\Gamma.$ This shows the existence of the desired
cover.
Now let $\norm{\ricci_{(M,g)}}{}{\infty}{}<\infty.$ Defining  
$k:=-2n\norm{\ricci_{(M,g)}}{}{\infty}{}$ we get $\ricci_{(M,g)}\geq kg,$ since
$\bet{\ricci_{(M,g)}(v,v)}\le 2n\bet{\ricci_{(M,g)}}_g\bet{v}_g^2$ for any $v\in TM.$ 
If one redefines $\Gamma:=\mbox{exp}\big(\sqrt{8(n-1)}\big)8^n,$ the second statement
follows.\eb\\
\notiz{uniformly locally finite covering comes from ounded geometry...}

\section{Interpolation inequalities}
For this section, apart from for \cref{theorem67}, we assume that $\Sigma$ is a 
d-dimensional Riemannian manifold and that $\g\in C^1_{c}(\Sigma)$ satisfies
$$0\le\g \le1,~~~~\bet{\na\g}\le G.$$ 
Except for \cref{korollar64} and \cref{theorem67} the interpolation inequalities from
\cite{KS02} carry over to the Riemannian setting, and are stated for convenience of the
reader.
%
%
\begin{lemma}\label{lemma61}
Let $\frac{1}{p}+\frac{1}{q}=\frac{1}{r},~1\le p,q,r<\infty$ and 
$\alpha+\beta=1,~\alpha,\beta\geq 0.$ For
$s\ge max\{\alpha q,\beta p\}$ and $-\frac{1}{p}\le t\le\frac{1}{q}$ we have 
\begin{eqnarray*}
\bigg(\is\bet{\na \phi}^{2r}\gamma^s\mm\bigg)^{\frac{1}{r}}&\le&
c\bigg(\iz\bet{\phi}^q\gamma^{s(1-tq)}\mm\bigg)^{\frac{1}{q}}
     \bigg(\is\bet{\na ^2\phi}^p\gamma^{s(1+tp)}\mm\bigg)^{\frac{1}{p}}\\
     &&\quad+c G s\bigg(\iz\bet{\phi}^q\gamma^{s-\alpha q}\mm\bigg)^{\frac{1}{q}}
     \bigg(\is\bet{\na \phi}^p\gamma^{s-\beta p}\mm\bigg)^{\frac{1}{p}}\\
\end{eqnarray*}
where c=c(d,r). 
\end{lemma}
\bb
Using integration by parts, we get
\begin{eqnarray*}
  \is\bet{\na \phi}^{2r}\g^s\mm
    &=&\is\eck{\phi,\na^*\ka\gamma^s\bet{\na\phi}^{2r-2}\na\phi\kz}\mm\\
  &=&-\is\eck{\phi,\gamma^s\bet{\na \phi}^{2r-2}\Delta\phi}\mm
     -s\is\eck{\phi,\gamma^{s-1}\na_{e_i}\gamma\bet{\na \phi}^{2r-2}\na_{e_i}\phi}\mm\\
  &&\quad-2(r-1)\is\eck{\phi,\na_{e_i}\phi}
     \eck{\na\phi,\na_{e_i}(\na\phi)}\bet{\na \phi}^{2r-4}\g^s\mm\\
  &\le&c(d,r)\is\bet{\phi}\gamma^s\bet{\na \phi}^{2r-2}\bet{\na ^2\phi}\mm
      +c(d)s G\is\bet{\phi}\gamma^{s-1}\bet{\na \phi}^{2r-1}\mm\\
  &\le& c\bigg(\iz\bet{\phi}^q\gamma^{s(1-tq)}\mm\bigg)^{\frac{1}{q}}
        \bigg(\is\bet{\na \phi}^{2r}\gamma^s\mm\bigg)^{\frac{r-1}{r}}
        \bigg(\is\bet{\na ^2\phi}\gamma^{s(1+tp)}\mm\bigg)^{\frac{1}{p}} \\
    &&\quad +c G s\bigg(\iz\bet{\phi}^q\gamma^{s-\alpha q}\mm\bigg)^{\frac{1}{q}}
        \bigg(\is\bet{\na \phi}^{2r}\gamma^s\mm\bigg)^{\frac{r-1}{r}}
	\bigg(\is\bet{\na \phi}^p\gamma^{s-\beta p}\mm\bigg)^{\frac{1}{p}},
\end{eqnarray*}
since $\frac{1}{q}+\frac{r-1}{r}+\frac{1}{p}=1.$
\eb

%
%
\begin{korollar}\label{korollar62}
For $2\le p<\infty$ and $s\ge p$ we have
\begin{eqnarray*}
\bigg(\is\bet{\na\phi}^p\gamma^s\mm\bigg)^{\frac{1}{p}}
   \le\e\bigg(\is\bet{\na^2\phi}^p\gamma^{s+p}\mm\bigg)^{\frac{1}{p}}
   +\frac{c}{\e}\bigg(\iz\bet{\phi}^p\gamma^{s-p}\mm\bigg)^{\frac{1}{p}},
\end{eqnarray*}
where $c=c(d,p,s,G).$\\
\end{korollar}

\bb
We take $p=q=2r,~\alpha=1,~\beta=0$ and $t=\frac{1}{s}$ in the previous lemma and obtain
\begin{eqnarray*}
    \bigg(\is\bet{\na \phi}^p\gamma^s\mm\bigg)^{\frac{1}{p}}&\le&
        \bigg(\iz\bet{\phi}^p\gamma^{s-p}\mm\bigg)^{\frac{1}{2p}}
	\bigg(\is\bet{\na ^2\phi}^p\gamma^{s+p}\mm\bigg)^{\frac{1}{2p}} \\
    &&\quad+c G s\bigg(\iz\bet{\phi}^p\gamma^{s-p}\mm\bigg)^{\frac{1}{2p}}
	\bigg(\is\bet{\na \phi}^p\gamma^s\mm\bigg)^{\frac{1}{2p}}.
\end{eqnarray*}
The claim follows using Young's inequality and absorption.
\eb
%
%
\begin{korollar}\label{korollar63}
  For $2\le p <\infty,~k\in\mathbb N,~s\ge kp$ and $c=c(d,p,s,k, G)$ we have
  \begin{eqnarray}
    \bigg(\is\np k^p\gamma^s\mm\bigg)^{\frac{1}{p}}\le
      \e\bigg(\is\np {k+1}^p\gamma^{s+p}\mm\bigg)^{\frac{1}{p}}
      +\frac{c}{\e}\bigg(\iz\bet{\phi}^p\gamma^{s-kp}\mm\bigg)^{\frac{1}{p}}.\label{62}
  \end{eqnarray}
\end{korollar}

\bb
(Induction over k). For $k=1$ the statement holds by \cref{korollar62}. Now let 
$k\ge 1.$ If one substitutes in the same Corollary $\phi $ by $\na^k\phi,$ we obtain
\begin{eqnarray*}
    \bigg(\is\np {k+1}^p\gamma^s\mm\bigg)^{\frac{1}{p}}\
       \le\e\bigg(\is\np {k+2}^p\gamma^{s+p}\mm\bigg)^{\frac{1}{p}}
       +\frac{c_k}{\e}\bigg(\is\np k^p\gamma^{s-p}\mm\bigg)^{\frac{1}{p}}.
\end{eqnarray*}
The inductive hypothesis, applied to $s-p$ instead of $s$, yields
\begin{eqnarray*}
    \bigg(\is\np k^p\gamma^{s-p}\mm\bigg)^{\frac{1}{p}}
       \le\tau\bigg(\is\np {k+1}^p\gamma^{s}\mm\bigg)^{\frac{1}{p}}
       +\frac{c}{\t}\bigg(\iz\bet{\phi}^p\gamma^{s-(k+1)p}\mm\bigg)^{\frac{1}{p}},
\end{eqnarray*}
where $c_k=c(d,p,s-p,k, G).$ Combining the inequalities proves the result.
\eb

%
%
\begin{korollar}\label{korollar64}
   For all $0\le k\le l-1$ and $s\geq 2(l-1)$ we have
      \begin{eqnarray}
         \is\bet{\na^{k}\!A}^{2}\gamma^{s}\mm
           \le \e\is\bet{\na^{l}\!A}^{2}\gamma^{s+2(l-k)}\mm
           +c(1+\e^{-1})\norm A220,\label{63}
      \end{eqnarray}
   where c=c(d,s,k,l,$ G$).
\end{korollar}
\bb
   (Induction over $l\in\mathbb N$). For $l=1$ and $k=0$ the assertion is clearly
  true. Now let $0\le k\le l$ and $s\geq 2l.$ It follows from \cref{korollar63}
  for $p:=2$ and $\phi:=A$ that
  $$\is\bet{ \na^{l}\!A}^{2}\gamma^{s}\mm
      \le\e\is\bet{ \na^{l+1}\!A}^{2}\gamma^{s+2}\mm
      +\frac{c}{\e}\norm A220.$$
    For $0\le k\le l-1,$ we get by induction hypothesis
   \begin{eqnarray*}
      \is\bet{ \na^{k}\!A}^{2}\gamma^{s}\mm
      &\le& \is\bet{ \na^{l}\!A}^{2}\gamma^{s+2(l-k)}\mm+c\norm A220\\
      &\le&\e\is\bet{\na^{l+1}\!A}^{2}\gamma^{s+2(l-k)+2}\mm
         +c(1+\e^{-1})\norm A220,
   \end{eqnarray*}
   where we used \cref{korollar63} in the last step for
   $\tilde p:=2;~\tilde k:=l;~\tilde s:=s+2(l-k)\geq2\tilde k$.\\
   \mbox{~}\eb
%
%
\begin{theorem}\label{theorem65}
For $k\in\mathbb N,~1\le i\le k$ and $s\ge 2k$ we have the inequality
\begin{eqnarray} 
   \bigg(\iz\np i^{\frac{2k}{i}}\gamma^s\mm\bigg)^{\frac{i}{2k}}
      \le c\norm {\phi}{1-\frac{i}{k}}{\infty}{0}\bigg(\bigg(\is\np
k^2\gamma^s\mm\bigg)^{\frac{1}{2}}
       +\norm{\phi}{}{2}{0}\bigg)^{\frac{i}{k}}\label{64}
\end{eqnarray}
where $c=c(d,s,k, G)$.
\end{theorem}

\bb
For $1\le i\le k$ define
$$a_i:=\bigg(\is\np
i^{\frac{2k}{i}}\gamma^s\mm\bigg)^{\frac{i}{2k}},~~a_0:=\norm{\phi}{}{\infty}{0}$$
$$b_i:=\bigg(\iz\bet{\phi}^{\frac{2k}{i}}\mm\bigg)^{\frac{i}{2k}},
   ~~b_0:=\norm{\phi}{}{\infty}{0}$$
Using \cref{lemma61} with 
$r=\frac{k}{i},~p=\frac{2k}{i+1},~q=\frac{2k}{i-1},~t=\alpha=0$ and $\beta=1$
we obtain for $s\ge 2k$
\begin{eqnarray*}
\bigg(\is\bet{\na\phi}^{\frac{2k}{i}}\gamma^s\mm\bigg)^{\frac{i}{k}}
&\le& c\bigg(\iz\bet{\phi}^{\frac{2k}{i-1}}\gamma^{s}\mm\bigg)^{\frac{i-1}{2k}}
  \bigg(\is\bet{\na^2\phi}^{\frac{2k}{i+1}}\gamma^{s}\mm\bigg)^{\frac{i+1}{2k}} \\
&&\quad+c G s\bigg(\iz\bet{\phi}^{\frac{2k}{i-1}}\gamma^{s}\mm\bigg)^{\frac{i-1}{2k}}
  \bigg(\is\bet{\na\phi}^{\frac{2k}{i+1}}
  \gamma^{s-\frac{2k}{i+1}}\mm\bigg)^{\frac{i+1}{2k}}\\
\end{eqnarray*}

With this definitions, we get from \cref{lemma61} if we substitute $\phi$ by
$\na^{i-1}\phi$
\begin{eqnarray*}
     a_i^2&\le& ca_{i-1}\bigg[a_{i+1}
     +\bigg(\is\np i^{\frac{2k}{i+1}}\gamma^{s-\frac{2k}{i+1}}
      \mm\bigg)^{\frac{i+1}{2k}}\bigg].
\end{eqnarray*}
On the the other hand, \cref{korollar63} implies for $s\geq 2k$ and $1\le i\le k-1$
\begin{eqnarray*}
\bigg(\is\np i^{\frac{2k}{i+1}}\gamma^{s-\frac{2k}{i+1}}\mm\bigg)^{\frac{i+1}{2k}}
     &\le& c\bigg(\is\np {i+1}^{\frac{2k}{i+1}}\gamma^s\mm\bigg)^{\frac{i+1}{2k}}
       +c\bigg(\iz\bet{\phi}^{\frac{2k}{i+1}}
       \gamma^{s-\frac{2k}{i+1}-i\frac{2k}{i+1}}\mm\bigg)^{\frac{i+1}{2k}}\\
     &\le& c(a_{i+1}+b_{i+1}),
\end{eqnarray*}
thus $a_i^2\le c(a_{i+1}+b_{i+1})$. We further get by an interpolation inequality for 
$L^p$-norms (see e.g. \cite{Eva} for the euclidean case), if $s\geq 2k$ and
$1\le i\le k-1$ $$(a_i+b_i)\le c(a_{i-1}+b_{i-1})(a_{i+1}+b_{i+1}).$$
Using a convexity argument (see \cite{HAM82}, Corollary 12.5), we can further estimate
\begin{eqnarray*}
a_i&\le& a_i+b_i \\
&\le& c(a_0+b_0)^{1-\frac{i}{k}}(a_k+b_k)^{\frac{i}{k}} \\
&\le& c\norm {\phi}{1-\frac{i}{k}}{\infty}{0}\bigg(\bigg(\is\np k^2\gamma^s
    \mm\bigg)^{\frac{1}{2}}
  +\norm{\phi}{}{2}{0}\bigg)^{\frac{i}{k}}.
\end{eqnarray*}
\eb

%
%
\begin{korollar}\label{korollar66}
   Let $0\le j_1,\ldots,j_r\le k,\;j_1+\ldots+j_r=2k,\;s\geq 2k$ and $r\geq 2.$ Then we 
   have
      \begin{eqnarray}
         \is\bet{\na^{j_1}\phi}\cdot\ldots\cdot\bet{\na^{j_r}\phi}\gamma^{s}\mm\label{65}
         \le c\norm\phi{r-2}\infty0\bigg(\is\bet{\na^{k}\phi}^{2}\gamma^{s}\mm
           +\norm\phi220\bigg),
      \end{eqnarray}
      where c=c(k,d,r,s,$ G$).
\end{korollar}

\bb
   Let $j_1,\ldots,j_l\geq 1$ and $j_{l+1},\ldots,j_r=0.$ Then it follows from H\"older's
   inequality and (\ref{64}) that 
\begin{eqnarray*}
   \lefteqn{\hspace{-4em}\is\bet{\na^{j_1}\phi}\cdot
     \ldots\cdot\bet{\na^{j_r}\phi}\gamma^{s}\mm}\\
     ~~~~~~~~~&\le&\norm\phi{r-l}\infty0\prod\limits^{l}_{i=1}
       \bigg(\is\bet{\na^{j_i}\phi}^{\frac{2k}{j_i}}
       \gamma^{s}\mm \bigg)^{\frac{j_i}{2k}} \\ 
   &\le&c\norm\phi{r-l}\infty0\prod\limits^{l}_{i=1}
     \bigg[\norm\phi{1-\frac{j_i}{k}}\infty0
     \bigg(\bigg(\is\bet{\na^{k}\phi}^{2}\gamma^{s}\mm\bigg)^{\frac{1}{2}}
     +\norm\phi{}20\bigg)^{\frac{j_i}{k}}\bigg] \\
   &\le& c\norm\phi{r-2}\infty0\bigg(\is\bet{\na^{k}\phi}^{2}\gamma^{s}\mm
     +\norm\phi220\bigg).
\end{eqnarray*}
\eb

\dipl{Theorem 6.7}
\begin{theorem}\label{theorem67}
  Let $(M^n,g)$ be a (possibly non-compact) Riemannian manifold with\\
  $\norm{ricci_{(M,g)}}{1/2}{\infty}{}+\inj(M,g)^{-1}=:\Lambda<\infty$
  and $f:(\Sigma^2,\tildeg)\rightarrow (M,g)$ an isometric $C^2-immersion$.
  For $u\in C^1_c(\Sigma),~2<p\le\infty,~1\le m\le\infty$
  and $0<\a< 1$ with
  $\frac{1}{\a}=(\frac{1}{2}-\frac{1}{p})m+1$ we have
  \begin{eqnarray}\label{66}
    \norm{u}{}{\infty}{}\le c\norm{u}{1-\a}{m}{}\left(\norm{\na u}{}{p}{}
    +\norm{uA}{}{p}{}+\Lambda\norm{u}{}{p}{} \right)^\a,
  \end{eqnarray}
  where $c=c(n,m,p).$
\end{theorem}

\bb We may assume that $u$ is non-negative and that
\begin{eqnarray}
  c_n\big(\norm{\na u}{}{p}{}+\norm{uA}{}{p}{}
     +\Lambda \norm{u}{}{p}{} \big)=1,\label{67}
\end{eqnarray}
where $c_n$ is the constant in the Michael-Simon Sobolev inequality (\ref{mssieqn}).
Letting $q=\frac{p}{p-1}\in[1,2)$ we infer for any $\t\geq 0$
\begin{eqnarray*}
  \norm{u^{1+\t}}{}{2}{}&\le& c_n\big(\norm{\na(u^{1+\t})}{}{1}{}
    +\norm{u^{1+\t}A}{}{1}{}+\Lambda \norm{u^{1+\t}}{}{1}{} \big) \\
  &\le& c_n\norm{u^\t}{}{q}{}\big((1+\t)\norm{\na u}{}{p}{}
    +\norm{uA}{}{p}{}+\Lambda \norm{u}{}{p}{} \big) \\
  &\le& (1+\t)\norm{u^\t}{}{q}{},
\end{eqnarray*}
where (\ref{67}) was used in the last step. With $k=\frac{2}{q}\in(1,2]$
we rewrite this as
\begin{eqnarray*}
  \norm{u}{}{k(1+\t)q}{}\le(1+\t)^{\frac{1}{1+\t}}\norm{u}{\frac{\t}{\t+1} }{\t q}{}.
\end{eqnarray*}
Putting $\t_0:=\frac{m}{q}\in(\frac{m}{2},m],~\t_{\nu+1}
:=k(1+\t_\nu),~\e_\nu:=\frac{\t_\nu}{\t_\nu+1}$ and
$c_\nu:=(1+\t_\nu)^{\frac{1}{1+\t_\nu}}$ we obtain for $\nu\in\N_0$
\begin{eqnarray}
  \norm{u}{}{\t_\nu q}{}\le c_\nu\norm{u}{\e_\nu}{\t_\nu q}{},\label{69}
\end{eqnarray} 
where
\begin{eqnarray}
  1+\t_\nu=k^\nu\t_0+\sum\limits_{\mu=0}^{\nu}k^\mu.\label{610}
\end{eqnarray}

By induction (\ref{69}) implies 
\begin{eqnarray}
  \norm{u}{}{\tau_{\nu}q}{}\le\bigg(\prod\limits_{\mu=0}^{\nu-1}c_{\mu}^{\e_{\mu+1}^{}
   \cdot\ldots\cdot\e_{\nu-1}^{}}  \bigg)
      \norm{u}{\e_{0}^{}\cdot\e_{1}^{}\cdot\ldots\cdot\e_{\nu-1}}{m}{}.\label{611}
\end{eqnarray}
Now (\ref{610}) yields
\begin{eqnarray}
 \frac{1}{c}k^{\nu}\le 1+\tau_{\nu}^{}\le ck^{\nu}\mbox{~~~   for $c=c(m,p)$,}\label{612}
\end{eqnarray}
and thus using $\e_\mu\le 1$ we can estimate
\begin{eqnarray}
  \log\prod\limits^{\nu}_{\mu=0}c_{\mu}^{\e_{\mu+1}^{}\cdot\ldots\cdot\e_{\nu}^{}}
  &\le& \sum\limits_{\mu=0}^{\nu}\frac{1}{1+\tau_{\mu}^{}}\log(1+\tau_{\mu})\no\\
  &\le&\sum\limits^{\infty}_{\mu=0}ck^{-\mu}(\log c
    +\mu\,\log k)=c(m,p)<\infty.\label{613}
\end{eqnarray}
Using $\tau_{\nu+1}^{}=k(1+\tau_{\nu}^{})$ we get from (\ref{610})
\begin{eqnarray}
  \prod\limits^{\nu}_{\mu=0}\e_{\mu}=k^{\nu}\frac{\tau_{0}^{}}{1+\tau_{\nu}^{}}
    \stackrel{\nu\rightarrow\infty}{\longrightarrow}\dfrac{\tau_{0}^{}}{\tau_{0}^{}
    +\frac{k}{k-1}}=1-\alpha.
\end{eqnarray}
Thus we may let $\nu\rightarrow\infty$ and conclude, using again (\ref{66})
\begin{eqnarray*}
  \norm{u}{}{\infty}{}&\le&c(m,p)\norm{u}{1-\alpha}{m}{}\\
   &=&c(m,p)c^{\,\alpha}_{n}\norm{u}{1-\alpha}{m}{}
   \big(\norm{\na u}{}{p}{}+\norm{uA}{}{p}{}
     +\Lambda \norm{u}{}{p}{}\big)^{\alpha}.
\end{eqnarray*}\eb
\end{appendix}

\end{document}